\documentclass{article}
\usepackage{microtype}
\usepackage[T1]{fontenc}
\usepackage{amsmath,amssymb,amsthm,mathtools,bm}
\usepackage{xcolor,booktabs,tabularx,enumitem}

\usepackage{graphicx,placeins}
\usepackage{hyperref}
\usepackage{bookmark}
\usepackage[preprint]{icml2026}
\makeatletter
\renewcommand{\Notice@String}{}
\makeatother

\definecolor{ink}{HTML}{183D55}
\hypersetup{colorlinks=true,linkcolor=ink,urlcolor=ink,citecolor=ink,
	pdfsubject={},
	pdftitle={Continuous-Time AdamW: Convergence and Nonconvergence}}
\newtheorem{theorem}{Theorem}
\newtheorem*{fenicheltheorem}{Fenichel's center-manifold theorem}
\newcommand{\restatementtitle}{}
\newtheorem*{restatedresult}{\restatementtitle}
\newenvironment{restatement}[2]{%
  \begingroup\renewcommand{\restatementtitle}{#1~\ref{#2} (restated)}%
  \begin{restatedresult}}{\end{restatedresult}\endgroup}
\newtheorem{proposition}[theorem]{Proposition}

\newtheorem{assumption}{Assumption}
\newtheorem{lemma}{Lemma}[section]

\newcounter{proofstep}

\theoremstyle{definition}
\newtheorem{definition}{Definition}
\newtheorem{example}{Example}
\theoremstyle{remark}
\newtheorem*{remark*}{Remark}
\newtheorem{assumptionremark}{Remark}[section]
\newcommand{\R}{\mathbb R}
\newcommand{\diag}{\operatorname{diag}}

\newcommand{\norm}[1]{\left\lVert#1\right\rVert}

\allowdisplaybreaks[1]
\icmltitlerunning{On the Convergence of Continuous-Time AdamW}
\begin{document}
\twocolumn[
\icmltitle{On the Convergence of Continuous-Time AdamW}
\begin{icmlauthorlist}
\icmlauthor{Bangjie Wang}{math}
\icmlauthor{Ning Ning}{stats}
\end{icmlauthorlist}
\icmlaffiliation{math}{Department of Mathematics, Texas A\&M University}
\icmlaffiliation{stats}{Department of Statistics, Texas A\&M University}
\icmlcorrespondingauthor{Ning Ning}{patning@tamu.edu}
\vskip 0.3in
]
\printAffiliationsAndNotice{}

\begin{abstract}
Continuous-time AdamW provides a framework for understanding whether an
optimizer widely used to train large-scale neural networks settles to
an equilibrium or sustains persistent motion. We study this long-time
behavior while retaining bias correction and fixed positive weight decay.
Under suitable assumptions, we establish that the dynamics are well
defined for all time and remain bounded when the decay is sufficiently
small. Our main result connects their limiting behavior to the geometry
of the objective: trajectories approach a locally invariant manifold
whose dimension equals the Hessian nullity along the associated critical
component. When critical points are isolated and have at most one such
direction, the trajectory converges to a single equilibrium. With two
degenerate directions, however, bounded persistent motion can occur even
for a smooth, strictly convex objective with a unique minimum. We also
quantify how close the long-time trajectory remains to the critical
set as weight decay decreases, without assuming convergence to a point.
These results identify objective geometry as a key to understanding when
continuous-time AdamW converges and when persistent motion is possible.

\end{abstract}

\section{Introduction}
\label{sec:introduction}

Adam \citep{KB2015} combines momentum with a coordinatewise normalization
based on exponential moving averages of squared gradients, providing a
mechanism for adapting step sizes across parameters. AdamW
\citep{LH2019} decouples weight decay from this adaptively rescaled
update. For adaptive methods, this differs from adding an $\ell_2$
penalty to the objective. AdamW is widely used to train large-scale
language, vision, and generative models
\citep{LiuDeepSeek2024,RaviSAM2024,Esser2024}.
Recent optimizer designs exploit matrix structure, as in Muon
\citep{JordanMuon2024}, and
extend it to tensors, as in TEON \citep{ZhangTEON2026} and TUON
\citep{TangTUON2026}. Although these methods can improve validation
loss in language-model pretraining, AdamW retains an advantage in
per-iteration computational cost. In the comparisons reported by
\citet[Table~4]{TangTUON2026}, AdamW has the lowest time per iteration
among AdamW, Muon, TEON, and TUON on all four tested GPT-2 and Qwen3
models. This computational advantage, together with its widespread
use, makes understanding AdamW's long-time dynamics important even as
more structured optimizers emerge.

Despite their practical success, the long-time dynamics of Adam-type methods remain substantially more complicated than those of ordinary gradient descent. The first- and second-moment variables introduce memory, coordinatewise normalization changes the geometry of the vector field, and bias correction makes the corresponding continuous-time dynamics explicitly time dependent near initialization. Moreover, convergence of Adam cannot be taken for granted: discrete Adam may fail to converge even in simple convex settings \citep{RKK2018}, and limit cycles and other nonconvergent behaviors have been identified for adaptive methods \citep{BockWeiss2019local,BockWeiss2019cycles}. Existing convergence results therefore rely on additional assumptions on the learning-rate schedule, memory parameters, gradient noise, or objective geometry
\citep{Zou2019,Defossez2020,ZhangAdam2022}.

These issues become particularly subtle for AdamW because a fixed positive weight-decay parameter changes the equilibrium geometry itself. Continuous-time AdamW provides a framework for studying how these
	mechanisms interact over long time horizons and whether they lead to an
	equilibrium or persistent motion. In the continuous-time system studied here, the position equation contains the decoupled drift $-\lambda x$. Consequently, an equilibrium position $p$ does not in general satisfy
$
\nabla F(p)=0.
$
Instead, it satisfies the coordinatewise balance
$
\nabla F(p)+\lambda D(p)p=0,
$
where $D(p)$ is determined by the second-moment field. Thus, even when the unregularized objective $F$ has a simple critical set, the long-time dynamics of AdamW need not reduce to convergence toward $\operatorname{Crit}(F)$. 

Despite AdamW's widespread use, the long-time dynamics of continuous-time
AdamW with fixed positive weight decay remain incompletely understood.
\citet{SuChenLiu2026Stability} address this gap through a local stability
analysis of autonomous continuous-time AdamW with exact gradients,
showing why convergence cannot be inferred from its equilibrium
conditions alone.
They identify an implicit objective associated with AdamW equilibria,
derive an exact one-dimensional Hurwitz stability criterion and a
sufficient higher-dimensional stability condition, and show that local
minima of the implicit objective can correspond to unstable equilibria.
Their examples also exhibit nonconvergent trajectories numerically.
These results establish that equilibrium stability cannot be taken for
granted. Our focus is the global long-time limit geometry of the
prescribed bias-corrected trajectory, including the compactness needed
to study its limit set. The second-moment assumptions differ: we allow
a state-dependent field $S$ and require its coordinates to be strictly
positive in the long-time theory, whereas their exact-gradient field
vanishes at critical points. The central question of this paper is therefore:

\begin{quote}
	\emph{For a fixed but small positive weight decay $\lambda$, what can a bias-corrected continuous-time AdamW trajectory approach as $t\to\infty$?}
\end{quote}

Under our hypotheses, the answer is governed by critical-set geometry,
particularly Hessian nullity. Local equilibrium stability does not
characterize the global limit set of a bounded trajectory. The spectrum
at a candidate equilibrium alone does not determine whether the
prescribed trajectory converges to one equilibrium or approaches a
positive-dimensional locally invariant manifold, possibly with
persistent motion. We establish global control first, then identify a
manifold containing the limit set and geometric conditions that force
point convergence. Figure~\ref{fig:adamw-long-time-geometry} illustrates a connected
critical component, the trapping region, and a schematic position
projection of the locally invariant manifold approached by the trajectory.

\paragraph{Nonautonomous continuous-time AdamW.}
Under a finite-memory scaling, the AdamW recursion in
Algorithm~\ref{alg:adamw} leads to a nonautonomous mean-moment ODE
that retains the time dependence induced by bias correction. Writing
$x(t)\in\mathbb R^d$ for the parameter vector and $m(t),v(t)\in\mathbb R^d$
for the corrected first and second moments, we study
\begin{gather*}
\dot{x}=-\frac{m}{\varepsilon+\sqrt{v}}-\lambda x,\\
\dot{m}=\frac{a}{A(t)}\bigl(\nabla F(x)-m\bigr),\\
\dot{v}=\frac{b}{B(t)}\bigl(S(x)-v\bigr),
\end{gather*}
where the operations in the position equation are coordinatewise,
$a,b>0$ are the memory rates, $\varepsilon>0$ is the stability parameter,
$\lambda>0$ is the weight-decay parameter, and $S$ is the coordinatewise
second-moment field. The bias-correction factors
$A(t)=1-e^{-at}$ and $B(t)=1-e^{-bt}$ make the system explicitly time
dependent and cause $a/A(t)$ and $b/B(t)$ to diverge as $t\downarrow0$.
Initializing the raw discrete moment filters at zero induces the
corrected initial values in~\eqref{eq:normalized-initial-value}, which
allow the solution to extend continuously to this singular initial
time; see Appendix~\ref{sec:moment-integral-derivation}. Our analysis
starts with this nonautonomous initial-value problem, establishes
global well-posedness and compact trapping, and then uses critical-set
geometry to determine when the long-time trajectory approaches a
locally invariant manifold or converges to a single equilibrium.

\subsection{Main contributions}
\label{sec:contributions}

This paper develops a global long-time geometric theory of
bias-corrected continuous-time AdamW. The main results form a sequence:
global control and localization, an invariant manifold containing the
limit set, a Hessian-nullity criterion for point convergence, and a
nullity-two construction with quantitative control of persistent motion.

\paragraph{Global control and localization.}
Theorem~\ref{thm:adamw-wp} establishes a unique global solution from the
prescribed corrected initialization despite the singular coefficients
at $t=0$. Under the long-time hypotheses, the trapping argument provides
a common compact region for sufficiently small positive decay, and
Theorem~\ref{thm:full-state-localization} localizes the long-time state
near the zero-decay equilibrium set as $\lambda\downarrow0$. These results
provide the compactness needed for long-time analysis beyond a
finite-time approximation.

\paragraph{Geometric description of the limit set.}
The conceptual centerpiece is Theorem~\ref{thm:invariant-manifold}.
Under the critical-manifold hypothesis and the long-time assumptions,
for every sufficiently small fixed $\lambda>0$ there is a locally
invariant manifold $\mathcal I_\lambda$ of the limiting autonomous
system such that
\[
\operatorname{dist}(z^\lambda(t),\mathcal I_\lambda)\longrightarrow0,
\qquad t\to\infty.
\]
For the selected critical component $\mathcal M$,
\[
\dim\mathcal I_\lambda=\dim\ker\nabla^2F(q),\qquad q\in\mathcal M.
\]
Hessian nullity thus determines the dimension of a manifold containing
the asymptotic dynamics. This is a statement about the trajectory's
limit geometry, beyond local stability of a candidate equilibrium.
Neither equality with the whole manifold nor convergence to a point
is asserted for the limit set.

\paragraph{Point convergence under identifiable geometry.}
Theorem~\ref{thm:corank-one-convergence} is the principal positive
point-convergence result. If the critical points are isolated and
$\dim\ker\nabla^2F(q)\le1$, then, under the stated long-time assumptions,
for every sufficiently small fixed positive decay the full state obeys
\[
z^\lambda(t)\longrightarrow
(p_\lambda,\nabla F(p_\lambda),S(p_\lambda)).
\]
The original bias-corrected AdamW dynamics themselves converge to a
single equilibrium. This complements the convergence guarantee for the
modified optimizer AdamMe in \citet{SuChenLiu2026Stability}; the two
results concern different dynamics and second-moment assumptions.

\paragraph{The nullity threshold and quantitative control.}
Example~\ref{ex:isolated-corank-two} shows that the geometric restriction
in the point-convergence theorem is substantive. A smooth, strictly
convex, coercive objective with a unique minimum of Hessian nullity two
admits a bounded nonconvergent corrected trajectory for every
sufficiently small positive $\lambda$. The functions $F,S$ and the
initial position are fixed independently of $\lambda$. State-dependent
second moments generate tangential motion in the degenerate plane.
Thus two degenerate directions can support persistent motion, although
nullity two alone does not force it. Under its local gradient and
distance inequalities, Theorem~\ref{thm:quantitative-localization}
quantifies how close the long-time motion remains to
$\operatorname{Crit}(F)$ as $\lambda\downarrow0$, whether or not a point
limit exists.

\begin{center}
\small
\renewcommand{\arraystretch}{1.15}
\begin{tabularx}{\linewidth}{@{}>{\raggedright\arraybackslash}X >{\raggedright\arraybackslash}X@{}}
\toprule
\textbf{Critical-point geometry} & \textbf{Long-time conclusion} \\
\midrule
Isolated critical points with Hessian nullity at most one
& Point convergence under Theorem~\ref{thm:corank-one-convergence}. \\[1mm]
A unique minimum with Hessian nullity two
& Persistent bounded motion can occur, as in
Example~\ref{ex:isolated-corank-two}. \\
\bottomrule
\end{tabularx}
\end{center}

\paragraph{Proof strategy.}
The analysis proceeds from global control to local geometry. We first resolve the singular behavior caused by bias correction at
initialization (Section~\ref{sec:main-wellposedness}) and establish compact
trapping and invariant-manifold localization (Section~\ref{sec:long-time-assumptions}).
We then study point convergence (Section~\ref{sec:main-equilibria}) and
quantitative localization (Section~\ref{sec:main-quantitative}). This separates two questions that are often conflated: whether a trajectory remains close to the critical geometry, and whether it actually converges to a single point. The distinction between these questions is central to the convergence and nonconvergence phenomena established in this paper.

\subsection{Relation to previous work}
\label{sec:related_work}

\paragraph{Discrete Adam and AdamW.}
Adam was introduced by \citet{KB2015} as an adaptive stochastic optimization method combining first-moment momentum and second-moment normalization. AdamW \citep{LH2019} decouples weight decay from this adaptive gradient step, resolving the inequivalence between $\ell_2$ regularization and weight decay for adaptive methods. The convergence theory of Adam is subtle. \citet{RKK2018} gave simple convex examples in which Adam fails to converge and proposed variants with improved long-term memory. \citet{BockWeiss2019local,BockWeiss2019cycles} studied local convergence as well as nonconvergence and limit cycles. Subsequent analyses established convergence under various assumptions on step sizes, memory parameters, and stochastic gradients
\citep{Zou2019,Defossez2020,ZhangAdam2022}.
Our focus differs from these discrete-time results: we study the global long-time dynamics of a bias-corrected continuous-time AdamW system with \emph{fixed positive} weight decay.

\paragraph{Continuous-time adaptive optimization.}
The closest continuous-time starting point is the nonautonomous Adam model of \citet{BB2021}. Their analysis establishes well-posedness and convergence properties for continuous-time Adam and connects the ODE with stochastic iterates. More general differential-equation frameworks for adaptive first-order algorithms were developed by \citet{BG2020}, while related momentum dynamics were analyzed by \citet{BBHS2021}. These works provide important foundations for continuous-time adaptive optimization, but they do not address the fixed positive decoupled position drift $-\lambda x$ considered here. That term changes both the equilibrium equation and the long-time geometry. 

\paragraph{AdamW stability and long-time limit geometry.}
\citet{SuChenLiu2026Stability} study autonomous AdamW with exact
gradients. They identify an implicit objective for its equilibria,
derive local stability criteria, and give examples of unstable
equilibria with numerical oscillations. They also introduce AdamMe,
which preserves the equilibria and admits a Lyapunov convergence
guarantee for bounded trajectories under their assumptions.

Our results concern the full bias-corrected AdamW trajectory, from
global well-posedness and compact trapping to geometric conditions
for point convergence and a bounded nonconvergent example.
The second-moment assumptions differ: our long-time theory requires
$S_i>0$, while squared exact gradients vanish at critical points.
Appendix~\ref{app:adam-flows} gives the specific stability and
convergence statements and compares the assumptions.

Continuous-time descriptions of adaptive optimization and normalization have also been studied by \citet{GT2024}, and different continuous-time scalings of Adam can produce different effective dynamics \citep{DJK2025}. The finite-memory scaling used here retains both moment variables at finite relaxation rates. This is essential for the phenomena studied in this paper: neither moment is replaced by its instantaneous target, and the resulting full-state dynamics can possess nontrivial center directions.

\paragraph{Critical manifolds and degenerate geometry.}
Nonisolated critical points naturally arise in overparameterized optimization problems. Morse--Bott and related error-bound geometries provide a useful framework for understanding such critical sets; see, for example, \citet{RB2025}. Our assumptions are deliberately more permissive than the Morse--Bott condition: along a critical component $\mathcal M$, we require constant Hessian nullity, but the Hessian kernel may contain normal directions in addition to $T_q\mathcal M$. This distinction is important because those additional neutral directions determine the dimension of the invariant manifold governing the AdamW limit set.

\paragraph{Persistent bias and quantitative localization.}
Our quantitative small-decay result is also related to work on optimization under persistent perturbations. \citet{TD2017} studied asymptotic bias in stochastic gradient search, while \citet{BLMP2025} analyzed inexact subgradient methods using geometric properties of the objective. Here the perturbation is generated internally by the coupled AdamW dynamics and the fixed decoupled decay. The resulting estimate controls the long-time distance of the trajectory to $\operatorname{Crit}(F)$ as $\lambda\downarrow0$, independently of whether the trajectory converges to a point.

\paragraph{Implicit bias versus finite equilibria.}
A separate line of work studies the implicit bias of adaptive optimizers, especially on homogeneous models, where parameter norms may diverge while normalized directions converge
\citep{WangMeng2021,XieLi2024,GronichVardi2026}.
Those results concern asymptotic directional selection and margin geometry. Our question is different: we study bounded full-state dynamics and finite equilibria under fixed positive weight decay. In particular, our invariant-manifold and nonconvergence results concern the actual position and moment variables rather than only their normalized directions.

	\section{The AdamW recursion and its continuous-time system}\label{sec:systems}
	
	We introduce the stochastic optimization setting
(Section~\ref{sec:objective-sampling}) and the constant-step AdamW recursion
(Section~\ref{sec:discrete-adamw}), derive the continuous-time system under
a finite-memory scaling (Section~\ref{sec:continuous-limit}), and identify
its autonomous limit and equilibria (Section~\ref{sec:autonomous-limit}). 
	
	\subsection{Objective and sampled gradients}\label{sec:objective-sampling}
	Let $\xi$ be a random input taking values in a measurable space $\Xi$,
	and let $f:\R^d\times\Xi\to\R$ be a measurable sample loss.
	We take the objective to be the expected sample loss
	\begin{equation}\label{eq:expected-objective}
		F(x):=\mathbb E[f(x,\xi)],\qquad x\in\R^d.
	\end{equation}
	In applications, $\xi$ may represent a data sample or a mini-batch.
	The algorithm observes independent copies $\xi_1,\xi_2,\ldots$ and uses
	sampled gradients to update the parameter $x$. Fix a deterministic
	initial position $x_0\in\R^d$.
	
	Assume that $f(\cdot,\xi)$ is continuously differentiable almost
	surely, and write 
	\[g(x,\xi):=\nabla_x f(x,\xi).\]
To justify differentiating the expectation, we use the local mean-square
regularity framework of Assumption~2.2 in~\cite{BB2021}.
Specifically, assume that the loss is integrable and the sampled gradient
is square-integrable at the initial position:
\[
\mathbb E|f(x_0,\xi)|<\infty
\quad\text{\normalfont and}\quad \mathbb E\norm{g(x_0,\xi)}^2<\infty.
\]
Assume also that the sampled gradient is locally Lipschitz in mean square:
for every compact $K\subset\R^d$, there is a constant $L_K>0$ such
that, for all $x,y\in K$,
\[
\mathbb E\norm{g(x,\xi)-g(y,\xi)}^2
\leq L_K^2\norm{x-y}^2.
\]
	The line-segment formula for $f$ and Fubini's theorem then show that
	$F$ is finite and continuously differentiable, with gradient equal to
	the mean sampled gradient. Define
	\begin{gather}
\nabla F(x)=\mathbb E[g(x,\xi)], \notag \\
S_i(x):=\mathbb E[g_i(x,\xi)^2],\quad i=1,\ldots,d. \label{eq:oracle-moments}
\end{gather}
	The vector $S$ records the coordinatewise second moments, which satisfy
	the following domination for $x\in\R^d$ and $i=1,\ldots,d$:
	\begin{equation}\label{eq:moment-domination}
S_i(x)=(\partial_i F(x))^2+\operatorname{Var}(g_i(x,\xi)) \geq (\partial_i F(x))^2.
\end{equation}
The second-moment field has a dynamical role beyond the mean gradient.
The exact-gradient choice $S=(\nabla F)^{\odot2}$ removes the variance
term, whereas a genuine second-moment field may remain positive and
vary along degenerate directions. Such variation is essential to the
rotational mechanism in Example~\ref{ex:isolated-corank-two}: it changes
the adaptive normalization along the Hessian kernel and can generate
persistent tangential motion. The example uses this freedom with
moment domination~\eqref{eq:moment-domination}; the general field $S$
therefore carries geometric information not encoded by $F$ alone.

	These sufficient conditions justify the stochastic construction. The
	results below state their hypotheses directly in terms of $F$, $\nabla
	F$, $S$, and the stochastic gradient $g$. The continuous system also
	makes sense for prescribed functions $\nabla F$ and $S$; in that case,
	we assume~\eqref{eq:moment-domination}. The stochastic results allow
	any jointly measurable stochastic gradient with the stated moments.
	
	\subsection{AdamW and bias correction}\label{sec:discrete-adamw}
	Adam maintains exponential averages of the sampled gradient and its
	coordinatewise square. We write $\widetilde m_n,\widetilde v_n$ for
	the raw filters and $m_n,v_n$ for their bias-corrected versions. Since
	the raw filters start at zero, their weights initially sum to less than
	one; using them without correction would therefore shrink the early
	moments toward zero. Dividing by the total averaging weights removes
	this initialization effect. The position $x_n$ is unchanged: only the
	moments are normalized. Throughout both the discrete and continuous
	analysis, a tilde marks a raw filter and an unadorned moment variable
	denotes its corrected version. Vector squares, square roots and
	quotients are coordinatewise, $u^{\odot2}=(u_i^2)_i$, and norms are
	Euclidean for vectors; matrix norms are the induced operator norms.

	\begin{algorithm}[t!]
		\caption{AdamW with constant step and bias correction}\label{alg:adamw}
		\setlength{\baselineskip}{1.12\baselineskip}
		\begin{algorithmic}[1]
			\REQUIRE Initial position $x_0\in\R^d$, number of steps $N\geq1$
			\REQUIRE Step $\eta>0$, memory factors $\alpha,\beta\in[0,1)$,
			offset $\varepsilon>0$, decay $\lambda\geq0$
			\STATE $\widetilde m_0\gets0$, $\widetilde v_0\gets0$
			\FOR{$n=0,\ldots,N-1$}
			\STATE Draw an independent sample $\xi_{n+1}$
			\STATE $g_{n+1}\gets\nabla_x f(x_n,\xi_{n+1})$
			\STATE $\widetilde m_{n+1}\gets\alpha \widetilde m_n+(1-\alpha)g_{n+1}$
			\STATE $\widetilde v_{n+1}\gets\beta \widetilde v_n+(1-\beta)g_{n+1}^{\odot2}$
			\STATE $m_{n+1}\gets \widetilde m_{n+1}/(1-\alpha^{n+1})$
			\STATE 
			$v_{n+1}\gets \widetilde v_{n+1}/(1-\beta^{n+1})$
			\STATE $x_{n+1}\gets x_n-\eta\bigl(m_{n+1}/(\varepsilon+
			\sqrt{v_{n+1}})+\lambda x_n\bigr)$
			\ENDFOR
			\STATE \textbf{return} $x_N$
		\end{algorithmic}
	\end{algorithm}
	Algorithm~\ref{alg:adamw} gives the constant-step version. The offset
$\varepsilon>0$ keeps the denominator positive. The decay parameter
$\lambda\geq0$ contributes a separate position increment
$-\eta\lambda x_n$; the moment filters use gradients of the sample
loss $f$~\citep{LH2019}. Setting $\lambda=0$ gives Adam, while
$\lambda>0$ gives AdamW with decoupled weight decay.
Initialize the position at $x_0$ and set
$\widetilde m_0=\widetilde v_0=0$. At each iteration
$n=0,\ldots,N-1$, draw an independent sample $\xi_{n+1}$,
evaluate its gradient at $x_n$, and update the raw moments:
\begin{equation}\label{eq:discrete-general-recursion}
 \begin{gathered}
 g_{n+1}=g(x_n,\xi_{n+1}),\\
 \widetilde m_{n+1}=\alpha\widetilde m_n+(1-\alpha)g_{n+1},\\
 \widetilde v_{n+1}=\beta\widetilde v_n+(1-\beta)g_{n+1}^{\odot2}.
 \end{gathered}
\end{equation}
Since the filters start at zero, their accumulated weights after this
update are
\begin{equation}\label{eq:discrete-correction-factors}
 A_{n+1}:=1-\alpha^{n+1}\quad\text{\normalfont and}\quad B_{n+1}:=1-\beta^{n+1}.
\end{equation}
Divide by these weights to obtain the corrected moments
\[
 m_{n+1}=\frac{\widetilde m_{n+1}}{A_{n+1}}\quad\text{\normalfont and}\quad
 v_{n+1}=\frac{\widetilde v_{n+1}}{B_{n+1}}.
\]
Finally, use the corrected moments to update the position, with
weight decay applied separately:
\[
 x_{n+1}=(1-\lambda\eta)x_n
 -\eta\frac{m_{n+1}}{\varepsilon+\sqrt{v_{n+1}}}.
\]

\subsection{Continuous-time limit}\label{sec:continuous-limit}
	To obtain finite memory rates in physical time, the memory factors
	must approach one as the step tends to zero. For each $\eta>0$,
use the same step $\eta$ and memory factors $\alpha=\alpha_\eta$,
$\beta=\beta_\eta$ at every iteration, with physical time $t_n:=n\eta$.
	Following the scaling in~\cite{BB2021}, choose $a,b>0$ and require, as $\eta\downarrow0$,
	\begingroup
\begin{equation}\label{eq:memory-scaling}
a_\eta:=\frac{1-\alpha_\eta}\eta \longrightarrow a,\qquad
b_\eta:=\frac{1-\beta_\eta}\eta \longrightarrow b.
\end{equation}
\endgroup

	Examples are $\alpha_\eta=e^{-a\eta}$, $\beta_\eta=e^{-b\eta}$, or
	$\alpha_\eta=1-a\eta$, $\beta_\eta=1-b\eta$ for small $\eta$.
	The rates $a,b$, offset $\varepsilon$ and decay $\lambda$ stay fixed.
For $n\ge1$, the corrected moments satisfy the exact identities
	\[
	\begin{gathered}
		m_{n+1}-m_n
		=\frac{1-\alpha_\eta}{A_{n+1}}
		(g_{n+1}-m_n),\\
		v_{n+1}-v_n
		=\frac{1-\beta_\eta}{B_{n+1}}
		(g_{n+1}^{\odot2}-v_n).
	\end{gathered}
	\]
	For the first identity, equations~\eqref{eq:discrete-general-recursion}
	and~\eqref{eq:discrete-correction-factors} give
	$\widetilde m_n=A_nm_n$,
	$\widetilde m_{n+1}=\alpha_\eta\widetilde m_n+(1-\alpha_\eta)g_{n+1}$,
	and $A_{n+1}-\alpha_\eta A_n=1-\alpha_\eta$. Hence
	\[
	\begin{aligned}
	 m_{n+1}-m_n
	 &=\frac{\alpha_\eta\widetilde m_n
	 +(1-\alpha_\eta)g_{n+1}}{A_{n+1}}-m_n\\
	 &=\frac{\alpha_\eta A_nm_n
	 +(1-\alpha_\eta)g_{n+1}-A_{n+1}m_n}{A_{n+1}}\\
	 &=\frac{1-\alpha_\eta}{A_{n+1}}(g_{n+1}-m_n).
	\end{aligned}
	\]
	The identity for $v_n$ follows similarly, using
	$\widetilde v_n=B_nv_n$ and the second-moment recursion in
	~\eqref{eq:discrete-general-recursion}.
	By~\eqref{eq:memory-scaling}, $(1-\alpha_\eta)/\eta=a_\eta$
	and $(1-\beta_\eta)/\eta=b_\eta$. 	Let $\mathcal F_n=\sigma(\xi_1,\ldots,\xi_n)$ be the sample filtration. Taking conditional expectations
	for the moment updates and combining them with the position update gives
	\[
	\begin{gathered}
	\frac{x_{n+1}-x_n}{\eta}
	=-\frac{m_{n+1}}{\varepsilon+\sqrt{v_{n+1}}}-\lambda x_n,\\
	\frac{\mathbb E[m_{n+1}-m_n\mid\mathcal F_n]}{\eta}
	=\frac{a_\eta}{A_{n+1}}(\nabla F(x_n)-m_n),\\
	\frac{\mathbb E[v_{n+1}-v_n\mid\mathcal F_n]}{\eta}
	=\frac{b_\eta}{B_{n+1}}(S(x_n)-v_n).
	\end{gathered}
	\]
	
	As $\eta\downarrow0$ with $n\eta\to t$, we also have
	$(n+1)\eta\to t$. The correction factors therefore satisfy
	\begin{align}
A_n&=1-\alpha_\eta^n =1-\big[(1-\eta a_\eta)^{1/(\eta a_\eta)}\big]^{n\eta a_\eta} \notag \\
&\longrightarrow A(t):=1-e^{-at}, \notag \\
B_n&=1-\beta_\eta^n =1-\big[(1-\eta b_\eta)^{1/(\eta b_\eta)}\big]^{n\eta b_\eta} \notag \\
&\longrightarrow B(t):=1-e^{-bt}. \label{eq:correction-factors}
\end{align}
	Since $a_\eta\to a$ and $b_\eta\to b$, the moment coefficients converge to
	$a_\eta/A_{n+1}\to a/A(t)$ and
	$b_\eta/B_{n+1}\to b/B(t)$. 
	Thus the
	limiting corrected state $z=(x,m,v)$ satisfies, for $t>0$,
	\begin{equation}\label{eq:normalized}
	\begin{gathered}
	\dot x=-\frac{m}{\varepsilon+\sqrt{v}}-\lambda x,\\
	\dot m=\frac a{A(t)}(\nabla F(x)-m),\\
	\dot v=\frac b{B(t)}(S(x)-v).
	\end{gathered}
	\end{equation}
	\citet{CompanionStoch} proved that, under their local
	sampling assumptions, the stochastic iterates approach solutions of this
	ODE on each fixed bounded time interval. The position path is
	approximated on bounded intervals including $t=0$; approximation of the
	full corrected state holds away from zero, where $a/A(t)$ and $b/B(t)$
	are nonsingular. At $\lambda=0$ this is continuous Adam; the separate
	drift $-\lambda x$ accounts for AdamW decay.
	
The zero initialization of the raw discrete filters determines the generally
nonzero initial values of the corrected moments. Along a continuous position path with
$x(0)=x_0$, for $t>0$,
\begin{gather}
m(t)=\frac{a}{A(t)}\int_0^t e^{-a(t-s)}\nabla F(x(s))\,ds, \notag \\
v(t)=\frac{b}{B(t)}\int_0^t e^{-b(t-s)}S(x(s))\,ds.
\label{eq:corrected-moment-representation}
\end{gather}
Proposition~\ref{prop:integral-formulation} in
Appendix~\ref{guide:appendix-3-1} derives these formulas by integrating
the equations for $A(t)m(t)$ and $B(t)v(t)$.
Continuity of $\nabla F$ and $S$ at $x_0$, together with
$A(t)=at+o(t)$ and $B(t)=bt+o(t)$, gives, as $t\downarrow0$,
\begin{gather*}
m(t)=\frac{at\nabla F(x_0)+o(t)}{at+o(t)}\longrightarrow\nabla F(x_0), \\
v(t)=\frac{btS(x_0)+o(t)}{bt+o(t)}\longrightarrow S(x_0).
\end{gather*}
	Thus the continuous extension of the corrected state has initial value
	\begin{equation}\label{eq:normalized-initial-value}
	z(0)=(x_0,\nabla F(x_0),S(x_0)).
\end{equation}
	Because the coefficients $a/A(t)$ and $b/B(t)$ are singular at zero,
	by a solution we mean a continuous trajectory
	$z\in C([0,\infty))$ with $v(t)\ge0$ that satisfies the ODE
	classically for every $t>0$ and has the initial value displayed above.
	We refer to~\eqref{eq:normalized} together with
	\eqref{eq:normalized-initial-value} as the \emph{bias-corrected initial-value problem}.

	\subsection{Autonomous limit and equilibrium sets}\label{sec:autonomous-limit}
Let $z=(x,m,v)$ take values in the closed state space $\overline U$,
with interior $U$, where
\begin{gather*}
U:=\R^d\times\R^d\times(0,\infty)^d, \\
\overline U:=\R^d\times\R^d\times[0,\infty)^d.
\end{gather*}
Since $A(t),B(t)\to1$ as $t\to\infty$, replacing the bias-correction
factors in~\eqref{eq:normalized} by their limits yields the autonomous system
\begin{equation}\label{eq:auto}
\dot z=f_\lambda(z):=
\begin{pmatrix}
-m/(\varepsilon+\sqrt v)-\lambda x\\
a(\nabla F(x)-m)\\
b(S(x)-v)
\end{pmatrix}.
\end{equation}
We use the invariant sets of~\eqref{eq:auto} to study the long-time
behavior of the bias-corrected initial-value
problem~\eqref{eq:normalized}--\eqref{eq:normalized-initial-value}.

A state $(p,m,v)$ is an equilibrium of~\eqref{eq:auto} precisely when
all three components of $f_\lambda$ vanish. Since $a,b>0$, the moment
equations require $m=\nabla F(p)$ and $v=S(p)$. Substituting these
identities into the position equation gives, for each $i=1,\ldots,d$,
\[
\frac{\partial_i F(p)}{\varepsilon+\sqrt{S_i(p)}}+\lambda p_i=0.
\]
The denominator is at least $\varepsilon>0$, even when $S_i(p)=0$.
Thus this equation is equivalent to
\[
\partial_i F(p)+\lambda\bigl(\varepsilon+\sqrt{S_i(p)}\bigr)p_i=0.
\]
Conversely, any position $p$ satisfying these coordinatewise conditions
determines the unique equilibrium $(p,\nabla F(p),S(p))$.
This motivates the diagonal matrix and equilibrium sets
\begin{equation}\label{eq:adamw-equilibrium-set}
\begin{gathered}
D(x):=\diag\bigl(\varepsilon+\sqrt{S_i(x)}\bigr),\\
E_\lambda:=\{p\in\R^d:\nabla F(p)+\lambda D(p)p=0\},\\
\mathcal E_\lambda:=\{(p,\nabla F(p),S(p)):p\in E_\lambda\}.
\end{gathered}
\end{equation}
Here $E_\lambda$ is the set of equilibrium positions in $\R^d$,
whereas $\mathcal E_\lambda$ is the set of corresponding equilibria
in the full state space $\overline U$.
If $S_i(p)>0$ for every coordinate and every
$p\in E_\lambda$, then $\mathcal E_\lambda\subset U$.
Throughout, convergence to an equilibrium set is measured using
Euclidean distance to $E_\lambda$ or $\mathcal E_\lambda$, as appropriate.
For $\lambda=0$, the equilibrium condition reduces to $\nabla F(p)=0$,
so the equilibrium positions are exactly the critical points of $F$.
For $\lambda>0$, the adaptive gradient $D(p)^{-1}\nabla F(p)$ balances
the decay term $\lambda p$. We write
	\begin{gather}
\operatorname{Crit}(F):=\{q\in\R^d:\nabla F(q)=0\}, \notag \\
\mathcal  E_0:=\{(q,0,S(q)):q\in\operatorname{Crit}(F)\}. \label{eq:adam-equilibrium-set}
\end{gather}
For $z\in U$, let $\Phi_\lambda(t,z)$ denote the solution of
$\dot z=f_\lambda(z)$ at time $t$ with initial state $z$, for times on
its maximal interval of existence in $U$.
	The critical set $\operatorname{Crit}(F)$ is the zero-decay equilibrium set in
	position space. The set $\mathcal E_0$ consists of the corresponding
	equilibria in the full state space.

	\section{Main results}\label{sec:main-results}
	Throughout, $a,b,\varepsilon>0$, $\lambda\ge0$, and
	moment domination~\eqref{eq:moment-domination} holds. All vector
	operations are coordinatewise unless stated otherwise. Distances and
	vector norms are Euclidean, matrix norms are the induced operator norms, and coercivity means $F(x)\to\infty$ as
	$\|x\|\to\infty$. Every small-decay threshold may depend on the fixed functions $F,S$,
	the coefficients $a,b,\varepsilon$, and the deterministic initial
	position.
Thresholds in different results need not be equal: we use
$\lambda_{\mathrm{man}}$ for manifold convergence,
$\lambda_{\mathrm{eq}}$ for equilibrium convergence,
$\lambda_{\mathrm{ql}}$ for quantitative localization, and
$\lambda_{\mathrm{ex}}$ for the nonconvergent example.
The constant $C_{\mathrm{ql}}$ belongs to the quantitative localization
bound. Auxiliary constants are local to the statement or proof in which
they are introduced, unless explicitly carried forward from another result;
reusing a symbol does not assert equality of those constants.
When several results are used together, the decay parameter must satisfy
all their thresholds, equivalently their minimum.
	\subsection{The singular initial time}\label{sec:main-wellposedness}
We first establish that the corrected initial values determine a unique
global trajectory despite the singular coefficients at $t=0$, even when
some initial second moments vanish. The proof is given in
Appendix~\ref{sec:wp-proof}.
	\begin{theorem}[Global well-posedness]\label{thm:adamw-wp}
		If $\nabla F$ and $\sqrt S:\R^d\to[0,\infty)^d$ are locally Lipschitz,
		then for every $x_0\in\R^d$, system~\eqref{eq:normalized} has a unique
		global solution $z=(x,m,v)$ with initial value~\eqref{eq:normalized-initial-value} and $v\ge0$.
		The solution is $C^1$ on $(0,\infty)$, and $x\in C^1([0,\infty))$.
	\end{theorem}

	\subsection{Convergence to invariant manifolds}
	\label{sec:long-time-assumptions}

The long-time results use the following assumptions.
\begin{assumption}\label{ass:long-time}
\leavevmode
\begin{enumerate}[label=(\roman*),leftmargin=*,itemsep=0pt,parsep=0pt,topsep=2pt]
\item The memory rates satisfy $0<b\le4a$.
\item The objective $F\in C^1(\R^d)$ is coercive, its gradient
$\nabla F$ is locally Lipschitz, and the second-moment map
$S:\R^d\to(0,\infty)^d$ is locally Lipschitz.
\item On an open neighborhood of $\operatorname{Crit}(F)$, the objective
$F$ is twice continuously differentiable and $S$ is continuously
differentiable.
\end{enumerate}
\end{assumption}
The roles and scope of these conditions are discussed in
Remark~\ref{rem:long-time-assumptions} in Appendix~\ref{sec:invariant-manifolds}.
	Throughout the long-time analysis, fix $x_0\in\R^d$. For each
	$\lambda\ge0$, let $z^\lambda$ solve~\eqref{eq:normalized} with the
	initial value~\eqref{eq:normalized-initial-value}.
	The system data $F,S,a,b,\varepsilon$ remain fixed while $\lambda$ varies.
	Small-decay thresholds may depend on these data and $x_0$.
	
	We first prove uniform long-time localization near the Adam equilibria.
	The proof constructs a common compact trapping region for small decay
	and is given in Appendix~\ref{sec:compact-limit-geometry}.
	
	\begin{theorem}[Small-decay localization]
		\label{thm:full-state-localization}
		Under Assumption~\ref{ass:long-time}(i)--(ii), if $F(\operatorname{Crit}(F))$ has empty
		interior, then
		\begin{equation}\label{eq:full-state-localization}
			\lim_{\lambda\downarrow0}\limsup_{t\to\infty}
			\operatorname{dist}(z^\lambda(t),\mathcal E_0)=0.
		\end{equation}
	\end{theorem}
For sufficiently small weight decay, the trajectory eventually stays near
the zero-decay equilibrium set. Reducing the decay makes this eventual
distance arbitrarily small; the time of entry may depend on the decay.
The critical-value assumption allows nonisolated critical points.

\begin{definition}[Locally invariant manifold]
\label{def:local-invariant-manifold}
Fix $\lambda\ge0$. A $C^1$ embedded manifold $\mathcal I_\lambda\subset U$
is locally invariant under~\eqref{eq:auto} if, for every point $p$
in the manifold interior of $\mathcal I_\lambda$, there exists $\tau_p>0$
such that
\[
\Phi_\lambda(t,p)\in\mathcal I_\lambda
\qquad\text{\normalfont for all }|t|<\tau_p.
\]
Thus the solution through each interior point remains on $\mathcal I_\lambda$
for a time interval around its initial time. Local invariance does not
require the solution to remain on $\mathcal I_\lambda$ for all time.
\end{definition}
	
	The next theorem concerns each fixed small $\lambda$ under the
	following geometric hypothesis.
	See Remark~\ref{rem:critical-manifold-assumption} in
	Appendix~\ref{sec:invariant-manifolds} for examples and the scope of this hypothesis.
	\begin{assumption}[Critical-manifold hypothesis]\label{ass:critical-manifold}
	The connected components $\mathcal M$ of $\operatorname{Crit}(F)$ are
	$C^1$ embedded manifolds without boundary. Every compact set meets
	only finitely many components, and the Hessian nullity is constant on each:
	\begin{equation}\label{eq:componentwise-nullity}
		\dim\ker\nabla^2F(q)=k_{\mathcal M}\qquad(q\in\mathcal M).
	\end{equation}
	The integers $k_{\mathcal M}$ may vary between components.
	\end{assumption}

\begin{figure}[t]
\centering
\includegraphics[width=\linewidth]{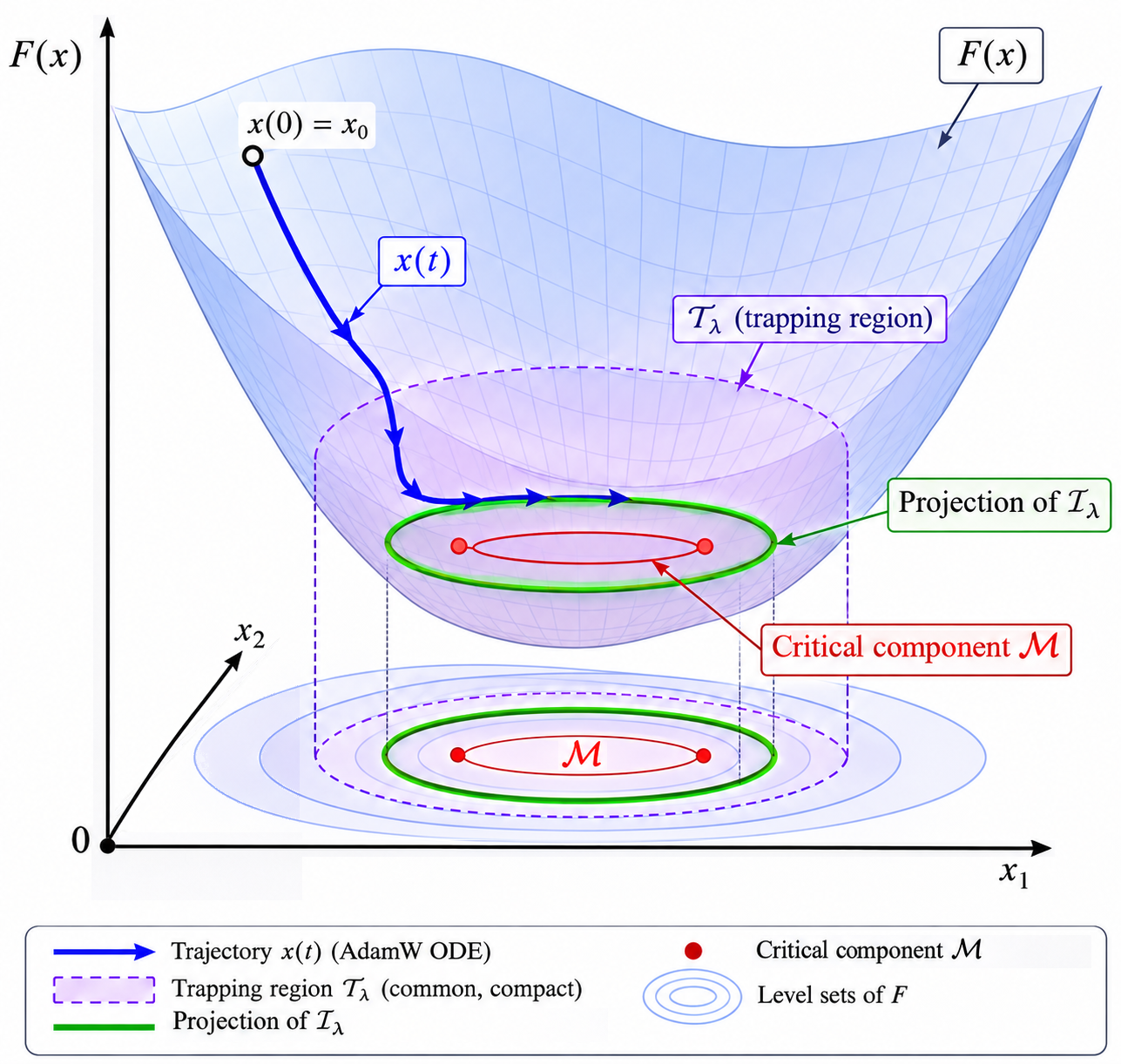}
\captionsetup{font=small,labelfont=it}
\hypersetup{allcolors=ink}
\caption{
Geometric illustration of the long-time behavior of continuous-time AdamW.
An illustrative objective is $F(x)=(x_1^2+x_2^2-1)^2$.
The blue curve depicts the position trajectory $x(t)$, the red circle
denotes a connected component $\mathcal M$ of the critical set,
and the purple region represents a compact trapping region.
The green curve shows a schematic projection of the locally invariant
manifold $\mathcal I_\lambda$ onto position space.
Under Assumption~\ref{ass:critical-manifold}, each connected component
$\mathcal M$ of $\operatorname{Crit}(F)$ is a $C^1$ embedded manifold with constant
Hessian nullity $k_{\mathcal M}$. Under the hypotheses of
Theorem~\ref{thm:invariant-manifold}, for sufficiently small fixed weight
decay $\lambda>0$, the full state $z^\lambda(t)$ remains in a compact region
and satisfies $\operatorname{dist}(z^\lambda(t),\mathcal I_\lambda)\to0$
as $t\to\infty$, where $\mathcal I_\lambda$ is a
$k_{\mathcal M}$-dimensional locally invariant manifold of the limiting
autonomous system. 
}
\label{fig:adamw-long-time-geometry}
\end{figure}

The proof of the following result is given in
Appendix~\ref{sec:manifold-limit-sets}.
\begin{theorem}[Manifold convergence]
		\label{thm:invariant-manifold}
		Under Assumption~\ref{ass:long-time}(i)--(iii) and
		Assumption~\ref{ass:critical-manifold}, there exists
		$\lambda_{\mathrm{man}}>0$ such that, for every $0<\lambda<\lambda_{\mathrm{man}}$, there are
		a connected component $\mathcal M$ of $\operatorname{Crit}(F)$ and a
		$k_{\mathcal M}$-dimensional $C^1$ locally invariant manifold
		$\mathcal I_\lambda$ for~\eqref{eq:auto} satisfying
		\begin{equation}\label{eq:invariant-localization}
			\operatorname{dist}\bigl(z^\lambda(t),\mathcal I_\lambda\bigr)
			\longrightarrow0\qquad(t\to\infty).
		\end{equation}
	\end{theorem}
By Theorem~\ref{thm:invariant-manifold}, the position and moments
approach a locally invariant manifold near a critical component.
Its dimension equals the Hessian nullity, which can be positive even
at an isolated critical point. The trajectory may settle at a point
or continue moving near the manifold, as illustrated in
Figure~\ref{fig:adamw-long-time-geometry}.

	\subsection{Convergence to equilibria}\label{sec:main-equilibria}
The next theorem gives a positive point-convergence result for the
original bias-corrected AdamW dynamics. When the critical points are
isolated and their Hessian nullity is at most one, the full state
$z^\lambda(t)$ converges to a single equilibrium for each sufficiently
small fixed positive decay under the stated assumptions. It strengthens
the invariant-manifold conclusion without modifying the optimizer.
The proof is given in Appendix~\ref{sec:equilibrium-proofs}.

	\begin{theorem}[Equilibrium convergence]
		\label{thm:corank-one-convergence}
		Suppose Assumption~\ref{ass:long-time}(i)--(iii) holds, the critical
points are isolated, and $\dim\ker\nabla^2F(q)\le1$ for every
$q\in\operatorname{Crit}(F)$.
		There is $\lambda_{\mathrm{eq}}>0$ such that, for each fixed $0<\lambda<\lambda_{\mathrm{eq}}$,
		some $p_\lambda\in E_\lambda$ satisfies
		\begin{equation}\label{eq:corank-one-equilibrium-limit}
			z^\lambda(t)\longrightarrow
			(p_\lambda,\nabla F(p_\lambda),S(p_\lambda))
			\qquad(t\to\infty).
		\end{equation}
	\end{theorem}
Under these conditions, the position and both moments settle at fixed
values for sufficiently small decay. The limiting position balances the
adaptive gradient update against weight decay:
$\nabla F(p_\lambda)+\lambda D(p_\lambda)p_\lambda=0$.
This includes objectives whose critical points are all nondegenerate.

	\par\addvspace{0.5\baselineskip}
\begin{remark*}[The role of the nullity restriction]
Example~\ref{ex:isolated-corank-two} shows that the nullity-at-most-one
condition cannot simply be dropped from the point-convergence theorem.
With Hessian nullity two, a smooth, strictly convex, coercive objective
with a unique minimum can support a bounded nonconvergent trajectory
for every sufficiently small positive decay, with $F,S$ and the initial
position fixed. The two-dimensional kernel leaves a degenerate plane
on which state-dependent second moments generate rotation. This is a
possible mechanism at nullity two, not a claim that every such objective
produces persistent motion.
\end{remark*}

\subsection{Quantitative localization without point convergence}\label{sec:main-quantitative}
\paragraph{Local quantitative geometry.}
Fix $0<\theta\le1/2$ and $0<\gamma\le1$. At a critical point $q$,
let $U_q$ be a neighborhood of $q$. For $x\in U_q$, consider the local bounds
\begin{gather}
|F(x)-F(q)|^{1-\theta}\le k_q\|\nabla F(x)\|, \notag \\
\operatorname{dist}(x,\operatorname{Crit}(F))\le M_q\|\nabla F(x)\|^\gamma, \label{eq:local-quantitative-geometry}
\end{gather}
Here $k_q$ and $M_q$ are positive constants that may depend on the
critical point $q$. The first bound in~\eqref{eq:local-quantitative-geometry}
is the gradient \L{}ojasiewicz inequality: it controls the local objective
gap through the gradient norm. The second bound is a distance error
bound: it controls the distance to $\operatorname{Crit}(F)$ through the
gradient norm. The parameters $\theta$ and $\gamma$, which determine the
powers $1-\theta$ and $\gamma$ in the respective inequalities, are required
to be the same at every critical point; the constants $k_q,M_q$ and
neighborhoods $U_q$ may vary with $q$. These bounds impose local
quantitative geometry without requiring isolated critical points.

\begin{theorem}[Quantitative localization]
  \label{thm:quantitative-localization}
  Suppose Assumption~\ref{ass:long-time}(i)--(ii) holds and that \eqref{eq:local-quantitative-geometry}
  holds near every $q\in\operatorname{Crit}(F)$. Then there are
  $C_{\mathrm{ql}},\lambda_{\mathrm{ql}}>0$ such that, for $0<\lambda<\lambda_{\mathrm{ql}}$,
  \begin{align}
&\limsup_{t\to\infty}\operatorname{dist}(x^\lambda(t),\operatorname{Crit}(F)) \notag \\
&\le C_{\mathrm{ql}}\left(\lambda^\gamma+ 
       \lambda^{\frac{\gamma}{(1+\gamma)(1-\theta)}}\right). \label{eq:quantitative-localization}
\end{align}
\end{theorem}
The proof of the above estimate is given in
Appendix~\ref{sec:quantitative-excursions}.
The bound quantifies how reducing weight decay brings the late-time
position closer to the original critical set, including trajectories
that continue to move.
The constants may depend on $F,S,a,b,\varepsilon,x_0$ and the geometric
bounds, but not on $\lambda$.

\begin{figure}[htbp]
\centering
\includegraphics[width=\linewidth]{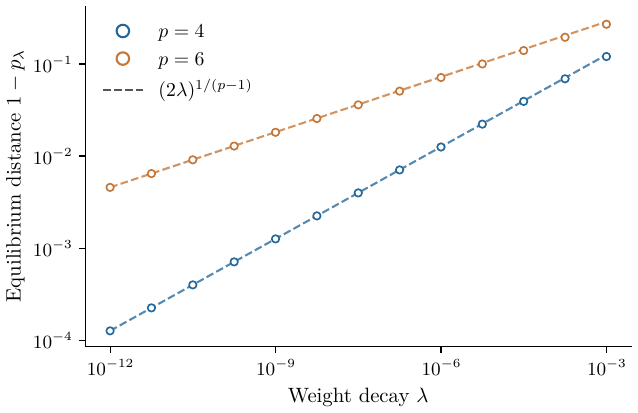}
\caption{The distance $1-p_\lambda$ from the AdamW equilibrium to the
minimizer in Example~\ref{ex:quantitative-polynomial-sharp}.
Markers are computed from the scalar equilibrium equation. Dashed lines
show $(2\lambda)^{1/(p-1)}$, with slopes $1/3$ for $p=4$ and $1/5$ for
$p=6$ on the log--log axes.}
\label{fig:localization-sharpness}
\end{figure}

\begin{remark*}[Sharpness of the bound]
For the polynomial objectives in Example~\ref{ex:quantitative-polynomial-sharp},
the limiting AdamW position satisfies
$1-p_\lambda\sim(2\lambda)^{1/(p-1)}$ as $\lambda\downarrow0$.
Both terms in Theorem~\ref{thm:quantitative-localization} have order
$\lambda^{1/(p-1)}$ for these examples, so the bound is sharp for
this family. Figure~\ref{fig:localization-sharpness} illustrates
the asymptotic law for $p=4$ and $p=6$.
\end{remark*}

\FloatBarrier

\section{Hessian nullity two: nonconvergence}
\label{sec:counterexample-visualization}
\label{sec:isolated-corank-two-counterexample}
To show why the nullity restriction in
Theorem~\ref{thm:corank-one-convergence} matters, we give fixed polynomial
functions $F,S$ and one fixed initial position for which the corrected
trajectory is bounded and nonconvergent for every sufficiently small
positive decay. The objective is smooth, strictly convex, and coercive,
with a unique minimum of Hessian nullity two. The degenerate plane
allows the second-moment field to generate tangential motion, despite
the uniqueness of the minimum.

Place the minimum at $c=(1,1,1)$ and let $n=(1,1,1)/\sqrt3$
be the unit vector along the diagonal.
For $x\in\mathbb R^3$, write
\begin{gather*}
\zeta:=n\cdot(x-c) =\frac{x_1+x_2+x_3-3}{\sqrt3}, \\
y:=x-c-\zeta n.
\end{gather*}
Thus $\zeta$ is the signed coordinate along $n$ and
$x-c=y+\zeta n$, where $y$ lies in the perpendicular plane
$n^\perp=\{u:u_1+u_2+u_3=0\}$.
The point $c+y$ is the orthogonal projection of $x$ onto $c+n^\perp$.

Define
\begin{equation}\label{eq:counterexample-objective}
F(x):=\frac14\|y\|^4+\frac12\zeta^2.
\end{equation}
The objective grows quartically in that plane and quadratically along
$n$.
The vector $n\times y$ is $y$ rotated by $90^\circ$ in $n^\perp$,
with orientation given by $n$.
For $i=1,2,3$, set
\begin{equation}\label{eq:visualization-second-moments}
S_i(x):=1+(\partial_iF(x))^2
       +3\bigl(1-3y_i+(n\times y)_i\bigr)^2.
\end{equation}
Appendix~\ref{app:counterexample-construction} verifies the assumptions
and shows that $F$ is strictly convex and coercive, with a unique
critical point $c$ and Hessian nullity two at $c$.

\begin{example}[Nonconvergence at corank two]
\label{ex:isolated-corank-two}
Let $F,S$ be defined in~\eqref{eq:counterexample-objective}
and~\eqref{eq:visualization-second-moments}. Set $a=b=1$ and fix
$0<\varepsilon\le1$. There exist constants $s,\lambda_{\mathrm{ex}}>0$,
which may depend on $\varepsilon$, such that the following holds.
Fix the single initial position
$x_0=c+s(1,-1,0)/\sqrt2$. For every $0<\lambda<\lambda_{\mathrm{ex}}$, the
solution of~\eqref{eq:normalized}--\eqref{eq:normalized-initial-value}
with this $x_0$ is bounded, but $x^\lambda(t)$ does not converge as
$t\to\infty$. The functions $F,S$ and the initial position are
independent of $\lambda$.
\end{example}
The proof of nonconvergence is given in Appendix~\ref{app:counterexample-proof}.
Figure~\ref{fig:counterexample-trajectories} shows the three-dimensional
position paths for three small decay values. After an initial transient,
the paths repeatedly trace nearly closed curves. Their distances to the
nearby equilibrium $p_\lambda$ oscillate and stay positive over the
displayed time interval, with time rescaled as $\tau=\lambda t$.
The simulations suggest attraction to periodic orbits.
Example~\ref{ex:isolated-corank-two} establishes nonconvergence.

\begin{figure}[t!]
	\centering
	\includegraphics[width=.84\linewidth]{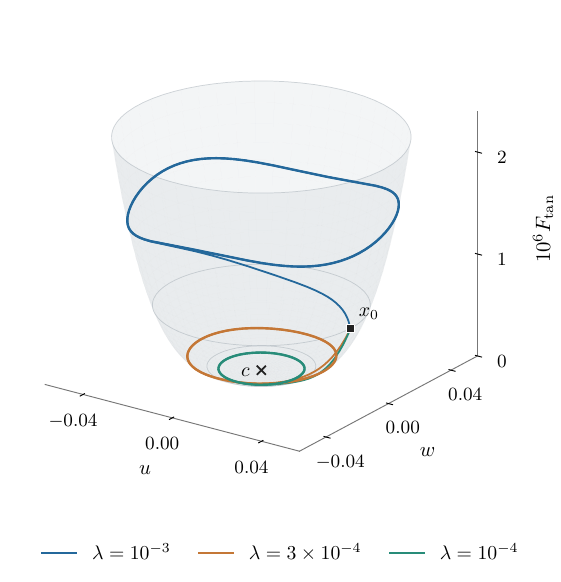}
	\caption{Position trajectories projected onto $c+n^\perp$ and drawn on
		the graph of $F_{\rm tan}$, with the vertical axis scaled by $10^6$.
		The square marks the common initial position, the cross marks $c$,
		and thicker lines indicate $10\le\tau\le20$.}
	\label{fig:counterexample-loss-3d}
	\label{fig:counterexample-projection}
\end{figure}

Figure~\ref{fig:counterexample-loss-3d} shows the same paths projected
onto $c+n^\perp$. The loops at different decay values are visible at
different heights on the graph of the restricted objective
$F_{\rm tan}(u,w)=(u^2+w^2)^2/4$, where $u,w$ are orthonormal
coordinates in that plane. Appendix~\ref{app:counterexample-simulation}
lists the simulation parameters and the nearby equilibrium positions.

\begin{figure}[t!]
	\centering
	\includegraphics[width=\linewidth]{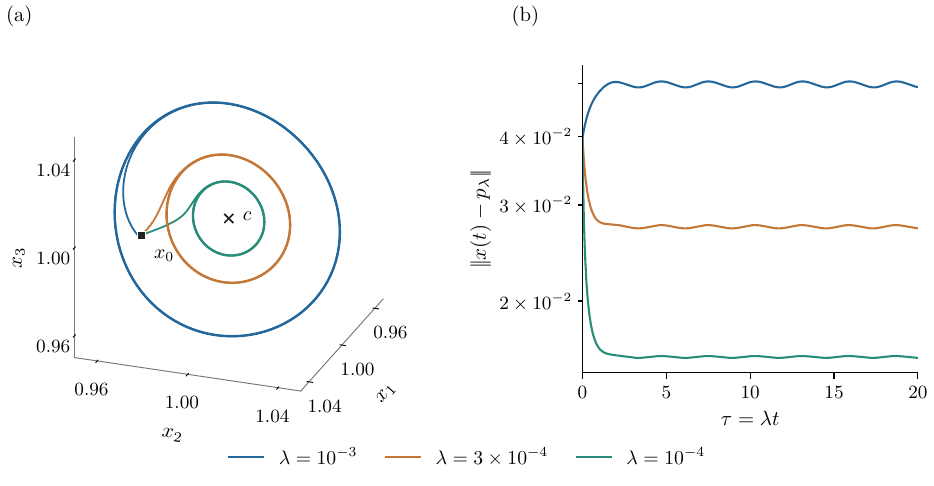}
	\caption{Numerical trajectories and distances to the nearby equilibrium.
		(a) Full position trajectories in $(x_1,x_2,x_3)$ with equal axis scales.
		The square marks the common initial position, the cross marks $c$,
		and thicker lines indicate $10\le\tau\le20$.
		(b) Distance $\|x(t)-p_\lambda\|$ on a logarithmic vertical scale for
		$0\le\tau\le20$.}
	\label{fig:counterexample-trajectories}
\end{figure}

\FloatBarrier

\section*{Acknowledgements}
This work was supported in part by a grant of access to OpenAI models through the ChatGPT for Academic Researchers program.

\clearpage
\label{sec:references}
\bibliographystyle{icml2026}
\bibliography{references}

@unpublished{CompanionStoch,
  author = {Wang, B. and Ning, N.},
  title = {Long-Time Behavior of Stochastic {AdamW}},
  year = {2026},
  note = {Manuscript},
}

@inproceedings{KB2015,
  author = {Kingma, Diederik P. and Ba, Jimmy},
  title = {{Adam}: A method for stochastic optimization},
  year = {2015},
  booktitle = {International Conference on Learning Representations},
  url = {https://arxiv.org/abs/1412.6980},
}

@inproceedings{LH2019,
  author = {Loshchilov, Ilya and Hutter, Frank},
  title = {Decoupled weight decay regularization},
  year = {2019},
  booktitle = {International Conference on Learning Representations},
  url = {https://arxiv.org/abs/1711.05101},
}

@inproceedings{RKK2018,
  author = {Reddi, Sashank J. and Kale, Satyen and Kumar, Sanjiv},
  title = {On the convergence of {Adam} and beyond},
  year = {2018},
  booktitle = {International Conference on Learning Representations},
  url = {https://arxiv.org/abs/1904.09237},
}

@article{BB2021,
  author = {Barakat, Anas and Bianchi, Pascal},
  title = {Convergence and dynamical behavior of the {ADAM} algorithm for nonconvex stochastic optimization},
  year = {2021},
  journal = {SIAM Journal on Optimization},
  volume = {31},
  pages = {244--274},
  url = {https://doi.org/10.1137/19M1263443},
  doi = {10.1137/19M1263443},
  number = {1},
}

@article{BG2020,
  author = {Belotto da Silva, Alexandre and Gazeau, Maxime},
  title = {A general system of differential equations to model first-order adaptive algorithms},
  year = {2020},
  journal = {Journal of Machine Learning Research},
  volume = {21},
  number = {129},
  pages = {1--42},
  url = {https://www.jmlr.org/papers/v21/18-808.html},
}

@article{BBHS2021,
  author = {Barakat, Anas and Bianchi, Pascal and Hachem, Walid and Schechtman, Sholom},
  title = {Stochastic optimization with momentum: Convergence, fluctuations, and traps avoidance},
  year = {2021},
  journal = {Electronic Journal of Statistics},
  volume = {15},
  pages = {3892--3947},
  url = {https://doi.org/10.1214/21-EJS1880},
  doi = {10.1214/21-EJS1880},
}

@misc{GT2024,
  author = {Gould, R. and Tanaka, H.},
  title = {Continuous-time analysis of adaptive optimization and normalization},
  year = {2024},
  howpublished = {arXiv preprint arXiv:2411.05746},
  url = {https://arxiv.org/abs/2411.05746},
}

@misc{DJK2025,
  author = {Dereich, Steffen and Jentzen, Arnulf and Kassing, Sebastian},
  title = {{ODE} approximation for the {Adam} algorithm: General and overparametrized setting},
  year = {2025},
  howpublished = {arXiv preprint arXiv:2511.04622},
  url = {https://arxiv.org/abs/2511.04622},
}

@article{TD2017,
  author = {Tadi{\'c}, Vladislav B. and Doucet, Arnaud},
  title = {Asymptotic bias of stochastic gradient search},
  year = {2017},
  journal = {The Annals of Applied Probability},
  volume = {27},
  pages = {3255--3304},
  url = {https://doi.org/10.1214/16-AAP1272},
  doi = {10.1214/16-AAP1272},
  number = {6},
}

@article{BLMP2025,
  author = {Bolte, J{\'e}r{\^o}me and Le, Tam and Moulines, {\'E}ric and Pauwels, Edouard},
  title = {Inexact subgradient methods for semialgebraic functions},
  year = {2026},
  url = {https://doi.org/10.1007/s10107-025-02245-w},
  journal = {Mathematical Programming},
  volume = {219},
  pages = {1--27},
  doi = {10.1007/s10107-025-02245-w},
}

@article{RB2025,
  author = {Rebjock, Quentin and Boumal, Nicolas},
  title = {Fast convergence to non-isolated minima: Four equivalent conditions for {$C^2$} functions},
  year = {2025},
  journal = {Mathematical Programming},
  volume = {213},
  pages = {151--199},
  url = {https://doi.org/10.1007/s10107-024-02136-6},
  doi = {10.1007/s10107-024-02136-6},
}

@incollection{B1999,
  author = {Bena\"im, M.},
  title = {Dynamics of stochastic approximation algorithms},
  year = {1999},
  booktitle = {S\'eminaire de probabilit\'es XXXIII},
  series = {Lecture Notes in Mathematics},
  volume = {1709},
  pages = {1--68},
  publisher = {Springer},
  url = {https://www.numdam.org/item/SPS_1999__33__1_0/},
}

@article{Fenichel1979,
  author = {Fenichel, N.},
  title = {Geometric singular perturbation theory for ordinary differential equations},
  year = {1979},
  journal = {Journal of Differential Equations},
  volume = {31},
  pages = {53--98},
  url = {https://doi.org/10.1016/0022-0396(79)90152-9},
  doi = {10.1016/0022-0396(79)90152-9},
}

@misc{LiuDeepSeek2024,
  author = {Liu, A. and Feng, B. and Xue, B. and Wang, B. and Wu, B. and Lu, C. and Zhao, C. and Deng, C. and Zhang, C. and Ruan, C. and others},
  title = {{DeepSeek-V3} technical report},
  year = {2024},
  howpublished = {arXiv preprint arXiv:2412.19437},
  url = {https://arxiv.org/abs/2412.19437},
}

@misc{JordanMuon2024,
  author = {Jordan, K. and Jin, Y. and Boza, V. and You, J. and Cesista, F. and Newhouse, L. and Bernstein, J.},
  title = {{Muon}: An optimizer for hidden layers in neural networks},
  year = {2024},
  howpublished = {Online article},
  url = {https://kellerjordan.github.io/posts/muon/},
}

@inproceedings{BockWeiss2019local,
  author = {Bock, Sebastian and Wei{\ss}, Martin},
  title = {A proof of local convergence for the {Adam} optimizer},
  year = {2019},
  booktitle = {International Joint Conference on Neural Networks (IJCNN)},
  publisher = {IEEE},
  pages = {1--8},
  url = {https://doi.org/10.1109/IJCNN.2019.8852239},
  doi = {10.1109/IJCNN.2019.8852239},
}

@inproceedings{BockWeiss2019cycles,
  author = {Bock, Sebastian and Wei{\ss}, Martin},
  title = {Non-convergence and limit cycles in the {Adam} optimizer},
  year = {2019},
  booktitle = {Artificial Neural Networks and Machine Learning---ICANN 2019: Deep Learning},
  series = {Lecture Notes in Computer Science},
  volume = {11728},
  publisher = {Springer},
  pages = {232--243},
  url = {https://doi.org/10.1007/978-3-030-30484-3_20},
  doi = {10.1007/978-3-030-30484-3_20},
}

@article{Defossez2020,
  author = {D{\'e}fossez, Alexandre and Bottou, L{\'e}on and Bach, Francis and Usunier, Nicolas},
  title = {A simple convergence proof of {Adam} and {Adagrad}},
  year = {2022},
  url = {https://openreview.net/forum?id=ZPQhzTSWA7},
  journal = {Transactions on Machine Learning Research},
}

@inproceedings{Zou2019,
  author = {Zou, Fangyu and Shen, Li and Jie, Zequn and Zhang, Weizhong and Liu, Wei},
  title = {A sufficient condition for convergences of {Adam} and {RMSProp}},
  year = {2019},
  booktitle = {Proceedings of the IEEE/CVF Conference on Computer Vision and Pattern Recognition},
  pages = {11127--11135},
}

@inproceedings{ZhangAdam2022,
  author = {Zhang, Yushun and Chen, Congliang and Shi, Naichen and Sun, Ruoyu and Luo, Zhi-Quan},
  title = {{Adam} can converge without any modification on update rules},
  year = {2022},
  booktitle = {Advances in Neural Information Processing Systems},
  volume = {35},
  pages = {28386--28399},
}

@inproceedings{GronichVardi2026,
  author = {Gronich, Eitan and Vardi, Gal},
  title = {The implicit bias of {Adam} and {Muon} on smooth homogeneous neural networks},
  year = {2026},
  booktitle = {Proceedings of the 43rd International Conference on Machine Learning},
  url = {https://proceedings.mlr.press/v306/gronich26a.html},
  series = {Proceedings of Machine Learning Research},
  volume = {306},
  pages = {36767--36811},
  publisher = {PMLR},
}

@inproceedings{WangMeng2021,
  author = {Wang, Bohan and Meng, Qi and Chen, Wei and Liu, Tie-Yan},
  title = {The implicit bias for adaptive optimization algorithms on homogeneous neural networks},
  year = {2021},
  booktitle = {Proceedings of the International Conference on Machine Learning},
  series = {Proceedings of Machine Learning Research},
  volume = {139},
  pages = {10849--10858},
  url = {https://proceedings.mlr.press/v139/wang21q.html},
}

@inproceedings{XieLi2024,
  author = {Xie, Shichun and Li, Zhiyuan},
  title = {Implicit bias of {AdamW}: $\ell_\infty$-norm constrained optimization},
  year = {2024},
  booktitle = {Proceedings of the International Conference on Machine Learning},
  series = {Proceedings of Machine Learning Research},
  volume = {235},
  pages = {54488--54510},
  url = {https://proceedings.mlr.press/v235/xie24e.html},
}

@article{zhou2024understanding,
  title = {Towards Understanding Convergence and Generalization of {AdamW}},
  author = {Zhou, Pan and Xie, Xingyu and Lin, Zhouchen and Yan, Shuicheng},
  journal = {IEEE Transactions on Pattern Analysis and Machine Intelligence},
  volume = {46},
  number = {9},
  pages = {6486--6493},
  year = {2024},
  doi = {10.1109/TPAMI.2024.3382294},
}

@inproceedings{cattaneo2024implicit,
  title = {On the Implicit Bias of {Adam}},
  author = {Cattaneo, Matias D. and Klusowski, Jason Matthew and Shigida, Boris},
  booktitle = {Proceedings of the 41st International Conference on Machine Learning},
  series = {Proceedings of Machine Learning Research},
  volume = {235},
  pages = {5862--5906},
  year = {2024},
  publisher = {PMLR},
}

@inproceedings{SuChenLiu2026Stability,
  author = {Su, Baiyu and Chen, Lizhang and Liu, Qiang},
  title = {A Stability Analysis of {AdamW}: Unstable Equilibria and Non-Convergence},
  booktitle = {Advances in Neural Information Processing Systems (NeurIPS)},
  year = {2026},
  url = {https://neurips.cc/virtual/2026/poster/151341},
}

@misc{RaviSAM2024,
  author = {Ravi, Nikhila and Gabeur, Valentin and Hu, Yuan-Ting and Hu, Ronghang and Ryali, Chaitanya and Ma, Tengyu and Khedr, Haitham and R{\"a}dle, Roman and Rolland, Chloe and Gustafson, Laura and others},
  title = {{SAM} 2: Segment Anything in Images and Videos},
  year = {2024},
  howpublished = {arXiv preprint arXiv:2408.00714},
  url = {https://arxiv.org/abs/2408.00714},
}

@inproceedings{Esser2024,
  author = {Esser, Patrick and Kulal, Sumith and Blattmann, Andreas and Entezari, Rahim and M{\"u}ller, Jonas and Saini, Harry and Levi, Yam and Lorenz, Dominik and Sauer, Axel and Boesel, Frederic and Podell, Dustin and Dockhorn, Tim and English, Zion and Rombach, Robin},
  title = {Scaling Rectified Flow Transformers for High-Resolution Image Synthesis},
  booktitle = {Proceedings of the 41st International Conference on Machine Learning},
  year = {2024},
  series = {Proceedings of Machine Learning Research},
  volume = {235},
  pages = {12606--12633},
  publisher = {PMLR},
  url = {https://proceedings.mlr.press/v235/esser24a.html},
}

@misc{ZhangTEON2026,
  author = {Zhang, Ruijie and Zhao, Yequan and Liu, Ziyue and Wang, Zhengyang and Li, Dongyang and Su, Yupeng and Liu, Sijia and Zhang, Zheng},
  title = {{TEON}: Tensorized Orthonormalization Beyond Layer-Wise {Muon} for Large Language Model Pre-Training},
  year = {2026},
  howpublished = {arXiv preprint arXiv:2601.23261},
  url = {https://arxiv.org/abs/2601.23261},
}

@misc{TangTUON2026,
  author = {Tang, Anda and Yang, Xuanlin and Wang, Junjie and Lin, Zhouchen},
  title = {{TUON}: A Tensor-Aware Optimizer Utilizing Tensor-{SVD} and Polarization for Large Language Model Pretraining},
  year = {2026},
  howpublished = {arXiv preprint arXiv:2610.35130},
  url = {https://arxiv.org/abs/2610.35130},
}

@book{Kato1995,
  author = {Kato, Tosio},
  title = {Perturbation Theory for Linear Operators},
  publisher = {Springer},
  address = {Berlin},
  year = {1995},
  doi = {10.1007/978-3-642-66282-9}
}

@book{Khalil2002,
  author = {Khalil, Hassan K.},
  title = {Nonlinear Systems},
  edition = {3},
  publisher = {Prentice Hall},
  year = {2002}
}

\clearpage
\onecolumn
\appendix

\section*{Guide to the appendices}
The appendices provide the background, derivations, and proofs supporting
the main results, together with details of the nonconvergent example.
\begin{itemize}[leftmargin=*,itemsep=3pt,parsep=0pt]
\item \textbf{Appendix~\ref{app:comparison}:} Comparison with existing continuous-time Adam and AdamW analyses.
\begin{itemize}[leftmargin=3.5em,itemsep=0pt,parsep=0pt,topsep=2pt]
\item[\ref{app:comparison-scope}] Scope of the comparison.
\item[\ref{guide:appendix-1-1}] Adam convergence and AdamW stability.
\item[\ref{guide:appendix-1-2}] Implicit bias and normalized directions.
\item[\ref{guide:appendix-1-3}] Constrained interpretations and equilibria.
\item[\ref{guide:appendix-1-4}] Time scaling.
\end{itemize}
\item \textbf{Appendix~\ref{sec:wp}:} Well-posedness of the bias-corrected initial-value problem.
\begin{itemize}[leftmargin=3.5em,itemsep=0pt,parsep=0pt,topsep=2pt]
\item[\ref{guide:appendix-3-1}] Moment integral formulas, initial values, and the position equation.
\item[\ref{guide:appendix-3-2}] Boundedness.
\item[\ref{guide:appendix-3-3}] Lipschitz continuity.
\item[\ref{guide:appendix-3-4}] Proof of Theorem~\ref{thm:adamw-wp}.
\end{itemize}
\item \textbf{Appendix~\ref{sec:invariant-manifolds}:} Localization and convergence to invariant manifolds.
\begin{itemize}[leftmargin=3.5em,itemsep=0pt,parsep=0pt,topsep=2pt]
\item[\ref{sec:long-time-assumption-scope}] Assumptions and their scope.
\item[\ref{guide:appendix-4-1}] Energy estimates.
\item[\ref{guide:appendix-4-2}] Proof of Theorem~\ref{thm:full-state-localization}.
\item[\ref{guide:appendix-4-3}] Local invariant manifolds.
\item[\ref{guide:appendix-4-4}] Proof of Theorem~\ref{thm:invariant-manifold}.
\end{itemize}
\item \textbf{Appendix~\ref{sec:equilibrium-convergence}:} Convergence to equilibria.
\begin{itemize}[leftmargin=3.5em,itemsep=0pt,parsep=0pt,topsep=2pt]
\item[\ref{sec:approximate-convergence-criterion}] A convergence criterion for approximate trajectories.
\item[\ref{guide:appendix-5-1}] Proof of Theorem~\ref{thm:corank-one-convergence}.
\end{itemize}
\item \textbf{Appendix~\ref{sec:quantitative-localization}:} Quantitative bounds on the distance to the critical set.
\begin{itemize}[leftmargin=3.5em,itemsep=0pt,parsep=0pt,topsep=2pt]
\item[\ref{guide:appendix-6-1}] Uniform geometric bounds.
\item[\ref{guide:appendix-6-2}] Energy estimates.
\item[\ref{guide:appendix-6-3}] Proof of Theorem~\ref{thm:quantitative-localization}.
\item[\ref{guide:appendix-6-4}] Examples.
\end{itemize}
\item \textbf{Appendix~\ref{app:detailed-proofs}:} Proofs of auxiliary lemmas.
\begin{itemize}[leftmargin=3.5em,itemsep=0pt,parsep=0pt,topsep=2pt]
\item[\ref{guide:appendix-7-1}] Proof of Lemma~\ref{lem:finite-time-trap}.
\item[\ref{guide:appendix-7-2}] Proof of Lemma~\ref{lem:corrected-tracking}.
\item[\ref{guide:appendix-7-3}] Proof of Lemma~\ref{lem:compact-limit-geometry}.
\item[\ref{guide:appendix-7-4}] Proof of Lemma~\ref{lem:uniform-center-reduction}.
\item[\ref{guide:appendix-7-5}] Proof of Lemma~\ref{lem:approximate-trajectory-convergence}.
\item[\ref{guide:appendix-7-6}] Proof of Lemma~\ref{lem:quantitative-energy}.
\item[\ref{guide:appendix-7-7}] Proof of Lemma~\ref{lem:quantitative-band}.
\item[\ref{guide:appendix-7-8}] Proof of Lemma~\ref{lem:quantitative-excursion}.
\end{itemize}
\item \textbf{Appendix~\ref{app:counterexample-details}:} Construction, simulations, and proof of the nonconvergent example.
\begin{itemize}[leftmargin=3.5em,itemsep=0pt,parsep=0pt,topsep=2pt]
\item[\ref{guide:appendix-8-1}] Objective and second moments.
\item[\ref{guide:appendix-8-2}] Simulation setup.
\item[\ref{guide:appendix-8-3}] Proof of Example~\ref{ex:isolated-corank-two}.
\end{itemize}
\end{itemize}

\clearpage
\section{Comparison with existing continuous-time Adam and AdamW analyses}
\label{app:comparison}

We compare the long-time results of this paper with continuous-time Adam convergence, local AdamW stability, and optimization interpretations of adaptive methods.

\subsection{Scope of the comparison}\label{app:comparison-scope}
The comparisons below concern the full position--moment dynamics with
bias correction, fixed $\varepsilon>0$, and sufficiently small fixed
positive weight decay $\lambda$.

\subsection{Adam convergence and AdamW stability}\label{guide:appendix-1-1}
\label{app:adam-flows}

\citet{BB2021} derived a nonautonomous Adam system retaining both moment variables and bias correction. They established well-posedness and convergence results toward the critical set, with stronger conclusions under geometric assumptions such as a \L{}ojasiewicz inequality. \citet{BG2020} developed convergence results for a broader class of nonautonomous adaptive systems. For AdamW, fixed positive decay changes the equilibrium equation to
\[
\nabla F(p)+\lambda D(p)p=0,
\qquad
D(p)=\operatorname{diag}\bigl(\varepsilon+\sqrt{S_1(p)},\ldots,
\varepsilon+\sqrt{S_d(p)}\bigr).
\]
The original critical set therefore serves as reference geometry as $\lambda\downarrow0$, while the equilibria at positive decay satisfy this perturbed equation.

\citet{SuChenLiu2026Stability} study autonomous AdamW with exact
gradients, so $S_i=(\partial_iF)^2$, and replace the bias-correction
factors by their large-time limits. Their Theorem~3.1 identifies a
constrained and regularized implicit objective whose stationary points
are precisely the AdamW equilibrium positions. This first-order
equivalence does not ensure dynamical stability: even a strict local
minimum of that implicit objective may correspond to an unstable
AdamW equilibrium. Their Proposition~4.3 gives an exact one-dimensional
Hurwitz criterion via the cubic characteristic polynomial, and
Proposition~5.2 supplies a sufficient higher-dimensional condition
through a Lyapunov analysis of the linearized system.

Their Examples~4.5--4.6 establish instability at specified parameter
values and display nonconvergent oscillations numerically. Their
AdamMe modification preserves the AdamW equilibria but rescales the
weight-decay term. Under the assumptions of their Theorem~4.8,
including $b\le4a$ and initialization in the stated forward-invariant
domain, every bounded AdamMe trajectory approaches the equilibrium
set; isolated limiting equilibria give convergence to one equilibrium.
This positive result is for a modified optimizer and is conditional
on boundedness.

Our analysis follows the full bias-corrected trajectory from its
prescribed initial values and establishes compact trapping before
studying its limit set. The second-moment assumptions also differ.
Our long-time theory requires $S_i>0$, while the exact-gradient choice
$S_i=(\partial_iF)^2$ vanishes at critical points.
In Example~\ref{ex:isolated-corank-two}, the variation of $S$ along the
Hessian kernel generates the tangential motion responsible for
nonconvergence.

For discrete AdamW, \citet{zhou2024understanding} obtain stochastic-gradient complexity and convergence results under a decaying weight-decay schedule, using a dynamically regularized loss. Their finite-iteration estimates complement the fixed-decay limit-set analysis considered here.

\subsection{Implicit bias and normalized directions}\label{guide:appendix-1-2}
\label{app:implicit-bias}

\citet{WangMeng2021} and \citet{GronichVardi2026} study implicit bias in homogeneous classification problems and characterize limiting parameter directions through margin problems. In that setting, $x(t)/\|x(t)\|$ may converge while $\|x(t)\|$ diverges. Here we study bounded trajectories and the limiting behavior of the full state $(x,m,v)$. We also keep $\varepsilon>0$ fixed; it enters both the adaptive denominator and the equilibrium equation.

\citet{cattaneo2024implicit} use backward-error analysis to study implicit regularization induced by finite-step effects in Adam and RMSProp. Their modified equations describe discretization effects, whereas our results concern the long-time dynamics of the stated continuous-time system.

\subsection{Constrained interpretations and equilibria}\label{guide:appendix-1-3}
\label{app:adamw-constrained}

\citet{XieLi2024} show that, if discrete AdamW converges under their assumptions, its limit satisfies the KKT conditions of an $\ell_\infty$-constrained problem. For exact gradients and fixed $\varepsilon>0$, \citet[Theorem~3.1]{SuChenLiu2026Stability} identify the equilibrium equation with stationarity of an objective containing a separable barrier regularizer on $\|\lambda x\|_\infty<1$. Thus a positive stability constant is compatible with an optimization interpretation, with an additional regularization term.
For the general second-moment field $S$ considered here, the equilibrium equation involves $D(x)x$, whose components depend on $S_i(x)$. This vector field need not be a gradient, so the exact-gradient barrier representation does not extend directly. We work with the equilibrium equation itself and establish conditions under which the full trajectory converges to an equilibrium satisfying this equation.

\subsection{Time scaling}\label{guide:appendix-1-4}
\label{app:time-scaling}

The continuous-time limit depends on how the memory parameters scale with the step size. In our scaling, they approach one so that both moment relaxation rates remain finite, and the limiting system retains $(x,m,v)$.
\citet{GT2024} use continuous-time descriptions to study adaptive updates and normalization. \citet{DJK2025} keep the momentum parameters fixed as the step size vanishes, obtaining a fast--slow problem and an averaged limiting vector field. These different scalings can produce different equilibrium equations and long-time dynamics; the results here apply to the finite-rate moment equations specified in Section~\ref{sec:systems}.

\clearpage
\section{Well-posedness}\label{sec:wp}
	\label{sec:integral-formulation}\label{sec:wp-estimates}
	
We prove Theorem~\ref{thm:adamw-wp} by reducing the bias-corrected
initial-value problem to an integral equation for the position path.
Throughout this section, write
$\norm{h}_{[0,t]}=\sup_{0\leq r\leq t}\norm{h(r)}$.

\subsection{Integral formulation}\label{guide:appendix-3-1}
\label{sec:moment-integral-derivation}
We derive the moment integral formulas and their initial values along a
prescribed continuous position path, then obtain the integral equation
for position.
	
	For a continuous path $x:[0,T]\to\R^d$, define the moment functions for $t>0$ by
	\begingroup
\begin{equation}\label{eq:moment-integrals}
m_x(t)=\frac a{A(t)}\int_0^t e^{-a(t-r)}\nabla F(x(r))\,dr,\qquad
v_x(t)=\frac b{B(t)}\int_0^t e^{-b(t-r)}S(x(r))\,dr,
\end{equation}
\endgroup
	and put
	\begin{equation}\label{eq:normalized-path-moments}
		R_x(t)=\frac{m_x(t)}{\varepsilon+\sqrt{v_x(t)}}.
	\end{equation}
	For $t>0$, the identity
	$\int_0^t a e^{-a(t-r)}\,dr=A(t)$ and the triangle inequality give
	\begin{align*}
\norm{m_x(t)-\nabla F(x(0))} &=\norm{\int_0^t\frac{a e^{-a(t-r)}}{A(t)}
			\bigl[\nabla F(x(r))-\nabla F(x(0))\bigr]\,dr} \\
&\leq\int_0^t\frac{a e^{-a(t-r)}}{A(t)}
		\norm{\nabla F(x(r))-\nabla F(x(0))}\,dr \\
&\leq\sup_{0\leq r\leq t}\norm{\nabla F(x(r))-\nabla F(x(0))}.
\end{align*}
	Similarly,
	\[
	\norm{v_x(t)-S(x(0))}
	\leq\sup_{0\leq r\leq t}\norm{S(x(r))-S(x(0))}.
	\]
	By continuity, both suprema tend to zero as $t\downarrow0$.
Thus $m_x,v_x,R_x$ extend continuously to zero with values
	\begin{equation}\label{eq:entrance-moments}
		m_x(0)=\nabla F(x(0)),\qquad v_x(0)=S(x(0)),\qquad
		R_x(0)=\frac{\nabla F(x(0))}{\varepsilon+\sqrt{S(x(0))}}.
	\end{equation}

	The corresponding integral equation for position is
	\begin{equation}\label{eq:position-integral}
		x(t)=e^{-\lambda t}x_0-\int_0^t e^{-\lambda(t-r)}R_x(r)\,dr.
	\end{equation}
	
	\begin{proposition}[Integral formulation]
	\label{prop:integral-formulation}
	Fix a continuous path $x:[0,T]\to\mathbb R^d$ with $x(0)=x_0$.
	The formulas~\eqref{eq:moment-integrals}, extended to $t=0$ by
	\eqref{eq:entrance-moments}, give the unique continuous moments that
	satisfy the moment equations in~\eqref{eq:normalized} for $t>0$.
	With these moments, the position equation in~\eqref{eq:normalized}
	holds if and only if $x$ satisfies~\eqref{eq:position-integral}.
	Thus~\eqref{eq:moment-integrals} and~\eqref{eq:position-integral}
	form an integral formulation equivalent to~\eqref{eq:normalized}.
	\end{proposition}
	
	\begin{proof}
	Suppose first that $m$ satisfies the first moment equation in
\eqref{eq:normalized} for $t>0$ and extends continuously to zero.
The raw moment $\widetilde m(t):=A(t)m(t)$ then satisfies
$\widetilde m(t)\to0$ as $t\downarrow0$, consistent with zero
initialization of the raw filters.
	Because $A'(t)=a(1-A(t))$, the product rule gives
	\[
	\frac{d}{dt}\bigl(A(t)m(t)\bigr)
	=a\bigl(\nabla F(x(t))-A(t)m(t)\bigr),\qquad t>0.
	\]
	Multiplying by $e^{at}$ and integrating from $\tau>0$ to $t$ gives
	\begin{align*}
A(t)m(t)&=e^{-a(t-\tau)}A(\tau)m(\tau)+a\int_\tau^t e^{-a(t-r)}\nabla F(x(r))\,dr.
\end{align*}
	As $\tau\downarrow0$, the first term tends to zero because $A(\tau)$
	tends to zero and $m$ is continuous at zero. The integral tends to the same
	integral with lower endpoint zero. Dividing by $A(t)$ gives the first
	formula in~\eqref{eq:moment-integrals}.
Similarly, using $B'=b(1-B)$ and the integrating factor $e^{bt}$,
we obtain the second formula in~\eqref{eq:moment-integrals} by letting
$\tau\downarrow0$ and dividing by $B(t)$.

Conversely, differentiating the two formulas in~\eqref{eq:moment-integrals}
and using $A'=a(1-A)$ and $B'=b(1-B)$ gives the moment equations
in~\eqref{eq:normalized} for $t>0$. Together with the initial values
in~\eqref{eq:entrance-moments}, this proves existence and uniqueness
of the continuous moment paths in the proposition.
Equivalently, solving the raw equations
$\dot{\widetilde m}=a(\nabla F(x)-\widetilde m)$ and
$\dot{\widetilde v}=b(S(x)-\widetilde v)$ with zero initial data gives
$\widetilde m=A m_x$ and $\widetilde v=B v_x$.
Thus zero raw initialization determines the corrected initial values
in~\eqref{eq:entrance-moments}.

For the position equation, the same integrating-factor argument with
$e^{\lambda t}$ gives~\eqref{eq:position-integral}, using continuity
of $x$ and $R_x$ at zero.
Conversely, differentiating~\eqref{eq:position-integral} gives the
position equation for $t>0$. Continuity of $R_x$ at zero gives
$x\in C^1([0,T])$ and $x'(0)=-\lambda x_0-R_x(0)$.
	\end{proof}
	
	\subsection{Boundedness}\label{guide:appendix-3-2}
	
	Define $c(0)=1$ and, for $t>0$, set
\begin{equation}\label{eq:wp-c}
c(t)=\frac{a\sqrt{B(t)}}{\sqrt b\,A(t)}
		\left(\int_0^t e^{(b-2a)u}\,du\right)^{1/2}.
\end{equation}
	
	\begin{lemma}[Position update bound]\label{lem:wp-ratio}
		Assume $\nabla F$ and $S$ are continuous and satisfy
		\eqref{eq:moment-domination}. For every continuous path
		$x:[0,T]\to\R^d$,
		\begin{equation}\label{eq:wp-ratio-bound}
			|R_{x,i}(t)|\leq c(t),\qquad 0\leq t\leq T,\quad i=1,\ldots,d.
		\end{equation}
		The function $c$ is continuous on $[0,\infty)$.
	\end{lemma}
	
	\begin{proof}
		Cauchy--Schwarz with weights $e^{(b/2-a)(t-r)}$ and
		$\sqrt b\,e^{-b(t-r)/2}$, together with~\eqref{eq:moment-domination}, gives
		\[
|m_{x,i}(t)|^2\leq
		\frac{a^2B(t)}{bA(t)^2}
		\left(\int_0^t e^{(b-2a)u}\,du\right)v_{x,i}(t) =c(t)^2v_{x,i}(t).
\]
		If $v_{x,i}(t)>0$, divide the square root of this bound by
		$\varepsilon+\sqrt{v_{x,i}(t)}$; if $v_{x,i}(t)=0$, then
		$m_{x,i}(t)=0$. As $t\downarrow0$, $A(t)\sim at$, $B(t)\sim bt$, and the integral
		is asymptotic to $t$. Hence $c(t)\to1=c(0)$. At $t=0$, the bound
		follows from~\eqref{eq:moment-domination} and
		~\eqref{eq:entrance-moments}.
	\end{proof}
	
	Thus $\norm{R_x(t)}\leq\sqrt d\,c(t)$, and~\eqref{eq:position-integral} gives
	\begin{equation}\label{eq:wp-position}
		\norm{x(t)}\leq e^{-\lambda t}\norm{x_0}
		+\sqrt d\int_0^t e^{-\lambda(t-r)}c(r)\,dr.
	\end{equation}

	\subsection{Lipschitz continuity}\label{guide:appendix-3-3}
We next estimate how $R_x$ changes when the position path changes.
\begin{assumptionremark}[Regularity of the second moments]\label{rem:sqrt-regularity}
Under the sampling assumptions in Section~\ref{sec:systems}, the reverse
triangle inequality in $L^2$ gives
\begingroup
\[
\|\sqrt{S(x)}-\sqrt{S(y)}\|^2
\le\mathbb E\|g(x,\xi)-g(y,\xi)\|^2
\le L_K^2\|x-y\|^2,\qquad x,y\in K.
\]
\endgroup
Jensen's inequality gives the same Lipschitz bound for $\nabla F$.
Under Assumption~\ref{ass:long-time}(ii), local Lipschitz continuity of
$\sqrt S$ instead follows from positivity: on a compact ball $K$, set
$s_K=\min_{x\in K,i}S_i(x)>0$ and use
\[
|\sqrt{S_i(x)}-\sqrt{S_i(y)}|
=\frac{|S_i(x)-S_i(y)|}{\sqrt{S_i(x)}+\sqrt{S_i(y)}}
\le\frac{|S_i(x)-S_i(y)|}{2\sqrt{s_K}}.
\]
\end{assumptionremark}

\begin{lemma}[Path Lipschitz bound]
		\label{lem:wp-lipschitz}
		Assume $\nabla F$ and $\sqrt S$ are locally Lipschitz.
		For every compact ball $K\subset\R^d$, there is a constant $H_K>0$
		such that, for every $T>0$ and continuous paths $x,y:[0,T]\to K$,
		\begin{equation}\label{eq:wp-lipschitz}
			\norm{R_x(t)-R_y(t)}\leq H_K\norm{x-y}_{[0,t]},
			\qquad 0\leq t\leq T.
		\end{equation}
	\end{lemma}
	
	\begin{proof}
		Since $K$ is a compact ball, the assumed local Lipschitz continuity of
		$\nabla F$ and $\sqrt S$ gives a common Lipschitz constant $L_K$
		for these two maps on $K$. Set
		$M_K=\sup_K\norm{\nabla F}<\infty$.
		For $t>0$, we have
		\begin{align*}
\norm{m_x(t)-m_y(t)} &=\left\|\frac{a}{A(t)}\int_0^t e^{-a(t-r)}
		\bigl(\nabla F(x(r))-\nabla F(y(r))\bigr)\,dr\right\| \\
&\leq\frac{a}{A(t)}\int_0^t e^{-a(t-r)}
		\norm{\nabla F(x(r))-\nabla F(y(r))}\,dr \\
&\leq L_K\norm{x-y}_{[0,t]}
		\frac{a}{A(t)}\int_0^t e^{-a(t-r)}\,dr =L_K\norm{x-y}_{[0,t]},
\end{align*}
		\begin{align*}
\norm{\sqrt{v_x(t)}-\sqrt{v_y(t)}}^2 &=\sum_{i=1}^d\left|
		\left(\frac{b}{B(t)}\int_0^t e^{-b(t-r)}S_i(x(r))\,dr\right)^{1/2}
		-\left(\frac{b}{B(t)}\int_0^t e^{-b(t-r)}S_i(y(r))\,dr\right)^{1/2}
		\right|^2 \\
&\leq\frac{b}{B(t)}\int_0^t e^{-b(t-r)}
		\norm{\sqrt{S(x(r))}-\sqrt{S(y(r))}}^2\,dr \\
&\leq L_K^2\norm{x-y}_{[0,t]}^2
		\frac{b}{B(t)}\int_0^t e^{-b(t-r)}\,dr =L_K^2\norm{x-y}_{[0,t]}^2.
\end{align*}
		Since $y(r)\in K$ and the first-moment kernel has total mass one,
		for each coordinate
		\[
|m_{y,i}(t)|\leq\frac{a}{A(t)}\int_0^t e^{-a(t-r)}
		|\partial_iF(y(r))|\,dr \leq M_K\frac{a}{A(t)}\int_0^t e^{-a(t-r)}\,dr =M_K.
\]
		For each coordinate, subtracting the ratios and using
		$\varepsilon+\sqrt{v_{x,i}(t)}\geq\varepsilon$ gives
		\begin{gather*}
R_{x,i}(t)-R_{y,i}(t) =\frac{m_{x,i}(t)-m_{y,i}(t)}{\varepsilon+\sqrt{v_{x,i}(t)}}
		+\frac{m_{y,i}(t)\bigl(\sqrt{v_{y,i}(t)}-\sqrt{v_{x,i}(t)}\bigr)}
		{\bigl(\varepsilon+\sqrt{v_{x,i}(t)}\bigr)
		\bigl(\varepsilon+\sqrt{v_{y,i}(t)}\bigr)}, \\
\begin{aligned}
\norm{R_x(t)-R_y(t)} &\leq\frac{\norm{m_x(t)-m_y(t)}}{\varepsilon}
		+\frac{M_K}{\varepsilon^2}
		\norm{\sqrt{v_x(t)}-\sqrt{v_y(t)}} \\
&\leq L_K\left(\frac{1}{\varepsilon}
		+\frac{M_K}{\varepsilon^2}\right)\norm{x-y}_{[0,t]} =:H_K\norm{x-y}_{[0,t]}.
\end{aligned}
\end{gather*}
The integral triangle inequality gives the square-root estimate even
when a second moment vanishes, and continuity gives the estimate at $t=0$.
	\end{proof}

	\subsection{Proof of Theorem~\ref{thm:adamw-wp}}\label{guide:appendix-3-4}
	\label{sec:wp-proof}

\begin{restatement}{Theorem}{thm:adamw-wp}
If $\nabla F$ and $\sqrt S:\R^d\to[0,\infty)^d$ are locally Lipschitz,
		then for every $x_0\in\R^d$, system~\eqref{eq:normalized} has a unique
		global solution $z=(x,m,v)$ with initial value~\eqref{eq:normalized-initial-value} and $v\ge0$.
		The solution is $C^1$ on $(0,\infty)$, and $x\in C^1([0,\infty))$.
\end{restatement}

	By Lemmas~\ref{lem:wp-ratio} and~\ref{lem:wp-lipschitz}, the position
map preserves a common ball and is Lipschitz on its path space.
An exponentially weighted norm will make this map a contraction on
any fixed bounded interval. Once the position is determined, the
moments follow from~\eqref{eq:moment-integrals}.
	
\begin{proof}[Proof of Theorem~\ref{thm:adamw-wp}]
		Fix any $T>0$ and set
		\begin{equation*}
R_T=\norm{x_0}+\sqrt d\int_0^T c(r)\,dr,\qquad K_T=\overline B(0,R_T),\qquad\text{\normalfont and}\qquad
\mathcal X_T=\{x\in C([0,T];K_T):x(0)=x_0\}.
\end{equation*}
		Let $\Phi x$ be the right-hand side of~\eqref{eq:position-integral}.
		Lemma~\ref{lem:wp-ratio} applies to every candidate path, so
		\[
\norm{\Phi x(t)}\leq e^{-\lambda t}\norm{x_0}
		+\sqrt d\int_0^t e^{-\lambda(t-r)}c(r)\,dr \leq R_T.
\]
		Thus $\Phi$ maps $\mathcal X_T$ into itself. Equip this nonempty closed
		path space with the norm
		\[
		\norm{h}_{\kappa,T}=\sup_{0\leq t\leq T}e^{-\kappa t}\norm{h(t)},
		\qquad \kappa>H_{K_T}.
		\]
		It is complete because this norm is equivalent to the uniform norm.
		Since $\norm{x-y}_{[0,r]}\leq e^{\kappa r}\norm{x-y}_{\kappa,T}$,
		Lemma~\ref{lem:wp-lipschitz} gives
		\[
		\begin{aligned}
e^{-\kappa t}\norm{\Phi x(t)-\Phi y(t)}
&\leq H_{K_T}e^{-\kappa t}
\int_0^t e^{-\lambda(t-r)}\norm{x-y}_{[0,r]}\,dr\\
&\leq H_{K_T}\norm{x-y}_{\kappa,T}
\int_0^t e^{-(\lambda+\kappa)(t-r)}\,dr\\
&=H_{K_T}\norm{x-y}_{\kappa,T}
\frac{1-e^{-(\lambda+\kappa)t}}{\lambda+\kappa}\\
&\leq \frac{H_{K_T}}{\lambda+\kappa}\norm{x-y}_{\kappa,T}
\leq\frac{H_{K_T}}{\kappa}\norm{x-y}_{\kappa,T}.
\end{aligned}
		\]
		Because $H_{K_T}/\kappa<1$, this is a contraction estimate.
		Banach's theorem therefore gives a unique fixed point on $[0,T]$.
		Every solution lies in $\mathcal X_T$ by~\eqref{eq:wp-position},
		so this proves uniqueness among all solutions. Consequently the
		solutions for different $T$ agree on overlapping intervals and define
		a global solution. By~\eqref{eq:moment-integrals}, both corrected moments are continuous
at zero and $C^1$ for $t>0$. The position equation gives $x\in C^1([0,\infty))$.
	\end{proof}
	
	\clearpage
	\section{Localization and invariant manifolds}\label{sec:invariant-manifolds}

We prove localization near the zero-decay equilibrium set and then
combine it with the local invariant manifolds to prove
Theorems~\ref{thm:full-state-localization} and~\ref{thm:invariant-manifold}.

\subsection{Assumptions and their scope}\label{sec:long-time-assumption-scope}
We recall the assumptions used in the long-time analysis.

\begin{restatement}{Assumption}{ass:long-time}
\leavevmode
\begin{enumerate}[label=(\roman*),leftmargin=*,itemsep=0pt,parsep=0pt,topsep=2pt]
\item The memory rates satisfy $0<b\le4a$.
\item The objective $F\in C^1(\R^d)$ is coercive, its gradient
$\nabla F$ is locally Lipschitz, and the second-moment map
$S:\R^d\to(0,\infty)^d$ is locally Lipschitz.
\item On an open neighborhood of $\operatorname{Crit}(F)$, the objective
$F$ is twice continuously differentiable and $S$ is continuously
differentiable.
\end{enumerate}
\end{restatement}

The following remark distinguishes the assumptions needed for
well-posedness, compact trapping, and the stronger convergence conclusions.

\par\addvspace{0.5\baselineskip}
\begin{assumptionremark}[Roles of the long-time assumptions]\label{rem:long-time-assumptions}
The localization results, Theorems~\ref{thm:full-state-localization}
	and~\ref{thm:quantitative-localization}, use Assumption~\ref{ass:long-time}(i)--(ii).
	The invariant-manifold and point-convergence results,
	Theorems~\ref{thm:invariant-manifold} and~\ref{thm:corank-one-convergence},
	also use Assumption~\ref{ass:long-time}(iii).
In Assumption~\ref{ass:long-time}, condition~(i) is the memory-rate
	inequality that supplies dissipation. Condition~(ii) makes $F$ coercive,
	which bounds its sublevels, and requires $S$ to be strictly positive;
	this positivity, together with the local Lipschitz regularity of $S$,
	avoids the square-root singularity. Condition~(iii) provides the
	near-critical differentiability needed for linearization.
	Condition~(i) is sufficient for the dissipation argument used here;
	we do not claim that it is necessary for convergence. Coercivity is
	a global restriction on the objective and does not follow from adding
	weight decay to the update. Strict positivity of $S_i$ is pointwise:
	no uniform positive lower bound on all of $\R^d$ is required, although
	continuity gives one on each compact set. It excludes a coordinate
	whose sampled gradient vanishes almost surely at some position.
	The additional smoothness in~(iii) is required only near the critical
	set, rather than globally.
\end{assumptionremark}

\begin{restatement}{Assumption}{ass:critical-manifold}
The connected components $\mathcal M$ of $\operatorname{Crit}(F)$ are
	$C^1$ embedded manifolds without boundary. Every compact set meets
	only finitely many components, and the Hessian nullity is constant on each:
	\begin{equation*}
		\dim\ker\nabla^2F(q)=k_{\mathcal M}\qquad(q\in\mathcal M).
	\end{equation*}
	The integers $k_{\mathcal M}$ may vary between components.
\end{restatement}

\par\addvspace{0.5\baselineskip}
\begin{assumptionremark}[Scope of the critical-manifold hypothesis]\label{rem:critical-manifold-assumption}
Two classes of examples motivate Assumption~\ref{ass:critical-manifold}.
	First, if all critical points are isolated, the assumption
	holds automatically: the closed critical set is finite on every compact
	set, its components are singletons, and constant nullity on each component
	imposes no restriction on the Hessian.
	Second, Morse--Bott critical manifolds motivate the assumption for
	nonisolated critical points. For a critical component $\mathcal M$,
	differentiating $\nabla F|_{\mathcal M}=0$ gives
	$T_q\mathcal M\subseteq\ker\nabla^2F(q)$; the Morse--Bott condition
	requires equality. Under Assumption~\ref{ass:critical-manifold}, the Hessian kernel has
	constant dimension along each component and may contain directions
	normal to $\mathcal M$ in addition to the tangent directions.
	Concrete examples of Morse--Bott minima arise locally in nonlinear
	least-squares problems with zero residual and a locally constant-rank
	residual Jacobian, including overparameterized models; see
	Section~1.2 of~\cite{RB2025}. These local examples illustrate the
	critical-manifold geometry; the global assumptions must still be
	checked. The hypothesis also excludes singular critical components,
	infinitely many components accumulating in a compact region, and
	changes of Hessian nullity along a component.
\end{assumptionremark}

\subsection{Energy estimates}\label{guide:appendix-4-1}
We study the trajectory through its limit set. If the trajectory
is contained in a compact subset of $U$, it approaches its limit set
in distance. Approximation by the autonomous flow will allow us to apply
the localization argument
in Appendix~\ref{app:compact-invariant}.
	
We use the $\lambda=0$ Adam energy estimate of \citet{BB2021},
written in bias-corrected coordinates, to construct a compact region
that the zero-decay flow maps strictly into its interior after a fixed
time. Positive decay adds $-\lambda x\cdot\nabla F(x)$ to the energy
derivative, a term that may have either sign. The interior margin gives
a common trapping region for sufficiently small positive decay.
	
	\begin{lemma}[Zero-decay energy]
		\label{lem:hp-zero-decay-energy}
		Under Assumption~\ref{ass:long-time}(i)--(ii) of
		Section~\ref{sec:long-time-assumptions}, let $z^0$ be the zero-decay
		solution of~\eqref{eq:normalized} given by Theorem~\ref{thm:adamw-wp}.
		For this estimate write $z^0=(x,m,v)$, suppressing the zero-decay
		superscript, so that $(x,m,v)(0)=(x_0,\nabla F(x_0),S(x_0))$. For $d_i(t)=\varepsilon+\sqrt{v_i(t)}$, define
		\[
		\mathcal H(t)=F(x(t))+
		\frac{A(t)}{2a}\sum_i\frac{m_i(t)^2}{d_i(t)}
		\quad\text{\normalfont and}\quad \mathcal H(0)=F(x_0).
		\]
		Then $\mathcal H$ is nonincreasing and
		$F(x(t))\le F(x_0)$ for all $t\ge0$. The normalized path
		$(x,m,v)$ has compact closure $K_0\subset U$.
	\end{lemma}
	
	\begin{proof}
		Theorem~\ref{thm:adamw-wp} applies because positive locally Lipschitz
		$S_i$ have locally Lipschitz square roots.
		The normalized moment equations are
		$\dot m=(a/A)(\nabla F-m)$ and $\dot v=(b/B)(S-v)$; their continuous
		extensions at zero are $\nabla F(x_0)$ and $S(x_0)$. For $t>0$,
		$\dot F=-\sum_i(\partial_i F)m_i/d_i$ and
		\[
		\frac d{dt}\frac{A m_i^2}{2ad_i}
		=\frac{A'm_i^2}{2ad_i}+\frac{m_i(\partial_i F-m_i)}{d_i}
		-\frac{Abm_i^2(S_i-v_i)}{4aB\sqrt{v_i}d_i^2}.
		\]
		Hence
		\begin{align*}
\dot{\mathcal H}
&=-\sum_i\frac{(\partial_iF)m_i}{d_i}
 +\sum_i\left[
 \frac{(\partial_iF)m_i-m_i^2}{d_i}
 +\frac{A'm_i^2}{2ad_i}
 -\frac{Abm_i^2(S_i-v_i)}{4aB\sqrt{v_i}d_i^2}
 \right] \\
&=\sum_i\left(\frac{A'}{2a}-1\right)\frac{m_i^2}{d_i}
 +\frac{bA}{4aB}\sum_i\frac{m_i^2(v_i-S_i)}{\sqrt{v_i}d_i^2}.
\end{align*}
		Since $A'=ae^{-at}$, this becomes
		\begingroup
\[
\dot{\mathcal H}
		=-c_t\sum_i\frac{m_i^2}{d_i}
		+d_t\sum_i\frac{m_i^2(v_i-S_i(x))}
		{\sqrt{v_i}d_i^2},\qquad
c_t=1-\tfrac12 e^{-at}\quad\text{\normalfont and}\quad
		d_t=\frac{bA(t)}{4aB(t)}.
\]
\endgroup
        By Lemma~\ref{lem:memory-rate-coefficients} below,
        \(
        d_t\le c_t.
        \)
        Using $S_i>0$, $d_i=\varepsilon+\sqrt{v_i}$ and $c_t\ge1/2$, we obtain
        \[
        \begin{aligned}
        \dot{\mathcal H}
        &\le-c_t\sum_i\frac{m_i^2}{d_i}
          +d_t\sum_i\frac{m_i^2\sqrt{v_i}}{d_i^2}\\
        &=-\sum_i\frac{m_i^2}{d_i^2}
          \bigl(c_td_i-d_t\sqrt{v_i}\bigr)\\
        &=-\sum_i\frac{m_i^2}{d_i^2}
          \bigl[c_t\varepsilon+(c_t-d_t)\sqrt{v_i}\bigr]
        \le-\frac\varepsilon2\sum_i\frac{m_i^2}{d_i^2}\le0.
        \end{aligned}
        \]
        Since $\mathcal H(t)\to F(x_0)$ as $t\downarrow0$,
        \[
F(x(t))\le\mathcal H(t)\le F(x_0),\qquad t\ge0.
\]
        By coercivity, $x(t)$ stays in the compact set
        $X_0=\{x\in\mathbb R^d:F(x)\le F(x_0)\}$.
        The kernels in~\eqref{eq:moment-integrals} have mass one, so
        \[
\|m(t)\|\le\max_{X_0}\|\nabla F\|
\qquad\text{\normalfont and}\qquad
0<\min_{X_0}S_i\le v_i(t)\le\max_{X_0}S_i<\infty.
\]
        Thus $(x,m,v)$ has compact closure in $U$.
	\end{proof}
	
    \begin{lemma}[Coefficient bound]
    \label{lem:memory-rate-coefficients}
    Let $0<b\le4a$. For $t>0$, define
    \[
    c_t:=1-\frac12e^{-at}\quad\text{\normalfont and}\quad
    d_t:=\frac{b(1-e^{-at})}{4a(1-e^{-bt})}.
    \]
    Then $d_t\le c_t$.
    \end{lemma}

    \begin{proof}
    For fixed $t>0$ and $r>0$,
    \[
\frac{d}{dr}\frac{r}{1-e^{-rt}}
    =\frac{1-(1+rt)e^{-rt}}{(1-e^{-rt})^2} >0,
\]
    since $e^{rt}>1+rt$. Thus $b\le4a$ gives, with $u=e^{-at}\in(0,1)$,
    \[
d_t=\frac{1-e^{-at}}{4a}\frac{b}{1-e^{-bt}} \le\frac{1-e^{-at}}{1-e^{-4at}} =\frac1{1+u+u^2+u^3} \le1-\frac u2 =c_t.
\]
The last inequality follows by multiplying by the positive denominator
$1+u+u^2+u^3$ and using the exact identity
\[
\left(1-\frac u2\right)(1+u+u^2+u^3)-1
=\frac u2\bigl(1+u+u^2(1-u)\bigr)\ge0.
\]
    \end{proof}

    Lemma~\ref{lem:hp-zero-decay-energy} concerns the bias-corrected zero-decay trajectory.
    To construct a fixed trapping region, we now use the autonomous
    zero-decay system obtained by replacing $A(t),B(t)$ with their limits $1$:
    \[
\dot x=-\frac{m}{\varepsilon+\sqrt v},\qquad
\dot m=a(\nabla F(x)-m),\qquad\text{\normalfont and}\qquad
\dot v=b(S(x)-v).
\]
    Replacing $A(t)$ by $1$ in the expression for $\mathcal H$ gives
    the corresponding energy function
	\begin{equation}\label{eq:adam-reference-energy}
		V(x,m,v)=F(x)+\frac1{2a}\sum_i\frac{m_i^2}{\varepsilon+\sqrt{v_i}}.
	\end{equation}
    Along this autonomous system, the same calculation with $A=B=1$
    gives, writing $d_i=\varepsilon+\sqrt{v_i}$ and using $S_i>0$ and $b\le4a$,
    \begingroup
\begin{equation}\label{eq:adam-energy-dissipation}
\dot V
    =-\sum_i\frac{m_i^2}{d_i}
      +\frac b{4a}\sum_i\frac{m_i^2(v_i-S_i)}{\sqrt{v_i}d_i^2}
\le-\sum_i\frac{m_i^2}{d_i^2}
      \left[\varepsilon+\left(1-\frac b{4a}\right)\sqrt{v_i}\right]
      \le-\varepsilon\sum_i\frac{m_i^2}{d_i^2}.
\end{equation}
\endgroup

\subsection{Proof of Theorem~\ref{thm:full-state-localization}}\label{guide:appendix-4-2}
	\label{sec:finite-time-trapping}
	\label{sec:compact-limit-geometry}

\begin{restatement}{Theorem}{thm:full-state-localization}
Under Assumption~\ref{ass:long-time}(i)--(ii), if $F(\operatorname{Crit}(F))$ has empty
		interior, then
		\begin{equation*}
			\lim_{\lambda\downarrow0}\limsup_{t\to\infty}
			\operatorname{dist}(z^\lambda(t),\mathcal E_0)=0.
		\end{equation*}
\end{restatement}

    We next use the zero-decay energy estimate to obtain a common compact
    region for trajectories with sufficiently small weight decay.
	\begin{lemma}[Small-decay trapping]
\label{lem:finite-time-trap}
Suppose Assumption~\ref{ass:long-time}(i)--(ii) holds and that
$F(\operatorname{Crit}(F))$ has empty interior. There are compact
sets $B,K\subset U$ with $B\subset\operatorname{int}K$ such that:
\leavevmode
\begin{enumerate}[label=(\alph*),leftmargin=*,itemsep=0pt,parsep=0pt,topsep=2pt]
\item The set $B$ is positively invariant under $f_0$, and
\begin{equation}\label{eq:finite-time-inward}
 \Phi_0(t,B)\Subset\operatorname{int}B\qquad(t>0).
\end{equation}
\item For every sufficiently large $T$, there is
$\lambda_0>0$ such that
\[
 \inf_{0\le\lambda\le\lambda_0}
       \operatorname{dist}(z^\lambda(T),B^c)>0,
\]
and, for $t\ge T$ and $0\le\lambda\le\lambda_0$,
\[
 z^\lambda(t)\in K.
\]
Moreover, $\{z^\lambda(t):t\ge0,\ 0\le\lambda\le\lambda_0\}$
has compact closure in $U$.
\end{enumerate}
The sets $B,K$ depend only on $F,S,a,b,\varepsilon$ and $x_0$;
the threshold $\lambda_0$ may also depend on $T$.
\end{lemma}
\begin{proof}
See Appendix~\ref{app:trapping-finite-length}.
\end{proof}

\begin{figure}[htbp]
\centering
\includegraphics[width=.5\linewidth]{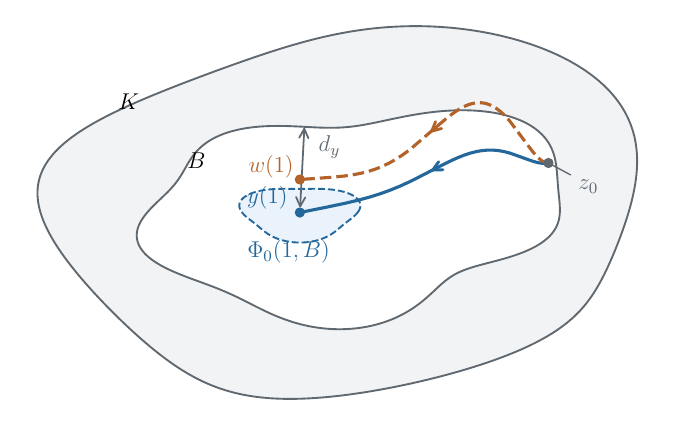}
\caption{One step of the trapping argument (Lemma~\ref{lem:unit-interval-trap}). Both paths start at
$z_0\in B$. The reference path $y(t)=\Phi_0(t,z_0)$ stays in $B$, while
the perturbed path $w$ may leave $B$ and stays in $K$ for $0\le t\le1$.
The double arrow marks $d_y=\operatorname{dist}(y(1),B^c)$.
The bound $\|w(1)-y(1)\|<d_y$ puts $w(1)$ in $\operatorname{int}B$,
so the comparison can start again there. The paths are schematic.}
\label{fig:trapping-geometry}
\end{figure}

For sufficiently small weight decay, the bias-corrected AdamW trajectory
eventually remains in one compact region, uniformly over the decay
parameter. Part~(b) uses the comparison in the next lemma.
Figure~\ref{fig:trapping-geometry} illustrates this comparison over one
unit of time.

\begin{lemma}[Confinement under small perturbations]
\label{lem:unit-interval-trap}
Let $B,K\subset U$ be nonempty compact sets with
$B\subset\operatorname{int}K$. Suppose the zero-decay flow satisfies
$\Phi_0([0,1],B)\subset B$ and
$\Phi_0(1,B)\Subset\operatorname{int}B$.
There is $\kappa>0$ such that every locally absolutely continuous
path $w:[0,\infty)\to U$ with $w(0)\in B$ and
$\|\dot w-f_0(w)\|\le\kappa$ almost everywhere while in $K$ satisfies
$w(t)\in K$ for all $t\ge0$ and $w(j)\in B$ for all integers $j\ge0$.
\end{lemma}

\begin{proof}
Set $d=\min\{\operatorname{dist}(B,K^c),
\operatorname{dist}(\Phi_0(1,B),B^c)\}>0$.
Choose a Lipschitz constant $L>0$ for $f_0$ on $K$ and take
$\kappa e^L<d$. With $y(t)=\Phi_0(t,w(0))$, Gronwall's inequality
gives, up to the first exit from $\operatorname{int}K$,
\[
\|w(t)-y(t)\|
\le\kappa\int_0^t e^{L(t-r)}\,dr\le\kappa e^L<d,
\qquad 0\le t\le1.
\]
Since $y(t)\in B$, even at a first exit time the bound holds by
continuity and puts $w(t)$ in $\operatorname{int}K$. Thus there is
no exit on $[0,1]$. At $t=1$, the same bound puts $w(1)$ in $B$
because $y(1)\in\Phi_0(1,B)$. Repeating on each unit interval gives
both conclusions.
\end{proof}

We now locate the limiting motion of trajectories that approach the autonomous flow.
A continuous path $w:[t_0,\infty)\to U$ tracks the autonomous flow
	on late time intervals if, for every $H>0$, that flow is defined throughout
	the interval for all sufficiently large $t$ and
	\begin{equation}\label{eq:deterministic-flow-tracking}
	 \sup_{0\le s\le H}
	 \|w(t+s)-\Phi_\lambda(s,w(t))\|\longrightarrow0
	 \qquad(t\to\infty). 
	\end{equation}
	Such a path is called an asymptotic pseudotrajectory.

\begin{figure}[H]
	\centering
	\includegraphics[width=.50\linewidth]{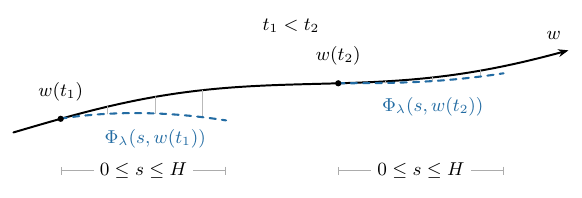}
	\caption{An asymptotic pseudotrajectory $w$ and autonomous solutions starting
		at $w(t_1)$ and $w(t_2)$, with $t_1<t_2$. Both comparisons use the same
		elapsed time interval $0\le s\le H$. Gray segments join points at equal $s$;
		the later comparison illustrates the decreasing tracking error. The curves are schematic.}
	\label{fig:apt-windows}
\end{figure}

The corrected trajectory satisfies this condition.
\begin{lemma}[Flow tracking]
\label{lem:corrected-tracking}
Under the hypotheses of Lemma~\ref{lem:finite-time-trap}, for every
sufficiently small $\lambda\ge0$, the corrected trajectory $z^\lambda$
satisfies \eqref{eq:deterministic-flow-tracking}.
\end{lemma}
\begin{proof}
See Appendix~\ref{app:corrected-tracking}.
\end{proof}

For a continuous path $w:[t_0,\infty)\to U$, write
\[
\omega(w):=\{q\in U:w(t_n)\to q\,\text{\normalfont  for some sequence }t_n\to\infty\}.
\]
For fixed $\lambda\ge0$, if the range of $w$ has compact closure in $U$
and $w$ satisfies
\eqref{eq:deterministic-flow-tracking}, then $L:=\omega(w)$ is nonempty,
compact and connected. For every $z\in L$, the solution through $z$
is defined for all real times and satisfies
\begin{equation}\label{eq:limit-set-invariance}
\Phi_\lambda(\mathbb R,z)\subset L.
\end{equation}
These facts are proved at the beginning of
Appendix~\ref{app:compact-invariant}. The following lemma bounds the
distance of these limit sets from $\mathcal E_0$ as $\lambda\downarrow0$.

	\begin{lemma}[Uniform localization]
		\label{lem:compact-limit-geometry}
Under the hypotheses of Theorem~\ref{thm:full-state-localization},
fix a nonempty compact set $K\subset U$. There exists
$\rho_K(\lambda)\to0$ as $\lambda\downarrow0$ such that, for every
sufficiently small $\lambda>0$ and every continuous path
$w:[0,\infty)\to U$ whose range has compact closure in $U$
and which satisfies \eqref{eq:deterministic-flow-tracking}
with $\omega(w)\subset K$,
\begin{equation}\label{eq:compact-uniform-localization}
\limsup_{t\to\infty}\operatorname{dist}(w(t),\mathcal E_0)
\le\rho_K(\lambda).
\end{equation}
The function $\rho_K$ and the decay threshold are independent of $w$.
	\end{lemma}
\begin{proof}
See Appendix~\ref{app:compact-invariant}.
\end{proof}
We now apply this lemma to the corrected trajectory.

\begin{proof}[Proof of Theorem~\ref{thm:full-state-localization}]
By Lemma~\ref{lem:finite-time-trap}(b), there are a compact set
$K\subset U$ and $\lambda_0>0$ such that, for every
$0\le\lambda\le\lambda_0$, the trajectory $z^\lambda$ has compact
closure in $U$ and belongs to $K$ for all sufficiently large times.
In particular, $\omega(z^\lambda)\subset K$.
After decreasing $\lambda_0$ if necessary, Lemma~\ref{lem:corrected-tracking}
gives, for every $H>0$,
\[
\sup_{0\le s\le H}
\|z^\lambda(t+s)-\Phi_\lambda(s,z^\lambda(t))\|
\longrightarrow0\qquad(t\to\infty).
\]
Applying Lemma~\ref{lem:compact-limit-geometry} with $w=z^\lambda$,
we obtain, for all sufficiently small $\lambda>0$,
\[
\limsup_{t\to\infty}\operatorname{dist}(z^\lambda(t),\mathcal E_0)
\le\rho_K(\lambda)\longrightarrow0\qquad(\lambda\downarrow0).
\]
This proves~\eqref{eq:full-state-localization}.

\end{proof}

	\subsection{Local invariant manifolds}\label{guide:appendix-4-3}
	\label{sec:invariant-manifold-statement}
We first check the critical-value condition under the hypotheses of
Theorem~\ref{thm:invariant-manifold}. Any two points of a connected
critical manifold can be joined by piecewise $C^1$ paths. Along such a
path $q(s)$,
\[
\frac{d}{ds}F(q(s))=\nabla F(q(s))\cdot q'(s) =0,
\]
so $F$ is constant on the component. The component is closed and lies
in a compact sublevel set by coercivity, hence is compact. Local
finiteness implies that each compact ball meets only finitely many
components. Covering the space by countably many balls therefore gives at most
countably many components and critical values. In particular, the
critical-value hypothesis of Theorem~\ref{thm:full-state-localization}
holds.
	
The next lemma shows that every complete orbit staying in a suitable
neighborhood of a critical component lies on a manifold whose dimension
equals the number of degenerate directions of the objective. For each critical
component $\mathcal M$, define the corresponding set of zero-decay
equilibria by
\begin{equation}\label{eq:lifted-critical-component}
\mathcal N_{\mathcal M}:=\{(q,0,S(q)):q\in\mathcal M\}.
\end{equation}

	\begin{lemma}[Uniform local reduction]
		\label{lem:uniform-center-reduction}
		Let $a,b,\varepsilon>0$ and let $\mathcal M$ be a nonempty compact
		connected $C^1$ embedded manifold without boundary, consisting of
		critical points of $F$. Assume $F\in C^2$ and $S\in C^1$, with
		$S_i>0$, near $\mathcal M$, and $\dim\ker\nabla^2F(q)=k$ on
		$\mathcal M$. Let $\mathcal N_{\mathcal M}$ be as in
		\eqref{eq:lifted-critical-component}.
		There are an open neighborhood $V$
		of $\mathcal N_{\mathcal M}$ with compact closure in $U$, a number
		$\lambda_{\mathcal M}>0$, and $k$-dimensional $C^1$ embedded local
		manifolds $\mathcal I_{\mathcal M,\lambda}\subset U$ for
        $|\lambda|<\lambda_{\mathcal M}$, locally invariant for $f_\lambda$,
        such that every solution $z:\mathbb R\to V$ of
        $\dot z=f_\lambda(z)$ satisfies
        $z(\mathbb R)\subset\mathcal I_{\mathcal M,\lambda}$.
	\end{lemma}
\begin{proof}
See Appendix~\ref{app:center-reduction}.
\end{proof}

	\subsection{Proof of Theorem~\ref{thm:invariant-manifold}}\label{guide:appendix-4-4}
	\label{sec:manifold-limit-sets}
	\label{sec:full-state-localization}

\begin{restatement}{Theorem}{thm:invariant-manifold}
Under Assumption~\ref{ass:long-time}(i)--(iii) and
		Assumption~\ref{ass:critical-manifold}, there exists
		$\lambda_{\mathrm{man}}>0$ such that, for every $0<\lambda<\lambda_{\mathrm{man}}$, there are
		a connected component $\mathcal M$ of $\operatorname{Crit}(F)$ and a
		$k_{\mathcal M}$-dimensional $C^1$ locally invariant manifold
		$\mathcal I_\lambda$ for~\eqref{eq:auto} satisfying
		\begin{equation*}
			\operatorname{dist}\bigl(z^\lambda(t),\mathcal I_\lambda\bigr)
			\longrightarrow0\qquad(t\to\infty).
		\end{equation*}
\end{restatement}

Each localized limit set lies on one of the invariant manifolds
constructed above.

	\begin{lemma}[Limit sets on invariant manifolds]
	\label{lem:manifold-limit-sets}
	Let $K\subset U$ be nonempty and compact.
		Under the full hypotheses of Theorem~\ref{thm:invariant-manifold},
		there are a finite collection
		$\mathfrak M_K$ of critical components and a common $\lambda_K>0$
		such that, for every $0<\lambda<\lambda_K$ and every continuous
path $w:[0,\infty)\to U$ whose range has compact closure in $U$
and which satisfies \eqref{eq:deterministic-flow-tracking}
with $\omega(w)\subset K$,
its limit set $\omega(w)\subset K$ is contained in one $\mathcal I_{\mathcal M,\lambda}$,
$\mathcal M\in\mathfrak M_K$.
	\end{lemma}

\begin{figure}[htbp]
\centering
\includegraphics[width=.68\linewidth]{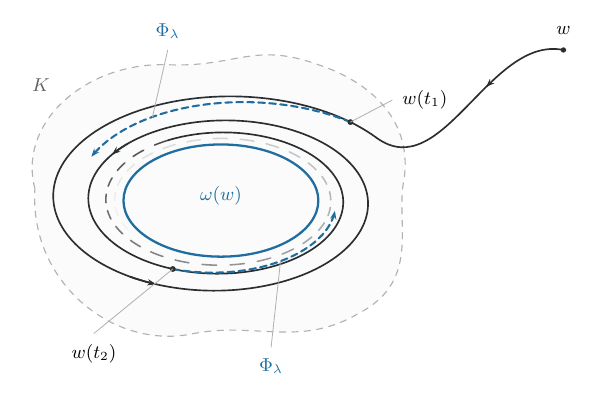}
\caption{Schematic illustration of Lemma~\ref{lem:manifold-limit-sets}.
The path $w$ enters the compact region $K$ and winds toward the invariant ellipse;
its fading dashed tail indicates continued motion.
Blue dashed arcs show the autonomous solutions starting at $w(t_1)$ and $w(t_2)$
with $t_1<t_2$, over the same interval $0\le s\le H$; the later comparison is closer to $w$.
The solid blue ellipse represents $\omega(w)=\mathcal I_{\mathcal M,\lambda}$
in this illustration. The lemma requires only
$\omega(w)\subset\mathcal I_{\mathcal M,\lambda}$.}
\label{fig:apt-manifold}
\end{figure}

	\begin{proof}
		Under the hypotheses of Theorem~\ref{thm:invariant-manifold}, the
critical components are compact and the critical values are countable,
as shown in Appendix~\ref{sec:invariant-manifold-statement}.
Thus Lemma~\ref{lem:compact-limit-geometry} applies. Since the range of each path $w$ in the statement has compact closure
in $U$, continuity of the distance function gives
\[
\max_{z\in\omega(w)}\operatorname{dist}(z,\mathcal E_0)
=\limsup_{t\to\infty}\operatorname{dist}(w(t),\mathcal E_0) \le\rho_K(\lambda).
\]
If $K\cap\mathcal E_0=\varnothing$, compactness of $K$ and closedness
of $\mathcal E_0$ give $d:=\operatorname{dist}(K,\mathcal E_0)>0$.
For sufficiently small $\lambda$, the previous estimate would imply
$d\le\rho_K(\lambda)<d$, since $\omega(w)$ is nonempty and contained
in $K$. Hence no such path exists. We may take
$\mathfrak M_K=\varnothing$ in this case.

Suppose now that $K\cap\mathcal E_0\ne\varnothing$.
We choose neighborhoods independently of the limit set. Let
\[
\mathfrak M_K:=\{\mathcal M:\mathcal N_{\mathcal M}\cap K\ne\varnothing\}.
\]
This family is finite because the position projection of $K$ is compact
and the critical components are locally finite. The sets
$\mathcal N_{\mathcal M}$ are compact and pairwise disjoint.
		By Lemma~\ref{lem:uniform-center-reduction}, each
$\mathcal N_{\mathcal M}$ has an open neighborhood $V_{\mathcal M}$
in which every complete orbit belongs to $\mathcal I_{\mathcal M,\lambda}$.
Since the components are compact and disjoint, we can shrink these
neighborhoods to make them pairwise disjoint. Every complete orbit staying in a smaller neighborhood also stays
in the original one, so it is still contained in
$\mathcal I_{\mathcal M,\lambda}$. We use the minimum
of the finitely many parameter thresholds.
		Set $W=\bigcup_{\mathcal M\in\mathfrak M_K}V_{\mathcal M}$, so
		$K\cap\mathcal E_0\subset W$.

		Next, all limit sets in the lemma lie in $W$ for sufficiently small
		decay. Otherwise there would be $\lambda_j\downarrow0$,
		$L_j=\omega(w_j)\subset K$ for paths satisfying the stated conditions,
        and $z_j\in L_j\setminus W$.
		Compactness gives a subsequence with $z_j\to z\in K\setminus W$,
		whereas~\eqref{eq:compact-uniform-localization} gives
		\begingroup
\[
\operatorname{dist}(z_j,\mathcal E_0)
\le\max_{z\in L_j}\operatorname{dist}(z,\mathcal E_0)
=\limsup_{t\to\infty}\operatorname{dist}(w_j(t),\mathcal E_0)
\le\rho_K(\lambda_j)\longrightarrow0.
\]
\endgroup
		Since $\mathcal E_0$ is closed, $z\in K\cap\mathcal E_0\subset W$,
		a contradiction.

		Fix a path $w$ as in the statement and put $L=\omega(w)$,
with $\lambda$ below both thresholds.
As recalled in Appendix~\ref{sec:compact-limit-geometry}, $L$ is
connected and satisfies~\eqref{eq:limit-set-invariance}.
Since the open sets $V_{\mathcal M}$ are disjoint,
we have $L\subset V_{\mathcal M}$ for one component.
For every $z\in L$, the complete orbit through $z$ stays in $L$. Hence
\[
\Phi_\lambda(\mathbb R,z)\subset L\subset V_{\mathcal M}.
\]
By Lemma~\ref{lem:uniform-center-reduction},
\[
z=\Phi_\lambda(0,z)\in\Phi_\lambda(\mathbb R,z)
\subset\mathcal I_{\mathcal M,\lambda}.
\]
Since $z\in L$ was arbitrary, $L\subset\mathcal I_{\mathcal M,\lambda}$.
	\end{proof}

	\begin{proof}[Proof of Theorem~\ref{thm:invariant-manifold}]
		The critical values are countable by the argument in Appendix~\ref{sec:invariant-manifold-statement}, so
Lemma~\ref{lem:finite-time-trap}(b) gives a common compact set $K\subset U$ and trajectories whose ranges
have compact closure in $U$.
Together with Lemma~\ref{lem:corrected-tracking}, this gives a nonempty limit set
		$L_\lambda=\omega(z^\lambda)\subset K$ for each
		small fixed decay. Lemma~\ref{lem:manifold-limit-sets} gives
		one finite collection of critical components. Choose a common threshold
$\lambda_{\mathrm{man}}>0$ below both its threshold $\lambda_K$ and
the trapping threshold
		such that $L_\lambda$ lies on one $\mathcal I_{\mathcal M,\lambda}$.
		For each fixed $0<\lambda<\lambda_{\mathrm{man}}$, choose such a component
		$\mathcal M$ and set $\mathcal I_\lambda:=\mathcal I_{\mathcal M,\lambda}$.
		The $C^1$ regularity, local invariance, and dimension identity
		$\dim\mathcal I_\lambda=k_{\mathcal M}$ follow from
		Lemma~\ref{lem:uniform-center-reduction}.
		
		Compactness gives $\operatorname{dist}(z^\lambda(t),L_\lambda)\to0$:
		otherwise a sequence at positive distance from $L_\lambda$ would
		have a limit point in $L_\lambda$. Since $L_\lambda\subset\mathcal I_\lambda$,
\[
\operatorname{dist}(z^\lambda(t),\mathcal I_\lambda)
\le\operatorname{dist}(z^\lambda(t),L_\lambda)\longrightarrow0,
\]
which proves~\eqref{eq:invariant-localization}.
	\end{proof}

\clearpage
	\section{Convergence to equilibria}\label{sec:equilibrium-convergence}
	\label{sec:approximate-trajectory-convergence}
	
We first state a convergence criterion for approximate trajectories near
isolated critical points with Hessian nullity at most one
(Appendix~\ref{sec:approximate-convergence-criterion}). We then apply
it to the bias-corrected system in Appendix~\ref{sec:equilibrium-proofs}
to prove Theorem~\ref{thm:corank-one-convergence}, using the decay of
the bias-correction terms.

\subsection{A convergence criterion for approximate trajectories}
\label{sec:approximate-convergence-criterion}
When critical points are isolated, a localized connected limit set lies
near one critical point. With at most one degenerate direction, a local
finite-length estimate gives convergence to a single equilibrium.
For the bias-corrected ODE studied here, we apply the following lemma with
$r(t)\equiv0$. The case of nonzero $r$ is prepared for the stochastic analysis in the
companion paper~\citep{CompanionStoch}.

\begin{lemma}[Approximate trajectory convergence]
\label{lem:approximate-trajectory-convergence}
Suppose Assumption~\ref{ass:long-time}(i)--(iii) holds and
the critical points are isolated with Hessian nullity at most one.
For every nonempty compact $K\subset U$, there is $\lambda_K>0$
such that, for each fixed $0<\lambda\le\lambda_K$, the following holds.
Let $w:[t_0,\infty)\to U$ be continuous with
$\operatorname{dist}(w(t),K)\to0$.
Assume that there is $t_1\ge t_0$ such that $w$ is locally absolutely
continuous on $[t_1,\infty)$ and, for measurable $C,r$, satisfies almost everywhere there
\begin{equation}
\dot w=(I+C(t))f_\lambda(w)+r(t),\qquad
C(t)\longrightarrow0,\qquad\text{\normalfont and}\qquad \int_{t_1}^\infty\|r(t)\|\,dt<\infty. \label{eq:deterministic-relative-perturbation}
\end{equation}
Then $w$ converges to a point of $K\cap\mathcal E_\lambda$.
\end{lemma}
\begin{proof}
See Appendix~\ref{app:approximate-trajectory-convergence}.
\end{proof}

\subsection{Proof of Theorem~\ref{thm:corank-one-convergence}}\label{guide:appendix-5-1}
\label{sec:equilibrium-proofs}

\begin{restatement}{Theorem}{thm:corank-one-convergence}
Suppose Assumption~\ref{ass:long-time}(i)--(iii) holds, the critical
points are isolated, and $\dim\ker\nabla^2F(q)\le1$ for every
$q\in\operatorname{Crit}(F)$.
		There is $\lambda_{\mathrm{eq}}>0$ such that, for each fixed $0<\lambda<\lambda_{\mathrm{eq}}$,
		some $p_\lambda\in E_\lambda$ satisfies
		\begin{equation*}
			z^\lambda(t)\longrightarrow
			(p_\lambda,\nabla F(p_\lambda),S(p_\lambda))
			\qquad(t\to\infty).
		\end{equation*}
\end{restatement}

\begin{proof}[Proof of Theorem~\ref{thm:corank-one-convergence}]
By hypothesis, the critical points are isolated and have Hessian nullity
at most one. Their critical values are countable,
so Lemma~\ref{lem:finite-time-trap}(b) gives a compact $K\subset U$
containing every trajectory with sufficiently small decay for all sufficiently large times.
The bias factors give a matrix perturbation that tends to zero:
for $t>0$, equations~\eqref{eq:normalized} and~\eqref{eq:auto} give
\[
\dot z^\lambda(t)
=\diag\left(I_d,\frac{1}{A(t)}I_d,\frac{1}{B(t)}I_d\right)
f_\lambda(z^\lambda(t)).
\]
By~\eqref{eq:correction-factors},
\begingroup
\[
\frac{1}{A(t)}
=\frac{1}{1-e^{-at}}=1+\frac{e^{-at}}{1-e^{-at}},\qquad
\frac{1}{B(t)}
=\frac{1}{1-e^{-bt}}=1+\frac{e^{-bt}}{1-e^{-bt}}.
\]
\endgroup
Hence
\begingroup
\begin{equation}\label{eq:bias-relative-perturbation}
\dot z^\lambda(t)=(I_{3d}+C(t))f_\lambda(z^\lambda(t)),\qquad
C(t):=\diag\left(0,
       \frac{e^{-at}}{1-e^{-at}}I_d,
       \frac{e^{-bt}}{1-e^{-bt}}I_d\right).
\end{equation}
\endgroup
Fix $t_1>0$. In \eqref{eq:deterministic-relative-perturbation}, take
$w=z^\lambda$ and $r(t)\equiv0$. Then
\[
C(t)\longrightarrow0\quad\text{\normalfont and}\quad
\int_{t_1}^{\infty}\|r(t)\|\,dt=0<\infty.
\]
Choose $\lambda_{\mathrm{eq}}>0$ below both the trapping threshold and
the threshold in Lemma~\ref{lem:approximate-trajectory-convergence}.
That lemma then applies for every $0<\lambda<\lambda_{\mathrm{eq}}$.
It gives convergence to some
$e\in\mathcal E_\lambda$. The equilibrium equations identify
$e=(p_\lambda,\nabla F(p_\lambda),S(p_\lambda))$ with
$p_\lambda\in E_\lambda$, proving~\eqref{eq:corank-one-equilibrium-limit}.
\end{proof}

\clearpage
\section{Distance to the critical set}
\label{sec:quantitative-localization}

We now bound the distance from the entire limiting AdamW motion to
the critical set. The proof uses the local inequalities for $F$ and
the compactness and localization results of Appendix~\ref{sec:invariant-manifolds}.
When the gradient is small, the distance inequality
in~\eqref{eq:local-quantitative-geometry} bounds the distance to the
critical set. On the remaining parts of an orbit in a limit set,
an energy estimate bounds the distance traveled.

\subsection{Uniform geometric bounds}\label{guide:appendix-6-1}
\label{sec:quantitative-geometry}
We first verify the critical-value condition needed for trapping
and obtain geometric bounds with constants uniform on compact sets.

\begin{lemma}[Uniform local geometry]\label{lem:quantitative-geometry}
  Let $F\in C^1(\R^d)$ and assume that the first inequality in
  \eqref{eq:local-quantitative-geometry} holds near every $q\in\operatorname{Crit}(F)$.
  Then $F(\operatorname{Crit}(F))$ is at most countable.

  Fix a compact set $X\subset\R^d$ with $X\cap\operatorname{Crit}(F)\ne\varnothing$.
  Its critical values form a finite set
  $F(\operatorname{Crit}(F)\cap X)=\{c_1,\ldots,c_J\}$.
  There are constants $s_0,k>0$ such that, for every $x\in X$ with
  $\|\nabla F(x)\|\le s_0$, there is an index
  $j=j(x)\in\{1,\ldots,J\}$ such that
  \begin{equation}\label{eq:uniform-quantitative-geometry}
|F(x)-c_j|\le k\|\nabla F(x)\|^{1/(1-\theta)}.
\end{equation}
  If the second inequality in~\eqref{eq:local-quantitative-geometry}
  also holds near every $q\in\operatorname{Crit}(F)\cap X$, then, after decreasing
  $s_0$ if necessary, there is $M>0$ such that, for $x\in X$ with $\|\nabla F(x)\|\le s_0$,
  \[
\operatorname{dist}(x,\operatorname{Crit}(F))\le M\|\nabla F(x)\|^\gamma.
\]
  The same constants work throughout $X$: they may depend on $X$,
  but not on $x$ or on which critical-point neighborhood is used.
\end{lemma}

\begin{proof}
  If $p\in\operatorname{Crit}(F)\cap U_q$, the first geometric inequality gives
  $F(p)=F(q)$. Cover the compact set $\operatorname{Crit}(F)\cap X$ by finitely
  many open neighborhoods $U_{q_1},\ldots,U_{q_N}$ on which that
  inequality holds. Every critical value in $X$ is one of the finitely
  many $F(q_i)$. The balls $\overline B(0,n)$, $n\ge1$, cover
  $\R^d$, so $F(\operatorname{Crit}(F))$ is a countable union of finite sets.
  It has empty interior because every nonempty open interval is uncountable.

  For some $s_0>0$,
\[
\{x\in X:\|\nabla F(x)\|\le s_0\}
\subset\bigcup_{i=1}^N U_{q_i}.
\]
Indeed, a sequence violating this inclusion as its gradients tend to
zero would have a critical limit in $X$ outside the union. For
$x\in U_{q_i}$, the first inequality gives
\[
|F(x)-F(q_i)|\le k_{q_i}^{1/(1-\theta)}\|\nabla F(x)\|^{1/(1-\theta)}.
\]
We therefore take $k=\max_i k_{q_i}^{1/(1-\theta)}$. If the second
  inequality is available, first intersect the two neighborhoods at each
  critical point and then take a finite subcover; taking
  $M=\max_i M_{q_i}$ proves the distance estimate as well.
\end{proof}

Under Assumption~\ref{ass:long-time}(i)--(ii) and the first geometric inequality at every critical
point, this lemma gives the empty-interior critical-value condition
in Lemmas~\ref{lem:finite-time-trap}
and~\ref{lem:compact-limit-geometry}. Fix any nonempty compact
$K\subset U$ and a closed ball $X$ whose interior contains its position
projection and meets $\operatorname{Crit}(F)$. Such a ball exists because coercivity
gives a global minimizer. Use $c_1,\ldots,c_J,s_0,k$ from
Lemma~\ref{lem:quantitative-geometry} on $X$, and also $M$ when the
distance inequality is assumed. The constants below may depend on
$K,X,F,S,a,b,\varepsilon$ and the geometric bounds. They are uniform
over the invariant sets in $K$ and all sufficiently small $\lambda$. For the continuous trajectory,
we will take $K$ to be the common compact region given by
Lemma~\ref{lem:finite-time-trap}(b).

The two inequalities in~\eqref{eq:local-quantitative-geometry}
control different errors. Near a critical point $q$, they give
\[
\begin{aligned}
|F(x)-F(q)|&\le k_q^{1/(1-\theta)}\|\nabla F(x)\|^{1/(1-\theta)}
&&\text{(objective-value error)},\\
\operatorname{dist}(x,\operatorname{Crit}(F))&\le M_q\|\nabla F(x)\|^\gamma
&&\text{(distance to the critical set)}.
\end{aligned}
\]
At a saddle, different directions can cancel in the objective-value
error. For example, $G(x_1,x_2):=x_1^2-x_2^2$ has
$\operatorname{Crit}(G)=\{(0,0)\}$, and, for $s\ne0$,
\begingroup
\[
G(s,s)-G(0,0)=s^2-s^2=0,\qquad
\operatorname{dist}((s,s),\operatorname{Crit}(G))=\sqrt{2}\,|s|>0.
\]
\endgroup
The distance inequality controls this positional error.

\subsection{Energy estimates}\label{guide:appendix-6-2}
\label{sec:quantitative-energy}
\label{sec:quantitative-band}
For the autonomous field $f_\lambda$, define the nonnegative residual $R$ by
\begin{equation}\label{eq:quantitative-residual}
 R(x,m,v)^2=\|\nabla F(x)\|^2+\|m\|^2+\|S(x)-v\|^2.
\end{equation}
It vanishes precisely on $\mathcal E_0$. The next lemma gives an
energy that decreases whenever $R$ is greater than a constant times $\lambda$.

\begin{lemma}[Compact energy estimate]\label{lem:quantitative-energy}
  Suppose $F\in C^1(\R^d)$, $\nabla F$ and $S$ are locally
  Lipschitz, $S_i>0$, and $0<b\le4a$. For every compact $K\subset U$
  there are a locally Lipschitz function $W:U\to\mathbb R$ and constants $c_0,B_0>0$,
  independent of $0\le\lambda\le1$, such that, along solutions of~\eqref{eq:auto} that remain in $K$,
  \begin{equation}\label{eq:quantitative-energy-estimates}
    \dot W\le-c_0R^2+B_0\lambda R\quad\text{\normalfont a.e.}\quad\text{\normalfont and}\quad
    |W-F(x)|\le B_0R^2.
  \end{equation}
  At every maximum or minimum of $W$ on a compact complete invariant
  set $L\subset K$, one has $R\le B_0\lambda/c_0$.
\end{lemma}
\begin{proof}
See Appendix~\ref{app:quantitative-proof}.
\end{proof}

Combining the $O(\lambda)$ residual bound at the extrema of $W$ with
the local gradient inequality gives the following energy-range estimate.

\begin{lemma}[Energy range]\label{lem:quantitative-band}
  Suppose Assumption~\ref{ass:long-time}(i)--(ii) holds, together with the first inequality in
  \eqref{eq:local-quantitative-geometry} near every critical point.
  Fix a nonempty compact $K\subset U$ and a closed ball $X$ whose
  interior contains its position projection and meets $\operatorname{Crit}(F)$.
  Take $W$ from Lemma~\ref{lem:quantitative-energy} on $K$.
  There are $C,\lambda_K>0$ such that, for every
  $0<\lambda\le\lambda_K$ and every continuous path
  $w:[0,\infty)\to U$ whose range has compact closure in $U$
  and which satisfies \eqref{eq:deterministic-flow-tracking}
  with $L:=\omega(w)\subset K$,
  \begin{equation}\label{eq:quantitative-band}
    \max_{z\in L}W(z)-\min_{z\in L}W(z)\le C\lambda^{1/(1-\theta)}.
  \end{equation}
  The constants are independent of $\lambda$ and the path $w$.
\end{lemma}
\begin{proof}
See Appendix~\ref{app:quantitative-band}.
\end{proof}

\subsection{Proof of Theorem~\ref{thm:quantitative-localization}}\label{guide:appendix-6-3}
\label{sec:quantitative-excursions}

\begin{restatement}{Theorem}{thm:quantitative-localization}
Suppose Assumption~\ref{ass:long-time}(i)--(ii) holds and that \eqref{eq:local-quantitative-geometry}
  holds near every $q\in\operatorname{Crit}(F)$. Then there are
  $C_{\mathrm{ql}},\lambda_{\mathrm{ql}}>0$ such that, for $0<\lambda<\lambda_{\mathrm{ql}}$,
  \begin{equation*}
\limsup_{t\to\infty}\operatorname{dist}(x^\lambda(t),\operatorname{Crit}(F))  
\le C_{\mathrm{ql}}\left(\lambda^\gamma+ 
       \lambda^{\frac{\gamma}{(1+\gamma)(1-\theta)}}\right). 
\end{equation*}
\end{restatement}

We introduce a threshold $u>0$ for the residual $R$.
For sufficiently small $u$, the distance error bound
in~\eqref{eq:local-quantitative-geometry} gives
$\operatorname{dist}(x,\operatorname{Crit}(F))\le Mu^\gamma$ when $R\le u$.
If $u$ is at least a sufficiently large multiple of $\lambda$, energy
decrease limits how far the position can move while $R>u$.
These requirements give the range $A\lambda\le u\le s_0$ in the next
lemma, where $A$ is sufficiently large and $s_0$ sufficiently small.

\begin{lemma}[Distance bound]
  \label{lem:quantitative-excursion}
  Under the hypotheses of Lemma~\ref{lem:quantitative-energy}, suppose
  also that the distance inequality in~\eqref{eq:local-quantitative-geometry}
  holds near every critical point with exponent $0<\gamma\le1$.
  Fix a compact $K\subset U$ and a closed ball $X$ containing its
  position projection and meeting $\operatorname{Crit}(F)$.
  Take $W$ from Lemma~\ref{lem:quantitative-energy} on $K$. There are
  $A,C,M,s_0>0$ such that, for every $0<\lambda\le1$, every nonempty
  compact complete invariant set $L\subset K$ of $f_\lambda$, and
  every $u$ with $A\lambda\le u\le s_0$,
  \begin{align}
\max_{(x,m,v)\in L}\operatorname{dist}(x,\operatorname{Crit}(F)) &\le M u^\gamma+\frac C{u}\left(\max_{z\in L}W(z)-\min_{z\in L}W(z)\right). \label{eq:quantitative-excursion}
\end{align}
  The constants $A,C,M,s_0$ depend only on the fixed compact sets,
  system data, and local geometric bounds; they are independent of
  $\lambda,L,u$.
\end{lemma}
\begin{proof}
See Appendix~\ref{app:quantitative-excursion}.
\end{proof}

We now combine Lemmas~\ref{lem:quantitative-band}
and~\ref{lem:quantitative-excursion} to prove
Theorem~\ref{thm:quantitative-localization}.

\begin{proof}[Proof of Theorem~\ref{thm:quantitative-localization}]
By Lemma~\ref{lem:quantitative-geometry}, the critical values are
at most countable and hence have empty interior.
Lemma~\ref{lem:finite-time-trap}(b) therefore gives a common compact
set $K\subset U$ containing $z^\lambda(t)$ for all sufficiently large
times and all sufficiently small $\lambda>0$; their ranges have
compact closure in $U$. Lemma~\ref{lem:corrected-tracking} gives
\eqref{eq:deterministic-flow-tracking}. Thus
$L:=\omega(z^\lambda)\subset K$ is nonempty and compact, and
\eqref{eq:limit-set-invariance} holds.

Choose a closed ball $X$ whose interior contains the position
projection of $K$ and which meets $\operatorname{Crit}(F)$.
Lemmas~\ref{lem:quantitative-band} and~\ref{lem:quantitative-excursion}
give constants $A,C,M,s_0>0$, independent of $\lambda$, such that, for $A\lambda\le u\le s_0$,
\begingroup
\begin{equation}\label{eq:quantitative-threshold-bound}
\limsup_{t\to\infty}\operatorname{dist}(x^\lambda(t),\operatorname{Crit}(F))
=\max_{(x,m,v)\in L}\operatorname{dist}(x,\operatorname{Crit}(F))
\le Mu^\gamma+C\frac{\lambda^{1/(1-\theta)}}{u}.
\end{equation}
\endgroup
The equality follows from compactness of the trajectory's closure
and continuity of the distance function.
To balance the two terms in~\eqref{eq:quantitative-threshold-bound}
while ensuring $u\ge A\lambda$, choose
\[
u=A\left(\lambda+\lambda^{\frac{1}{(1+\gamma)(1-\theta)}}\right).
\]
Then $u\ge A\lambda$ and $u\to0$ as $\lambda\downarrow0$, so
$u\le s_0$ for sufficiently small $\lambda$.
Since $0<\gamma\le1$,
\begingroup
\[
u^\gamma
\le A^\gamma\left(\lambda^\gamma+
\lambda^{\frac{\gamma}{(1+\gamma)(1-\theta)}}\right),\qquad
\frac{\lambda^{1/(1-\theta)}}{u}
\le A^{-1}\lambda^{\frac{1}{1-\theta}-\frac{1}{(1+\gamma)(1-\theta)}}
=A^{-1}\lambda^{\frac{\gamma}{(1+\gamma)(1-\theta)}}.
\]
\endgroup
Substituting into~\eqref{eq:quantitative-threshold-bound}
proves~\eqref{eq:quantitative-localization},
with $C_{\mathrm{ql}}=MA^\gamma+C/A$ and a sufficiently small common
threshold $\lambda_{\mathrm{ql}}>0$.
\end{proof}

For the stochastic application in \citet{CompanionStoch}, we record
the same bound for general paths satisfying
\eqref{eq:deterministic-flow-tracking}.

\begin{lemma}[Localization of tracking paths]
\label{lem:quantitative-path-localization}
Suppose Assumption~\ref{ass:long-time}(i)--(ii) holds and that both
inequalities in~\eqref{eq:local-quantitative-geometry} hold near every
$q\in\operatorname{Crit}(F)$, with common exponents $0<\theta\le1/2$ and
$0<\gamma\le1$. Then $F(\operatorname{Crit}(F))$ is at most countable.
For every nonempty compact $K\subset U$, there are $C_K,\lambda_K>0$
such that, for each fixed $0<\lambda\le\lambda_K$, every continuous
path $w=(w_x,w_m,w_v):[t_0,\infty)\to U$ satisfying
\eqref{eq:deterministic-flow-tracking} and
$\operatorname{dist}(w(t),K)\to0$ satisfies
\begin{equation}\label{eq:quantitative-path-localization}
\limsup_{t\to\infty}\operatorname{dist}(w_x(t),\operatorname{Crit}(F))
\le C_K\left(\lambda^\gamma+
\lambda^{\frac{\gamma}{(1+\gamma)(1-\theta)}}\right).
\end{equation}
The constants depend only on $K$, $F,S,a,b,\varepsilon$, and the
geometric bounds; they are independent of $\lambda$ and the path $w$.
\end{lemma}

\begin{proof}
Countability follows from Lemma~\ref{lem:quantitative-geometry}.
Since $\operatorname{dist}(w(t),K)\to0$, the path eventually lies in a
compact neighborhood of $K$ contained in $U$. Continuity on the initial
finite interval then gives compact closure of its range in $U$, so
$L:=\omega(w)$ is nonempty and compact, with $L\subset K$.
The tracking condition~\eqref{eq:deterministic-flow-tracking} also gives
\eqref{eq:limit-set-invariance}.

Choose a closed ball $X$ whose interior contains the position projection
of $K$ and a minimizer of $F$, which exists by coercivity.
Lemmas~\ref{lem:quantitative-band} and~\ref{lem:quantitative-excursion}
give constants $A,C,M,s_0>0$, independent of $w$ and of sufficiently
small $\lambda>0$, such that, for $A\lambda\le u\le s_0$,
\begingroup
\[
\limsup_{t\to\infty}\operatorname{dist}(w_x(t),\operatorname{Crit}(F))
=\max_{(x,m,v)\in L}\operatorname{dist}(x,\operatorname{Crit}(F))
\le Mu^\gamma+C\frac{\lambda^{1/(1-\theta)}}{u}.
\]
\endgroup
The equality follows from compactness and continuity of the distance
function. Taking $u=A(\lambda+\lambda^{\frac{1}{(1+\gamma)(1-\theta)}})$ and using
the calculation in the proof of Theorem~\ref{thm:quantitative-localization}
gives~\eqref{eq:quantitative-path-localization} with
$C_K=MA^\gamma+C/A$ and a sufficiently small $\lambda_K>0$.
All constants are uniform over the paths $w$.
\end{proof}

\subsection{Examples}\label{guide:appendix-6-4}
\label{sec:quantitative-examples}
Under the hypotheses of Theorem~\ref{thm:quantitative-localization},
we first compute the localization rate for nondegenerate and Morse--Bott
critical sets, then give a one-dimensional example attaining the bound
at a degenerate minimum.

\paragraph{Nondegenerate and Morse--Bott critical sets.}
If a critical component is Morse--Bott, its Hessian kernel equals its
tangent space and the Hessian is invertible in the normal directions.
Near a compact part of that component, take a nearest point $q$ and
write $x=q+h$, with $h$ normal. Taylor expansion, uniformly on this
compact part, gives
\[
\nabla F(x)=\nabla^2F(q)h+o(\|h\|)
\qquad\text{\normalfont and}\qquad |F(x)-F(q)|\le C\|h\|^2.
\]
The normal eigenvalues have absolute values bounded away from zero,
so $\|\nabla F(x)\|\ge c\|h\|$ after shrinking the neighborhood. Hence
\begingroup
\[
|F(x)-F(q)|^{1/2}\le\sqrt C\|h\|
\le\frac{\sqrt C}{c}\|\nabla F(x)\|,\qquad
\operatorname{dist}(x,\operatorname{Crit}(F))\le\|h\|
\le\frac1c\|\nabla F(x)\|.
\]
\endgroup
Thus~\eqref{eq:local-quantitative-geometry} holds with
$\theta=1/2$ and $\gamma=1$, giving a radius $O(\lambda)$.
The restriction of $\nabla^2F(q)$ to $(T_q\mathcal M)^\perp$ may
have both positive and negative eigenvalues. A nondegenerate critical point
is the zero-dimensional case. For a positive-dimensional critical
component, this bounds the normal distance of the limiting AdamW motion.

The next one-dimensional example attains the order in
Theorem~\ref{thm:quantitative-localization} for a degenerate minimum.

\begin{example}[Attaining the localization order]
\label{ex:quantitative-polynomial-sharp}
For an even integer $p\ge4$, set
\begin{equation*}
F(x)=\frac{(x-1)^p}{p},\qquad f(x,\xi)=F(x)+\xi x,\qquad\text{\normalfont and}\qquad
\mathbb P(\xi=1)=\mathbb P(\xi=-1)=\tfrac12.
\end{equation*}
Then
\begin{equation*}
\mathbb E f(x,\xi)=F(x),\qquad F'(x)=(x-1)^{p-1},\qquad\text{\normalfont and}\qquad
S(x)=\mathbb E[(\partial_x f(x,\xi))^2] =1+(x-1)^{2p-2}.
\end{equation*}
The two sample losses are smooth and coercive, and $S\ge1$.
Fix $a=b=\varepsilon=1$ and any initial position $x_0$, independently
of $\lambda$. The unique critical point is $1$, and
\[
|F(x)-F(1)|^{1-1/p}=p^{-(p-1)/p}|F'(x)|
\qquad\text{\normalfont and}\qquad |x-1|=|F'(x)|^{1/(p-1)}.
\]
Thus the geometric bounds hold with $\theta=1/p$ and
$\gamma=1/(p-1)$. Since $(1+\gamma)(1-\theta)=1$,
both terms in Theorem~\ref{thm:quantitative-localization} have the same
order, giving
\[
\limsup_{t\to\infty}|x^\lambda(t)-1|\le C\lambda^{1/(p-1)}.
\]

The Hessian nullity at $1$ is one, so
Theorem~\ref{thm:corank-one-convergence} gives
$x^\lambda(t)\to p_\lambda\in E_\lambda$ for every sufficiently
small fixed $\lambda>0$. Theorem~\ref{thm:quantitative-localization}
therefore gives $|p_\lambda-1|\le C\lambda^{1/(p-1)}\to0$.
The equilibrium equation is
\[
(p_\lambda-1)^{p-1}
=-\lambda p_\lambda\left(1+
\sqrt{1+(p_\lambda-1)^{2p-2}}\right).
\]
For small $\lambda$, $p_\lambda>0$, so this equation gives
$p_\lambda<1$ and
\[
\frac{(1-p_\lambda)^{p-1}}{\lambda}
=p_\lambda\left(1+\sqrt{1+(1-p_\lambda)^{2p-2}}\right)
\longrightarrow2.
\]
Consequently, as $\lambda\downarrow0$,
\begin{equation}\label{eq:quantitative-polynomial-sharp}
\lim_{t\to\infty}|x^\lambda(t)-1|
=1-p_\lambda =(2\lambda)^{1/(p-1)}(1+o(1)).
\end{equation}
For $p=4$, this is $(2\lambda)^{1/3}(1+o(1))$.
\end{example}

In this family, weight decay shifts the limiting AdamW position away
from the unregularized minimizer by order $\lambda^{1/(p-1)}$.
The exponent in Theorem~\ref{thm:quantitative-localization} is sharp
for this fixed sampling law and fixed initialization.

\clearpage
\section{Proofs of auxiliary lemmas}
\label{app:detailed-proofs}
We first prove the trapping and flow-tracking estimates
(Appendices~\ref{app:trapping-finite-length} and~\ref{app:corrected-tracking}),
then establish localization of compact limit sets
(Appendix~\ref{app:compact-invariant}), the local center reduction
(Appendix~\ref{app:center-reduction}), and finite-length convergence
(Appendix~\ref{app:approximate-trajectory-convergence}). The quantitative
arguments follow: energy dissipation (Appendix~\ref{app:quantitative-proof}),
the energy-range bound (Appendix~\ref{app:quantitative-band}), and the
excursion estimate (Appendix~\ref{app:quantitative-excursion}).
The nonconvergent example is developed and proved in
Appendix~\ref{app:counterexample-details}.

\subsection{Proof of Lemma~\ref{lem:finite-time-trap}}\label{guide:appendix-7-1}
\label{app:trapping-finite-length}
\begin{proof}
We first construct a strictly invariant region for the zero-decay flow,
then use perturbation estimates to trap the corrected trajectory.

\noindent\textbf{Step 1: Construction and strict invariance of the trapping region (part (a)).}
		Translate $F$ so that $F\ge0$, and let $K_0$ be the compact closure
of the corrected zero-decay trajectory from
Lemma~\ref{lem:hp-zero-decay-energy}. Recall the autonomous energy
\[
V(x,m,v)=F(x)+\frac1{2a}\sum_i\frac{m_i^2}{\varepsilon+\sqrt{v_i}}.
\]
Since $F(\operatorname{Crit}(F))$ has empty interior, choose
$c>\sup_{K_0}V$ outside this set and put
		$X=\{F\le c\}$, which is compact by coercivity. Choose
		\begingroup
\[
0<\ell<\min_{X,i}S_i,\qquad u>\max_{X,i}S_i,\qquad
B=\{V\le c,\ \ell\le v_i\le u\}.
\]
\endgroup
		This set is compact: $x\in X$, and its first moments satisfy
		$\|m\|^2\le2a(\varepsilon+\sqrt u)c$.
		For $z^0(t)\in K_0$ we have $F(x^0(t))\le V(z^0(t))<c$.
The second-moment kernel has mass one, so, also on taking the closure,
\[
V|_{K_0}\le\sup_{K_0}V<c,\qquad
\ell<\min_X S_i\le v_i|_{K_0}\le\max_X S_i<u.
\]
Thus $K_0\subset\operatorname{int}B$.
		The dissipation inequality~\eqref{eq:adam-energy-dissipation}
		gives $V(t)\le V(0)\le c$ while the solution is defined in $U$,
so $x(t)\in X$. The second-moment equation gives
\[
e^{-bt}v_i(0)+(1-e^{-bt})\min_X S_i
\le v_i(t)\le
e^{-bt}v_i(0)+(1-e^{-bt})\max_X S_i.
\]
Since $v_i(0)\in[\ell,u]$ and
$\ell<\min_X S_i\le\max_X S_i<u$, we obtain
\[
\ell\le e^{-bt}v_i(0)+(1-e^{-bt})\min_X S_i
\le v_i(t)
\le e^{-bt}v_i(0)+(1-e^{-bt})\max_X S_i\le u.
\]
Together with $V(t)\le c$, this gives $\Phi_0(t,B)\subset B$
for as long as the solutions are defined.
		By compactness, each solution extends to all positive times.
		
For every $z\in B$ and $t>0$, $V(\Phi_0(t,z))<c$.
		Indeed, equality would require $V(z)=c$ and $m\equiv0$ throughout
		$[0,t]$. Then $x$ is constant and $0=\dot m=a\nabla F(x)$, so
$F(x)=c$ is a critical value, contrary to the choice of $c$. For $t>0$, the bounds for $v_i(t)$ above are strict relative to
$\ell$ and $u$, since $1-e^{-bt}>0$ and
$\ell<\min_X S_i\le\max_X S_i<u$.
For each fixed $t>0$, compactness of $B$ and continuity of
$\Phi_0(t,\cdot)$ give~\eqref{eq:finite-time-inward}.
		
		\noindent\textbf{Step 2: Trapping the corrected trajectory (part (b)).}
		
Choose a compact neighborhood $K\subset U$ of $B$ with
$B\subset\operatorname{int}K$. Apply Lemma~\ref{lem:unit-interval-trap}
to obtain the perturbation tolerance $\kappa$.

For $t>0$, subtracting the autonomous zero-decay field gives
\[
\dot z^\lambda(t)-f_0(z^\lambda(t))
=\begin{pmatrix}
-\lambda x^\lambda(t)\\
\dfrac{ae^{-at}}{1-e^{-at}}(\nabla F(x^\lambda(t))-m^\lambda(t))\\
\dfrac{be^{-bt}}{1-e^{-bt}}(S(x^\lambda(t))-v^\lambda(t))
\end{pmatrix}.
\]
On $K$, once $e^{-at},e^{-bt}\le1/2$, this gives
\[
\|\dot z^\lambda(t)-f_0(z^\lambda(t))\|
\le C_K(\lambda+e^{-at}+e^{-bt}).
\]
Choose $\lambda_0>0$ with $C_K\lambda_0\le\kappa/4$.
The sets $B,K$ and the tolerance $\kappa$ are now fixed.
Since $K_0\subset\operatorname{int}B$, put
$d_0=\operatorname{dist}(K_0,B^c)>0$.
Fix any $T>1$ sufficiently large that
$e^{-aT},e^{-bT}\le1/2$ and
$C_K(e^{-aT}+e^{-bT})\le\kappa/4$.
Since $z^0(T)\in K_0$, its distance from $B^c$ is at least $d_0$.

Continuous dependence on decay holds up to this time.
		Choose a closed ball $D$ with
$d_T:=\operatorname{dist}(x^0([0,T]),D^c)>0$. On $D$,
		Lemma~\ref{lem:wp-lipschitz} bounds the path-dependent ratio by a
		Lipschitz constant $L_T$. If $B_T$ bounds positions there, integration of
$\dot x^\lambda=-R_{x^\lambda}-\lambda x^\lambda$ gives, up to first exit,
\begingroup
\[
\sup_{s\le t}\|x^\lambda(s)-x^0(s)\|
\le \lambda B_Tt
 +L_T\int_0^t\sup_{u\le r}\|x^\lambda(u)-x^0(u)\|\,dr
\le\lambda B_Tt e^{L_Tt},\qquad 0\le t\le T,
\]
\endgroup
where the second line follows from Gronwall's inequality.
		At a first boundary hit $\tau\le T$ the estimate would give
\[
d_T\le\|x^\lambda(\tau)-x^0(\tau)\|
\le\lambda B_T T e^{L_TT}<d_T
\]
for small enough $\lambda$. Thus the position estimate holds on
$[0,T]$. The moment kernels have mass one, so with $L_D$ the sum of
the Lipschitz constants of $\nabla F$ and $S$ on $D$,
\[
\sup_{0\le t\le T}
\bigl(\|m^\lambda(t)-m^0(t)\|+\|v^\lambda(t)-v^0(t)\|\bigr)
\le L_D\sup_{0\le t\le T}\|x^\lambda(t)-x^0(t)\|.
\]
At $t=0$ the moment differences are zero.
		Decrease $\lambda_0$ so that $\|z^\lambda(T)-z^0(T)\|<d_0/2$.
		Then
\begingroup
\[
\operatorname{dist}(z^\lambda(T),B^c)
\ge d_0-\|z^\lambda(T)-z^0(T)\|>d_0/2.
\]
\endgroup For $t\ge T$, the corrected trajectory differs from $f_0$ by at
		most $\kappa/2$ while in $K$. Applying Lemma~\ref{lem:unit-interval-trap} to
		$w(s)=z^\lambda(T+s)$ gives
		$z^\lambda(t)\in K$ for all $t\ge T$.
		The sets and $\kappa$ were chosen before $T$; only the decay
threshold is reduced after $T$ is fixed.

By the finite-time continuous dependence proved above, the paths $x^\lambda,m^\lambda,v^\lambda$
are uniformly bounded on $[0,T]$ after decreasing $\lambda_0$ if necessary.
The positions remain in the ball $D$ chosen above. Since $S_i>0$ on $D$,
		\[
v_i^\lambda(t)
		=\frac b{B(t)}\int_0^t e^{-b(t-r)}S_i(x^\lambda(r))\,dr \ge\min_{x\in D}S_i(x) >0\qquad(0<t\le T).
\]
		The same lower bound holds at $t=0$, since $v_i^\lambda(0)=S_i(x_0)$.
		Thus the trajectories on $[0,T]$, together with their tails in $K$,
		have compact closure in $U$.

	\end{proof}

\subsection{Proof of Lemma~\ref{lem:corrected-tracking}}\label{guide:appendix-7-2}
\label{app:corrected-tracking}
We compare the corrected trajectory with the autonomous flow.

\begin{proof}
By Lemma~\ref{lem:finite-time-trap}(b), choose $T,\lambda_0>0$ and a
compact set $K\subset U$ such that
\[
z^\lambda(t)\in K\qquad(t\ge T,\ 0\le\lambda\le\lambda_0),
\]
and the trajectories have compact closure in $U$. Fix
$0\le\lambda\le\lambda_0$.
Choose $T_0\ge T$ such that $e^{-aT_0},e^{-bT_0}\le1/2$.
		For $t\ge T_0$, subtracting the autonomous field from the corrected
equation gives, with $(x,m,v)=z^\lambda(t)\in K$,
\begin{gather*}
\dot z^\lambda(t)-f_\lambda(z^\lambda(t))
=\left(0,
\frac{ae^{-at}}{1-e^{-at}}(\nabla F(x)-m),
\frac{be^{-bt}}{1-e^{-bt}}(S(x)-v)\right),\\
\|\dot z^\lambda(t)-f_\lambda(z^\lambda(t))\|
\le C_K e^{-ct}\quad\text{\normalfont and}\quad c=\min\{a,b\}.
\end{gather*}
		Multiply $f_\lambda$ by a smooth cutoff supported in $U$ and equal
		to one on a neighborhood of $K$. The resulting field extends by zero
		to a globally Lipschitz field; write $\widetilde\Phi_\lambda$ for its
		flow and choose a Lipschitz constant $L_K>0$.
		For fixed $H>0$, integration of the two equations gives, for $0\le s\le H$,
		\[
		\begin{aligned}
		\|z^\lambda(t+s)-\widetilde\Phi_\lambda(s,z^\lambda(t))\|
		&\le\int_0^s\left[
		L_K\|z^\lambda(t+r)-\widetilde\Phi_\lambda(r,z^\lambda(t))\|
		+C_K e^{-c(t+r)}\right]dr.
		\end{aligned}
		\]
		Hence Gronwall's inequality gives
		\[
\sup_{0\le s\le H}
		\|z^\lambda(t+s)-\widetilde\Phi_\lambda(s,z^\lambda(t))\|
		\le C_K e^{-ct}\int_0^H e^{L_K(H-r)}\,dr =C_K e^{-ct}\frac{e^{L_KH}-1}{L_K}\longrightarrow0.
\]
		Let $O$ be the open neighborhood of $K$ where the cutoff equals one.
        Set $d=\operatorname{dist}(K,O^c)>0$. For all sufficiently large
$t$ and $0\le s\le H$, we have
\[
\operatorname{dist}(\widetilde\Phi_\lambda(s,z^\lambda(t)),K)
\le\|\widetilde\Phi_\lambda(s,z^\lambda(t))-z^\lambda(t+s)\|
\le C_K e^{-ct}\frac{e^{L_KH}-1}{L_K}<d.
\]
Thus $\widetilde\Phi_\lambda(s,z^\lambda(t))\in O$ throughout $[0,H]$.
		It therefore agrees with $\Phi_\lambda$ throughout $[0,H]$.
		Thus $z^\lambda$ is an asymptotic pseudotrajectory whose range
has compact closure in $U$.
\end{proof}

\subsection{Proof of Lemma~\ref{lem:compact-limit-geometry}}\label{guide:appendix-7-3}
\label{app:compact-invariant}

We pass to a Hausdorff limit of invariant sets as decay vanishes,
then use the zero-decay energy to identify that limit.

\begin{proof}
A compact invariant set is internally chain transitive (ICT) if, for every $\epsilon,H>0$,
any two of its points can be joined by finitely many orbit segments
inside the set, each of duration at least $H$, with successive jumps
of size less than $\epsilon$.
By Theorem~5.7(i) of~\cite{B1999}, if a path is an asymptotic
pseudotrajectory of a continuous flow and its range has compact closure,
its limit set is internally chain transitive. By Proposition~6.4 of~\cite{B1999}, if a continuous
Lyapunov function is constant along orbits in a compact invariant set
$A$, strictly decreases along every orbit outside $A$, and has a set of
values on $A$ with empty interior, every internally chain
transitive set is contained in $A$.

Fix a nonempty compact set $K\subset U$ and $\lambda\ge0$.
Let $w:[0,\infty)\to U$ be continuous, with compact closure of its
range in $U$, and put $L=\omega(w)\subset K$. The tracking condition is
\[
\sup_{0\le s\le H}\|w(t+s)-\Phi_\lambda(s,w(t))\|\longrightarrow0
\qquad(t\to\infty)
\]
for every $H>0$, with the flow defined on each such interval for all
sufficiently large $t$. Choose a smooth cutoff supported in $U$ and
equal to one on a neighborhood of the compact closure of the path. The extended field is globally Lipschitz.
For each fixed $H$, the tracking condition keeps its comparison
orbits in that neighborhood for all sufficiently large times, so $w$
is an asymptotic pseudotrajectory of the extended flow.
By Theorem~5.7(i), $L$ is compact, invariant and ICT.
The extended and original fields agree near $L$. By uniqueness and
compactness, the solutions of the original equation through points of $L$
are defined for all real times and remain in $L$.
Also,
\[
L=\bigcap_{s\ge0}\overline{w([s,\infty))}
\]
is nonempty and connected, since the sets on the right are nested,
nonempty, compact and connected.
Keeping $K$ fixed, we first extract a compact invariant limit as decay
vanishes, then establish its internal chain transitivity, and finally
use the zero-decay energy to localize it to the equilibrium set.

\noindent\textbf{Step 1: Compactness and invariance of the limit.}
		Suppose $\lambda_j\downarrow0$ and
		$L_j=\omega(w_j)\subset K$, where each $w_j$ is a continuous path
whose range has compact closure in $U$ and which tracks $f_{\lambda_j}$.
For nonempty compact sets $E,G\subset K$, the Hausdorff distance is
\[
d_H(E,G):=\max\left\{
\sup_{z\in E}\operatorname{dist}(z,G),\,
\sup_{z\in G}\operatorname{dist}(z,E)\right\}.
\]
The nonempty compact subsets of the compact set $K$ form a compact
metric space for $d_H$. Hence, after taking a subsequence,
$L_j\to L$ in Hausdorff distance for a nonempty compact $L\subset K$.
		The limit is connected: two separated compact pieces of $L$ would
		give disjoint neighborhoods covering and separating $L_j$ for large $j$.
		For $p\in L$ choose $p_j\in L_j$ with $p_j\to p$ and set
$z_j(t)=\Phi_{\lambda_j}(t,p_j)$ for $t\in\mathbb R$. These orbits stay
in $K$, where their speeds are uniformly bounded. Apply Arzel\'a--Ascoli
successively on $[-N,N]$, $N=1,2,\ldots$, taking nested subsequences.
Selecting the $j$th element of the $j$th subsequence gives, after relabeling,
a single subsequence and a continuous $z:\mathbb R\to K$ such that
\[
\sup_{|t|\le N}\|z_j(t)-z(t)\|\longrightarrow0
\qquad\text{for every }N\ge1.
\]
For every $t\in\mathbb R$,
\[
\begin{aligned}
\operatorname{dist}(z(t),L)&\le\|z(t)-z_j(t)\|+d_H(L_j,L)\longrightarrow0,\\
z(t)&=\lim_{j\to\infty}\left(p_j+\int_0^t f_{\lambda_j}(z_j(s))\,ds\right)
=p+\int_0^t f_0(z(s))\,ds.
\end{aligned}
\]
The integral converges by uniform convergence of $f_{\lambda_j}$ on
$K$ and uniform continuity of $f_0$ there. By uniqueness, this is the
orbit through $p$ for all real times. Thus $\Phi_0(t,L)=L$.

\noindent\textbf{Step 2: Internal chain transitivity.}

We subdivide long chain segments so that the same finite-time comparison
applies to every jump. Fix endpoints in
		$L$, a tolerance $\epsilon>0$ and a minimum duration $H>0$.
		Choose a compact convex box in $U$ containing $K$ in its interior,
		and a smooth cutoff equal to one near $K$ and supported inside the
		box. The fields have common Lipschitz bounds on this box: the
		coefficients are locally Lipschitz and the $v_i$ are bounded away
		from zero. Multiplying by the cutoff and extending by zero gives
		globally Lipschitz fields $\widetilde f_\lambda$ with a common
		constant $L_*\ge1$. If
		$\delta_j=\sup\|\widetilde f_{\lambda_j}-\widetilde f_0\|=O(\lambda_j)$,
		Gronwall gives
		\[
		\|\widetilde\Phi_{\lambda_j}(t,p)-\widetilde\Phi_0(t,q)\|
\le e^{L_*|t|}\|p-q\|+\delta_j\frac{e^{L_*|t|}-1}{L_*}.
\]
Here $\widetilde\Phi_\lambda$ denotes the extended flow. It agrees
with $\Phi_\lambda$ on the corresponding invariant sets.
		Put $d_j=d_H(L_j,L)$ and
		\[
		e_j=\sup_{0\le t\le2H,\ z\in K}
\|\widetilde\Phi_{\lambda_j}(t,z)-\widetilde\Phi_0(t,z)\|
\longrightarrow0.
		\]
		
Choose an $\epsilon/3$-chain in $L_j$ between points within $d_j$
		of the prescribed endpoints. A segment of duration $T\ge H$ splits
		into $n=\lfloor T/H\rfloor$ equal pieces. Since
        $nH\le T<(n+1)H\le 2nH$, their duration $T/n$ lies in $[H,2H]$.
        The inserted nodes lie on the original orbit, so their jumps are zero. Replace every node by a point
		of $L$ within $d_j$, using the prescribed endpoints exactly. If $C_H$ is a Lipschitz constant of
$\widetilde\Phi_0(t,\cdot)$ for $0\le t\le2H$,
		write $p_i$ for the original nodes and $y_i$ for their replacements.
		For every subdivided duration $T_i$, the triangle inequality gives
		\[
		\begin{aligned}
		\|\Phi_0(T_i,y_i)-y_{i+1}\|
		&\le C_H\|y_i-p_i\|+e_j
		+\|\Phi_{\lambda_j}(T_i,p_i)-p_{i+1}\|+\|p_{i+1}-y_{i+1}\|\\
		&\le\epsilon/3+e_j+(C_H+1)d_j<\epsilon.
		\end{aligned}
		\]
		This includes both endpoints, whose replacement errors are at most
		$d_j$. The bound is independent of the number of segments,
		and invariance keeps the new flow segments in $L$.
		
		\noindent\textbf{Step 3: Localization by the zero-decay energy.}

Recall
\[
V(x,m,v)=F(x)+\frac1{2a}\sum_i\frac{m_i^2}{\varepsilon+\sqrt{v_i}}
\quad\text{\normalfont and}\quad
\mathcal E_0=\{(q,0,S(q)):\nabla F(q)=0\}.
\]
We show that $V$ strictly decreases on every orbit in $L$ outside
$\mathcal E_0$.
If $m=0$ and $\nabla F(x)=0$ at a point of $L$, its complete orbit
		has constant $x$, zero $m$, and
		\[
		v(t)=S(x)+(v(0)-S(x))e^{-bt}.
		\]
		Since $v(t)$ remains bounded as $t\to-\infty$, we must have
$v(0)=S(x)$. Thus the point belongs to $\mathcal E_0$. At every other point of $L$, either $m\ne0$ or
		$\dot m=a\nabla F(x)\ne0$ when $m=0$.
		In the latter case, $m(t)=ta\nabla F(x(0))+o(t)$ is nonzero for all
		sufficiently small positive $t$; in the former, continuity suffices.
		Integrating~\eqref{eq:adam-energy-dissipation}, with
$d_i(s)=\varepsilon+\sqrt{v_i(s)}$, gives
\[
V(\Phi_0(t,z))-V(z)
\le-\varepsilon\int_0^t\sum_i\frac{m_i(s)^2}{d_i(s)^2}\,ds <0
\]
for every $t>0$ and $z\in L\setminus\mathcal E_0$.
By uniqueness, such an orbit cannot reach an equilibrium in finite time.
		
A minimizer of $V$ on $L$ must therefore be an equilibrium.
		Consequently $A=L\cap\mathcal E_0$ is nonempty, compact and invariant,
		and $V(A)\subset F(\operatorname{Crit}(F))$. The latter set has empty
interior by hypothesis, so $V(A)$ does as well.
		Proposition~6.4 of~\cite{B1999}, applied to the flow restricted to $L$
		with Lyapunov function $V|_L$, gives
		$L\subset\mathcal E_0$.
		
If the limit sets did not approach $\mathcal E_0$ uniformly over the
        paths satisfying these conditions, there would be
$\epsilon>0$, $\lambda_j\downarrow0$, paths $w_j$ as above, and $z_j\in L_j=\omega(w_j)$
with $\operatorname{dist}(z_j,\mathcal E_0)\ge\epsilon$.
Along a subsequence, $L_j\to L$ in Hausdorff distance and $z_j\to z\in L$.
Since $L\subset\mathcal E_0$,
\[
\epsilon\le\operatorname{dist}(z_j,\mathcal E_0) \le\|z_j-z\|\longrightarrow0,
\]
a contradiction. Applying Proposition~6.4 of~\cite{B1999}
directly to $\omega(w)$ for a path tracking $f_0$ proves the
		claim for zero decay. If $K\cap\mathcal E_0=\varnothing$, the
		same subsequence argument shows that, for every sufficiently small
        $\lambda\ge0$, no such path has its limit set in $K$.
For sufficiently small $\lambda$, define $\rho_K(\lambda)$ as the
supremum of $\max_{z\in\omega(w)}\operatorname{dist}(z,\mathcal E_0)$
over all continuous $w:[0,\infty)\to U$ with compact closure of their
ranges in $U$, satisfying~\eqref{eq:deterministic-flow-tracking} and
$\omega(w)\subset K$, with value zero if there are none.
The argument above shows that this supremum is finite and that
$\rho_K(\lambda)\to0$ and $\rho_K(0)=0$.
Finally, take $t_j\to\infty$ along which
$\operatorname{dist}(w(t_j),\mathcal E_0)$ tends to its limit superior.
By compactness, along a subsequence $w(t_j)\to z_*\in\omega(w)$, so
\[
\limsup_{t\to\infty}\operatorname{dist}(w(t),\mathcal E_0)
=\operatorname{dist}(z_*,\mathcal E_0)
\le\max_{z\in\omega(w)}\operatorname{dist}(z,\mathcal E_0).
\]
Conversely, for each $z\in\omega(w)$ there are $s_j\to\infty$ with
$w(s_j)\to z$. Taking limits gives
$\operatorname{dist}(z,\mathcal E_0)\le
\limsup_{t\to\infty}\operatorname{dist}(w(t),\mathcal E_0)$.
Thus
\[
\limsup_{t\to\infty}\operatorname{dist}(w(t),\mathcal E_0)
=\max_{z\in\omega(w)}\operatorname{dist}(z,\mathcal E_0),
\]
which proves \eqref{eq:compact-uniform-localization}.

	\end{proof}

\subsection{Proof of Lemma~\ref{lem:uniform-center-reduction}}\label{guide:appendix-7-4}
\label{app:center-reduction}
The proof below and Appendix~\ref{app:approximate-trajectory-convergence}
use the following consequence of Theorem~9.1(i) of
\citet[p.~74]{Fenichel1979}, with $r=1$.

\phantomsection\label{app:fenichel-center-theorem}
\begin{fenicheltheorem}
Let $\mathcal P$ be a smooth finite-dimensional manifold and let
$X_\lambda$ be a jointly $C^1$ family of vector fields on $\mathcal P$,
with $\lambda$ in an open interval about zero. Let $\mathcal A$ be a
$C^1$ embedded submanifold consisting of equilibria of $X_0$, and let
$K_*\subset\mathcal A$ be compact. Assume that the operator induced by
$\mathrm D X_0(p)$ on
\[
T_p\mathcal P/T_p\mathcal A
\]
has constant numbers of eigenvalues in the left half-plane, on the
imaginary axis, and in the right half-plane as $p$ varies in $K_*$,
counted with algebraic multiplicity.

Then the extended system
\[
\dot p=X_\lambda(p),\qquad \dot\lambda=0
\]
has a $C^1$ embedded locally invariant center manifold
$\widehat{\mathcal C}$ containing $K_*\times\{0\}$. If
$\widehat E^c_p$ is the generalized center space of the extended
linearization at $(p,0)$, then
\[
T_{(p,0)}\widehat{\mathcal C}=\widehat E^c_p
\qquad(p\in K_*).
\]
There is an open neighborhood $\mathcal V$ of $K_*\times\{0\}$ such
that every complete extended orbit whose closure is contained in
$\mathcal V$ lies in $\widehat{\mathcal C}$.
\end{fenicheltheorem}

\begin{proof}
Let $\mathcal M$ be the nonempty compact connected $C^1$ embedded
manifold of critical points, without boundary. Its Hessian nullity is constant:
\[
\dim\ker\nabla^2F(q)=k\quad(q\in\mathcal M)\quad\text{\normalfont and}\quad
\mathcal N_{\mathcal M}=\{(q,0,S(q)):q\in\mathcal M\}.
\]
We write $\mathrm D S(q)$ for the Jacobian of $S$ at $q$.
Since $S_i>0$, $F\in C^2$ and $S\in C^1$ near $\mathcal M$,
the field $f_\lambda$ is jointly $C^1$ near $\mathcal N_{\mathcal M}\times\{0\}$. The map
		$q\mapsto(q,0,S(q))$ identifies $\mathcal N_{\mathcal M}$ with the compact
		$C^1$ manifold $\mathcal M$.
The proof has three steps: linearization, parameter extension, and
construction of the invariant manifolds.

\noindent\textbf{Step 1: The linearization.}
		
At $z_q=(q,0,S(q))$, the first moment vanishes, so the derivative
with respect to $v$ of the position equation is zero. Set
		$H(q)=\nabla^2F(q)$ and
		$D(q)=\diag(\varepsilon+\sqrt{S_i(q)})$. Direct differentiation gives
		\[
		J(q)=D_z f_0(z_q)=
		\begin{pmatrix}
			0&-D(q)^{-1}&0\\
			a H(q)&-a I&0\\
			b \mathrm D S(q)&0&-b I
		\end{pmatrix},
		\]
		and
		\[
		\det(rI-J(q))=(r+b)^d\det(r^2I+a rI+a D(q)^{-1}H(q)).
		\]
		
		To compute the determinant, write $D=D(q)$ and $H=H(q)$.
For $r\ne0$, adding $(a/r)H$ times the first block row to the second gives
\[
\begin{pmatrix}I&0\\(a/r)H&I\end{pmatrix}
\begin{pmatrix}rI&D^{-1}\\-aH&(r+a)I\end{pmatrix}
=\begin{pmatrix}rI&D^{-1}\\0&(r+a)I+(a/r)HD^{-1}\end{pmatrix}.
\]
The leftmost matrix has determinant one, so
\begingroup
\[
\begin{aligned}
\det\begin{pmatrix}rI&D^{-1}\\-aH&(r+a)I\end{pmatrix}
&=r^d\det\left((r+a)I+\frac ar HD^{-1}\right)\\
&=\det(r^2I+arI+aHD^{-1})
=\det(r^2I+arI+aD^{-1}H).
\end{aligned}
\]
\endgroup
The last equality uses $D^{-1}(HD^{-1})D=D^{-1}H$. The remaining
diagonal block contributes $(r+b)^d$. Both sides are polynomials,
so the identity also holds at $r=0$.

The matrix $D^{-1}H$ is similar to the real symmetric matrix
		$D^{-1/2}HD^{-1/2}$. For a real eigenvalue $\mu$ of that matrix, the
		corresponding two eigenvalues of $J$ solve $r^2+a r+a\mu=0$. If $\mu>0$
		both have negative real part; if $\mu<0$ they are real and have opposite
		signs; and if $\mu=0$ they are $0,-a$. Thus zero is the only possible
		imaginary-axis eigenvalue of $J$. Its algebraic multiplicity is $k$,
		and solving the kernel equations gives
\[
E^c(q):=\ker J(q)
=\{(u,0,\mathrm D S(q)u):u\in\ker H(q)\}.
\]
Its geometric multiplicity is also $k$, so its Jordan blocks at zero
are all of size one.

The real parts of all nonzero eigenvalues of $J(q)$ stay a uniform
positive distance from zero as $q$ varies in $\mathcal M$. Indeed, the ordered
		eigenvalues of $D^{-1/2}HD^{-1/2}$ are continuous, and constant nullity
		prevents a nonzero eigenvalue from tending to zero on this compact set.
		The corresponding quadratic roots stay away from the imaginary axis,
		since $r=i\nu$ would give $a\nu=0$ and then $\mu=0$.
		Together with $-a,-b$, this gives a common gap for $J$.
		Since the critical component is connected, the numbers of eigenvalues
with negative and positive real parts are constant, counted with
algebraic multiplicity.
		By continuity of $J(q)$ and the uniform spectral gap, the stable,
center, and unstable spectral projections $P_s(q),P_c(q),P_u(q)$
depend continuously on $q$; see
\citet[Chapter~II, \S5.1, Theorem~5.1]{Kato1995}, applied locally to
the distinct eigenvalues and with the corresponding projections summed.
These projections extract the three components in the direct sum
$\mathbb R^{3d}=E^s(q)\oplus E^c(q)\oplus E^u(q)$; they need not be
orthogonal projections. Their ranges
\[
E^j(q):=\operatorname{Ran}P_j(q),\qquad
\dim E^j(q)=\operatorname{rank}P_j(q)=\text{constant},
\quad j\in\{s,c,u\},
\]
therefore form continuous vector bundles.
The finite-dimensional exponential estimates on the stable and unstable
subspaces are uniform on the compact set $\mathcal M$: for some
$C,\sigma>0$,
\[
\|e^{tJ(q)}P_s(q)\|\le Ce^{-\sigma t},\qquad
\|e^{-tJ(q)}P_u(q)\|\le Ce^{-\sigma t}
\qquad(q\in\mathcal M,\ t\ge0).
\]
On $E^c(q)=\ker J(q)$ the linear flow is the identity.
For an absent stable or unstable subspace, the corresponding projection
is zero.

\noindent\textbf{Step 2: The parameter extension.}
		Extend $f_\lambda$ in~\eqref{eq:auto} to negative $\lambda$ and consider
\[
\begin{cases}
\dot z=f_\lambda(z),\\
\dot\lambda=0.
\end{cases}
\]
The compact connected manifold $Y=\mathcal N_{\mathcal M}\times\{0\}$ consists of equilibria
		of this extended field. Its linearization is
		\[
		\widehat J(q)=\begin{pmatrix}J(q)&h(q)\\0&0\end{pmatrix}
		\quad\text{\normalfont and}\quad h(q)=(-q,0,0).
		\]
		Let $P_c$ be the center projection, $P_h=I-P_c$, and
$J_h=J|_{E^s\oplus E^u}$. For a vector $(v+w,\tau)$ with
$v\in E^c$ and $w\in E^s\oplus E^u$, the component in $E^s\oplus E^u$
after one application of $\widehat J(q)$ is $J_hw+P_hh(q)\tau$.
Any further applications act on this component by $J_h$. Since
$J_h$ is invertible, the component must vanish for the vector to
belong to the generalized eigenspace for the zero eigenvalue of
$\widehat J(q)$. Solving for $w$ gives
\begin{equation}\label{eq:extended-center-space}
			\widehat E^c(q)=
			\{(v-J_h^{-1}P_hh(q)\tau,\tau):v\in E^c(q),\ \tau\in\mathbb R\}.
		\end{equation}
		For $v\in E^c(q)$ and $\tau\in\mathbb R$, we have
\[
\begin{aligned}
\widehat J(q)\binom{v-J_h^{-1}P_hh(q)\tau}{\tau}
&=\binom{-P_hh(q)\tau+h(q)\tau}{0}
=\binom{P_ch(q)\tau}{0},\\
\widehat J(q)^2\binom{v-J_h^{-1}P_hh(q)\tau}{\tau}
&=\binom{J(q)P_ch(q)\tau}{0}=0.
\end{aligned}
\]
The parametrization in~\eqref{eq:extended-center-space} is injective, so
\[
\dim\widehat E^c(q)=\dim E^c(q)+1=k+1,\qquad
\widehat E^s(q)=E^s(q)\times\{0\},\quad
\widehat E^u(q)=E^u(q)\times\{0\}.
\]
If $E^s(q)\oplus E^u(q)=\{0\}$, the correction term is zero.
These bundles are continuous, have constant ranks, and satisfy
$T_{(z_q,0)}Y\subset\widehat E^c(q)$.
		
		\noindent\textbf{Step 3: The invariant manifold and its slices.}

We apply \hyperref[app:fenichel-center-theorem]{the Fenichel theorem stated above}
with $X_\lambda=f_\lambda$ and
$K_*=\mathcal N_{\mathcal M}$.
		The family extends to positive and negative $\lambda$ near zero.
		For $v\in T_q\mathcal M$, differentiating $\nabla F|_{\mathcal M}=0$
gives $H(q)v=0$. Thus
\[
T_{z_q}\mathcal N_{\mathcal M}
=\{(v,0,\mathrm D S(q)v):v\in T_q\mathcal M\}\subset\ker J(q).
\]
Put $\ell=\dim\mathcal M$. Since this tangent space lies in
$E^c(q)=\ker J(q)$, passing to the quotient removes only $\ell$
zero-eigenvalue directions. The stable and unstable subspaces retain
their dimensions, and
\[
\dim\bigl(E^c(q)/T_{z_q}\mathcal N_{\mathcal M}\bigr)=k-\ell.
\]
For the extended system, $\dim\widehat E^c(q)=k+1$ and
$\dim T_{(z_q,0)}Y=\ell$, so
\[
\dim\bigl(\widehat E^c(q)/T_{(z_q,0)}Y\bigr)=k+1-\ell,
\qquad \ell+(k+1-\ell)=k+1.
\]
The stable and unstable eigenvalues remain separated from the imaginary
axis by the uniform bound proved above.
The theorem gives a $(k+1)$-dimensional $C^1$
		embedded manifold $\widehat{\mathcal I}$ containing $Y$, tangent to
		$\widehat E^c$ there, and locally invariant under the extended flow.

Let $\pi(z,\mu)=\mu$. Taking $v=0$ and $\tau=1$
in~\eqref{eq:extended-center-space} gives
\[
\binom{-J_h^{-1}P_hh(q)}{1}\in T_{(z_q,0)}\widehat{\mathcal I},
\qquad
D\pi_{(z_q,0)}\binom{-J_h^{-1}P_hh(q)}{1}=1.
\]
Thus $\operatorname{rank}D(\pi|_{\widehat{\mathcal I}})=1$ on $Y$.
By continuity and compactness, the rank remains one in a neighborhood
of $Y$. The submersion theorem gives, at each point of a level there,
\[
T(\widehat{\mathcal I}\cap\{\pi=\lambda\})
=\ker\bigl(D\pi|_{T\widehat{\mathcal I}}\bigr),
\qquad
\dim\ker\bigl(D\pi|_{T\widehat{\mathcal I}}\bigr)
=(k+1)-1=k.
\]
The slices
		\[
		\mathcal I_{\mathcal M,\lambda}
		=\{z:(z,\lambda)\in\widehat{\mathcal I}\}
		\]
		in this neighborhood are therefore $k$-dimensional $C^1$ embedded
		manifolds. By the submersion theorem, local coordinates on $\widehat{\mathcal I}$
near any point of $Y$ can be chosen with $\lambda$ as one coordinate.
Thus these slices are nonempty for every sufficiently small $|\lambda|$.
		Since $\widehat{\mathcal I}$ is locally invariant and $\dot\lambda=0$,
each slice is locally invariant under $f_\lambda$.
At $\lambda=0$, the slice contains $\mathcal N_{\mathcal M}$ and its tangent spaces are
		$\widehat E^c\cap\{\tau=0\}=E^c\times\{0\}$. Thus
\[
\mathcal N_{\mathcal M}\subset\mathcal I_{\mathcal M,0}\quad\text{\normalfont and}\quad
T_{(q,0,S(q))}\mathcal I_{\mathcal M,0}
=\{(v,0,\mathrm D S(q)v):v\in\ker\nabla^2F(q)\}.
\]

The statement about complete orbits in Fenichel's theorem gives a neighborhood
$\widehat V$ of $Y$. Intersect it with the ambient
		neighborhood used for the submersion restriction above. Every complete
		extended orbit whose closure lies in $\widehat V$ is then contained in
		the restricted $\widehat{\mathcal I}$. Compactness allows us to choose
		$V$ with compact closure in $U$ and
		$\lambda_{\mathcal M}>0$ such that
		\[
		\overline V\times[-\lambda_{\mathcal M},\lambda_{\mathcal M}]
		\subset\widehat V,
		\]
		inside the submersion neighborhood and the interval of nonempty
		slices. For $|\lambda|<\lambda_{\mathcal M}$ and a solution
$z:\mathbb R\to V$ of $\dot z=f_\lambda(z)$,
\[
\frac{d}{dt}(z(t),\lambda)=(f_\lambda(z(t)),0),\qquad
\overline{\{(z(t),\lambda):t\in\mathbb R\}}
\subset\overline V\times\{\lambda\}\subset\widehat V.
\]
Fenichel's complete-orbit conclusion therefore gives
\[
\{(z(t),\lambda):t\in\mathbb R\}\subset\widehat{\mathcal I}
\quad\Longrightarrow\quad
z(\mathbb R)\subset\{z:(z,\lambda)\in\widehat{\mathcal I}\}
=\mathcal I_{\mathcal M,\lambda}.
\]
This proves the lemma.
	\end{proof}

\subsection{Proof of Lemma~\ref{lem:approximate-trajectory-convergence}}\label{guide:appendix-7-5}
\label{app:approximate-trajectory-convergence}
Steps 1--4 prove local convergence; Steps 5--6 show that the path
eventually remains in one of these neighborhoods.

\begin{proof}
\noindent\textbf{Step 1: Local coordinates.}
Fix a critical point $q$. By hypothesis,
$\dim\ker\nabla^2F(q)\le1$. Write $e^0=(q,0,S(q))$.
In translated coordinates, put
\[
X_\lambda(z):=f_\lambda(e^0+z),\quad J:=DX_0(0).
\]
The field is jointly $C^1$ near $(0,0)$ and $X_0(0)=0$.
By the spectral calculation in Appendix~\ref{app:center-reduction},
Step~1, the only possible imaginary-axis
eigenvalue of $J$ is zero, with algebraic multiplicity at most one.
We construct a neighborhood $V$ of zero and a parameter threshold,
independent of the path, such that a locally absolutely continuous path
eventually staying in $V$ and satisfying almost everywhere
\[
\dot w=(I+C(t))X_\lambda(w)+r(t),\quad
C(t)\to0,\quad \int^\infty\|r(t)\|\,dt<\infty,
\]
has finite length and converges to an equilibrium.
Here $w$ denotes the translated path, and $C,r$ are measurable.

If zero is an eigenvalue, its algebraic multiplicity is one, so its
generalized eigenspace equals $\ker J$; this is what semisimple means
here. We apply
\hyperref[app:fenichel-center-theorem]{Fenichel's theorem as stated in Appendix~\ref{app:center-reduction}} with
both the equilibrium submanifold of $X_0$ and its compact subset equal
to $\{0\}$. The extended compact set is $\{(0,0)\}$, so the required numbers of
eigenvalues are constant on it. This gives a $C^1$ local center manifold
tangent to the generalized eigenspace for the zero eigenvalue of the
extended linearization, computed below.
Write $E^c=\ker J$. Let $E^s$ and $E^u$ be the real invariant
subspaces associated with the generalized eigenspaces of $J$ whose
eigenvalues have negative and positive real parts, respectively. Set
\[
E^h:=E^s\oplus E^u,\qquad \mathbb R^{3d}=E^c\oplus E^h.
\]
Thus $E^h$, called the hyperbolic spectral subspace, contains all
directions associated with eigenvalues of nonzero real part:
contracting directions in $E^s$ and expanding directions in $E^u$.
Let $P_c,P_h$ be the projections onto $E^c,E^h$. In these coordinates,
we use the product Euclidean norm, so the coordinate projections have
norm at most one. Since zero is not an eigenvalue on $E^h$,
$A:=J|_{E^h}$ is invertible. Put
$h_q:=\left.\partial_\lambda X_\lambda(0)\right|_{\lambda=0}$.
		The extended generalized zero eigenspace is
		\[
		\{(u-A^{-1}P_hh_q\,\tau,\tau):u\in E^c,\ \tau\in\R\}.
		\]
		Indeed, the extended linearization sends this vector to
		$(P_ch_q\,\tau,0)$ and then to zero. Projection onto $(u,\lambda)$ is an isomorphism on this space.
By the inverse function theorem, the center manifold has the form
$v=h(u,\lambda)$, with $h\in C^1$ and $D_uh(0,0)=0$.
		For fixed $\lambda$, set
		$\psi(u)=P_cX_\lambda(u,h(u,\lambda))$. Since the vector field is tangent to the graph $v=h(u,\lambda)$, we have
		\begin{equation}\label{eq:one-center-tangency}
		X_\lambda(u,h(u,\lambda))=(\psi(u),D_uh(u,\lambda)\psi(u)).
		\end{equation}
		When $E^c=\{0\}$ the slice is a single equilibrium and the scalar
		terms below are omitted. Empty hyperbolic blocks are omitted as well.
		
		Fix $M,\sigma>0$ from the hyperbolic block $A$ so that, for $t\ge0$,
\[
\|e^{tA}P_s\|\le Me^{-\sigma t},
\qquad
\|e^{-tA}P_u\|\le Me^{-\sigma t}.
\]
Thus $C_H:=2M/\sigma$ depends only on $A$ and can be fixed before
choosing the coordinate neighborhood. Choose small balls in $E^c$ and
$E^h$ and a neighborhood $V$ whose closure lies in their product.
Shrinking these balls and the parameter interval makes
		$DX_\lambda-J$ and $D_uh$ uniformly small. With
		$\eta=v-h(u,\lambda)$, the fundamental theorem of calculus along
$(u,h(u,\lambda)+\tau\eta)$, $0\le\tau\le1$, gives
\begingroup
\[
X_\lambda(u,v)-X_\lambda(u,h(u,\lambda))-J(0,\eta)
=\int_0^1
\bigl[DX_\lambda(u,h(u,\lambda)+\tau \eta)-J\bigr](0,\eta)\,d\tau.
\]
\endgroup
The right-hand side has norm at most $\kappa\|\eta\|$. Since
$J(0,\eta)=(0,A\eta)$, the two coordinate projections give
\begingroup
\begin{equation}\label{eq:one-center-vertical-estimates}
P_cX_\lambda(u,v)=\psi(u)+O(\kappa\|\eta\|),\qquad
P_hX_\lambda(u,v)=D_uh\,\psi(u)+A\eta+O(\kappa\|\eta\|),
\end{equation}
\endgroup
		where both error norms are at most $\kappa\|\eta\|$ and
		$\|D_uh\|\le\kappa\le1$. The fixed change of linear coordinates
		preserves $C(t)\to0$ and integrability of $r$; we use the same notation.
		Take the initial time late enough that $\|C(t)\|\le\epsilon\le1$.
        Put $L_A=\|A\|+2$ and $L_0=1+L_A$. The chain rule gives
\begin{equation}\label{eq:one-center-relative-coordinates}
\dot u=(1+c(t))\psi(u)+\rho(t)\qquad\text{\normalfont and}\qquad
\dot \eta=A\eta+\zeta(t),
\end{equation}
where the perturbation terms satisfy
\begingroup
\[
|c(t)|\le2\epsilon,\qquad
|\rho(t)|\le L_0\|\eta\|+\|r\|,\qquad
\|\zeta(t)\|\le(2\kappa+2\epsilon L_A)\|\eta\|
+4\epsilon|\psi(u)|+2\|r\|.
\]
\endgroup
		Indeed, $\dot \eta=\dot v-D_uh\,\dot u$ gives
\[
\dot \eta=(P_h-D_uh P_c)\bigl((I+C)X_\lambda(u,v)+r\bigr).
\]
By~\eqref{eq:one-center-tangency}--\eqref{eq:one-center-vertical-estimates},
\begingroup
\[
\|X_\lambda(u,h(u,\lambda))\|\le2|\psi(u)|,\qquad
\|X_\lambda(u,v)-X_\lambda(u,h(u,\lambda))\|\le L_A\|\eta\|.
\]
\endgroup
The three terms in $\zeta$ satisfy
\begin{align*}
\|(P_h-D_uh P_c)X_\lambda(u,v)-A\eta\|&\le(\kappa+\kappa^2)\|\eta\|\le2\kappa\|\eta\|,\\
\|(P_h-D_uh P_c)C X_\lambda(u,v)\|&\le2\epsilon\bigl(2|\psi(u)|+L_A\|\eta\|\bigr),\\
\|(P_h-D_uh P_c)r\|&\le2\|r\|.
\end{align*}
Hence the estimate on $\zeta$ holds. For the scalar equation, with
$u=u(t)$, set
\[
c(t):=P_c C(t)(1,D_uh(u,\lambda)).
\]
By~\eqref{eq:one-center-tangency},
\[
P_c C(t)X_\lambda(u,h(u,\lambda))=c(t)\psi(u),
\qquad
|c(t)|\le\epsilon(1+\|D_uh(u,\lambda)\|)\le2\epsilon.
\]
The remaining scalar terms are
\[
\begin{aligned}
\rho(t):={}&P_cX_\lambda(u,v)-\psi(u)\\
&+P_cC(t)\bigl[X_\lambda(u,v)-X_\lambda(u,h(u,\lambda))\bigr]
+P_cr(t).
\end{aligned}
\]
Thus
\[
|\rho|\le(\kappa+\epsilon L_A)\|\eta\|+\|r\| \le L_0\|\eta\|+\|r\|.
\]
		
		\noindent\textbf{Step 2: An integral estimate for $\eta$.}
On any sufficiently late finite interval $[T,T']$, integrate the stable
equation forward from $T$ and the unstable equation backward from $T'$,
as illustrated in Figure~\ref{fig:one-center-integration}:
		\[
		\begin{gathered}
			P_s\eta(t)=e^{A(t-T)}P_s\eta(T)+\int_T^t e^{A(t-s)}P_s\zeta(s)\,ds,\\
			P_u\eta(t)=e^{-A(T'-t)}P_u\eta(T')-\int_t^{T'} e^{-A(s-t)}P_u\zeta(s)\,ds.
		\end{gathered}
		\]

\begin{figure}[htbp]
\centering
\includegraphics[width=.63\linewidth]{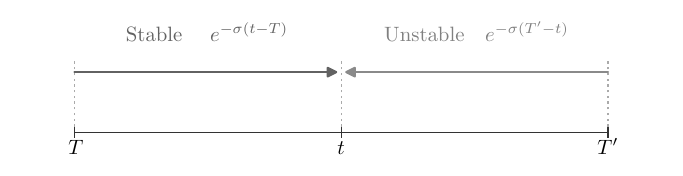}
\caption{Endpoint contributions to the deviation $\eta=v-h(u,\lambda)$ at $t\in[T,T']$.
The stable component is integrated forward from $T$, and the unstable
component backward from $T'$. The labeled exponential factors, multiplied by $M$, bound the norms
of the linear maps acting on the respective endpoint values.}
\label{fig:one-center-integration}
\end{figure}

The endpoint terms have integrals bounded independently of $T'$: if
		$\|\eta(T)\|,\|\eta(T')\|\le C_V$ on this coordinate neighborhood, then each contributes at most
		$MC_V\int_0^{T'-T}e^{-\sigma s}\,ds\le MC_V/\sigma$.
		For the stable convolution, Fubini gives
		\begingroup
\[
\int_T^{T'}\int_T^t Me^{-\sigma(t-s)}\|\zeta(s)\|\,ds\,dt
=M\int_T^{T'}\|\zeta(s)\|\int_s^{T'} e^{-\sigma(t-s)}\,dt\,ds
\le\frac M\sigma\int_T^{T'}\|\zeta(s)\|\,ds;
\]
\endgroup
		the backward unstable convolution satisfies
\[
\int_T^{T'}\int_t^{T'} Me^{-\sigma(s-t)}\|\zeta(s)\|\,ds\,dt
=\int_T^{T'}\|\zeta(s)\|\int_T^s Me^{-\sigma(s-t)}\,dt\,ds
\le \frac M\sigma\int_T^{T'}\|\zeta(s)\|\,ds.
\]
The two endpoint terms sum to at most $C_0:=2MC_V/\sigma$.
The convolution terms sum to at most $C_H\int_T^{T'}\|\zeta\|$.

\noindent\textbf{Step 3: The center integral.}
To estimate the scalar motion, fix any point $s_0$ in the scalar
coordinate interval and define the Lipschitz function
$\Psi(s)=\int_{s_0}^s\operatorname{sgn}(\psi(v))\,dv$, with
$\operatorname{sgn}(0)=0$. Since $|\operatorname{sgn}(\psi)|\le1$,
\[
|\Psi(u(T'))-\Psi(u(T))|
\le |u(T')-u(T)|\le\operatorname{diam}V,
\]
where the diameter is measured in the product coordinates chosen above.
Where $\psi(u)\ne0$, its sign is locally constant, so
\begingroup
\[
(\Psi\circ u)'=\operatorname{sgn}(\psi(u))\bigl((1+c(t))\psi(u)+\rho(t)\bigr)
\ge(1-2\epsilon)|\psi(u)|-|\rho|.
\]
\endgroup
		Where $\psi(u)=0$, the equation gives $\dot u=\rho$, and the
		Lipschitz chain rule gives $|(\Psi\circ u)'|\le|\rho|$ almost
		everywhere. The same lower bound therefore holds throughout.
		Integrating the lower bound and using the estimate on $\rho$ gives
\begingroup
\[
(1-2\epsilon)\int_T^{T'}|\psi(u(t))|\,dt
\le \Psi(u(T'))-\Psi(u(T))+\int_T^{T'}|\rho(t)|\,dt
\le\operatorname{diam}V+L_0\int_T^{T'}\|\eta\|
+\int_T^{T'}\|r\|.
\]
\endgroup
For $\epsilon\le1/4$, we have $1-2\epsilon\ge1/2$.

\noindent\textbf{Step 4: Finite length and convergence.}
Shrink the coordinate neighborhood and parameter interval so that $2C_H\kappa\le1/4$.
Choose $T$ late enough that $\|C(t)\|\le\epsilon$ for $t\ge T$, with
$\epsilon\le1/4$ and $C_H\epsilon(2L_A+8L_0)\le1/4$. Set
\[
R_{T'}=\int_T^{T'}\|\eta\|,\qquad P_{T'}=\int_T^{T'}|\psi(u)|,\qquad
E_{T'}=\int_T^{T'}\|r\|.
\]
The bounds from Steps 2 and 3, with the estimate for $\zeta$, give
\begin{align}
R_{T'} &\le C_0+C_H(2\kappa+2\epsilon L_A)R_{T'}
          +4C_H\epsilon P_{T'}+2C_HE_{T'},\label{eq:one-center-normal-integral}\\
P_{T'} &\le 2\operatorname{diam}V+2L_0R_{T'}+2E_{T'}.\label{eq:one-center-scalar-integral}
\end{align}
Substituting the second bound into the first gives
\[
\tfrac12 R_{T'}
\le \bigl[1-2C_H\kappa-C_H\epsilon(2L_A+8L_0)\bigr]R_{T'}
\le C_0+8C_H\epsilon\operatorname{diam}V+(2C_H+8C_H\epsilon)E_{T'}.
\]
Since $E_{T'}\le\int_T^\infty\|r\|$, both $R_{T'}$ and $P_{T'}$ are bounded
uniformly in $T'$. If either block is absent, set its integral to zero
and use only the other inequality.

Letting $T'\to\infty$ gives $\eta,\psi(u)\in L^1$. The coordinate equations give
\[
 |\dot u|\le(1+2\epsilon)|\psi(u)|+L_0\|\eta\|+\|r\|
 \quad\text{\normalfont and}\quad
 \|\dot \eta\|\le C\bigl(\|\eta\|+|\psi(u)|+\|r\|\bigr).
\]
Both velocities are integrable. Since $\dot v=D_uh\,\dot u+\dot \eta$
and $D_uh$ is bounded, the fixed coordinate change gives
\[
\int_T^\infty\|\dot w\|\,dt
\le C\int_T^\infty\bigl(|\dot u|+\|\dot \eta\|\bigr)\,dt <\infty.
\]
Moreover,
\[
\sup_{s\ge t}\|w(s)-w(t)\|
\le\int_t^\infty\|\dot w\|\,ds\longrightarrow0.
\]
Thus $w$ converges, and its coordinates have limits $u_*,\eta_*$.
Continuity gives $\psi(u(t))\to\psi(u_*)$. A function with a nonzero
limit cannot be integrable:
\[
f(t)\to f_*\ne0
\quad\Longrightarrow\quad
\|f(t)\|\ge\tfrac12\|f_*\|\quad(t\ge T_*)
\quad\Longrightarrow\quad
\int_{T_*}^\infty\|f(t)\|\,dt=\infty.
\]
Applying this to $\eta$ and, when the center is one-dimensional,
$\psi(u)$ gives $\eta_*=0$ and $\psi(u_*)=0$.
Equation~\eqref{eq:one-center-tangency} then gives
$X_\lambda(u_*,h(u_*,\lambda))=0$, so the limit is an equilibrium. The limit lies in $\overline V$, which is contained in the domain of $X_\lambda$.

\noindent\textbf{Step 5: Flow tracking.}
Return to the original coordinates and fix a path $w$ satisfying
\eqref{eq:deterministic-relative-perturbation} and the other hypotheses
of the lemma.

The condition $\operatorname{dist}(w(t),K)\to0$ implies that $w$
eventually lies in a compact neighborhood of $K$ contained in $U$.
Together with continuity on the initial finite interval, this gives
compact closure in $U$ and $L:=\omega(w)\subset K$.
Let $M$ bound $\|f_\lambda\|$ on that compact neighborhood. By
\eqref{eq:deterministic-relative-perturbation}, for every fixed $H>0$,
\[
\begin{aligned}
\int_t^{t+H}\|\dot w(u)-f_\lambda(w(u))\|\,du
&\le\int_t^{t+H}\bigl(\|C(u)\|\,\|f_\lambda(w(u))\|+\|r(u)\|\bigr)\,du\\
&\le HM\sup_{u\ge t}\|C(u)\|+\int_t^\infty\|r(u)\|\,du
\longrightarrow0.
\end{aligned}
\]
As in the proof of Lemma~\ref{lem:corrected-tracking}, we extend $f_\lambda$
by a cutoff equal to one on a larger neighborhood. Write
$\widetilde\Phi_\lambda$ for the extended flow and let $L_*>0$ be a
Lipschitz constant for its field. Gronwall's inequality gives
\[
\sup_{0\le s\le H}\|w(t+s)-\widetilde\Phi_\lambda(s,w(t))\|
\le e^{L_*H}\left(HM\sup_{u\ge t}\|C(u)\|
                  +\int_t^\infty\|r(u)\|\,du\right)\longrightarrow0.
\]
For large $t$, the comparison orbit remains where the cutoff equals
one, so the estimate holds for the original flow. This proves
\eqref{eq:deterministic-flow-tracking}.

\noindent\textbf{Step 6: Localization and convergence.}
The critical set is closed and consists of isolated points, so it is
locally finite and countable. Thus the critical values have empty
interior, and Lemma~\ref{lem:compact-limit-geometry} gives
\[
\max_{z\in L}\operatorname{dist}(z,\mathcal E_0)
=\limsup_{t\to\infty}\operatorname{dist}(w(t),\mathcal E_0) \le\rho_K(\lambda).
\]
As recalled in Appendix~\ref{sec:compact-limit-geometry}, $L$ is
nonempty and connected, and satisfies~\eqref{eq:limit-set-invariance}.
If $K\cap\mathcal E_0=\varnothing$, then
$\operatorname{dist}(K,\mathcal E_0)>0$, contradicting this estimate
for sufficiently small $\lambda$. Hence no such path exists in that case.
Otherwise, write
\[
K\cap\mathcal E_0=\{e_1^0,\ldots,e_N^0\}\quad\text{\normalfont and}\quad
 e_j^0=(q_j,0,S(q_j)).
\]
Apply Steps 1--4 at each $q_j$ and let $V_j$ be the resulting
neighborhoods in the original coordinates. Choose disjoint open balls
$B_j$ about $e_j^0$ with $\overline{B_j}\subset V_j$.
Put $D_K=K\setminus\bigcup_j B_j$. If $D_K$ is nonempty, it is
compact and disjoint from $\mathcal E_0$, so for small enough $\lambda$,
\[
\max_{z\in L}\operatorname{dist}(z,\mathcal E_0)
\le\rho_K(\lambda)<\operatorname{dist}(D_K,\mathcal E_0).
\]
Hence $L\subset\bigcup_j B_j$, also immediate when $D_K$ is empty.
Connectedness puts $L$ in one $B_j$.
Choose $\lambda_K$ below this threshold and all the finitely many local
thresholds; these choices depend on $K$, not on $w$.
Since the range of $w$ has compact closure in $U$,
$\operatorname{dist}(w(t),L)\to0$, so
$\overline{B_j}\subset V_j$ implies that $w(t)\in V_j$ for all
sufficiently large $t$.

Step 4, applied in $V_j$, gives finite length after some time and
convergence to an equilibrium. Since $\operatorname{dist}(w(t),K)\to0$
and $K$ is closed, the limit belongs to $K$.
\end{proof}

\subsection{Proof of Lemma~\ref{lem:quantitative-energy}}\label{guide:appendix-7-6}
\label{app:quantitative-proof}
We differentiate the energy terms along time, using local Lipschitz
bounds for $\nabla F$ and $S$.

\begin{proof}
Fix a compact set $K\subset U$ and $0\le\lambda\le1$. We estimate
solutions of $\dot z=f_\lambda(z)$ that remain in $K$. Recall the residual
\[
R(x,m,v)^2=\|\nabla F(x)\|^2+\|m\|^2+\|S(x)-v\|^2.
\]
For a sufficiently small $\delta>0$, to be chosen below, set
  \begin{equation}\label{eq:quantitative-energy-function}
    W=F(x)+\frac1{2a}\sum_i\frac{m_i^2}{\varepsilon+\sqrt{v_i}}
       -\delta\langle \nabla F(x),m\rangle+\delta\|S(x)-v\|^2.
  \end{equation}

  Enlarge $K$ slightly within $U$ and take a compact ball containing
  its position projection. The maps $\nabla F,S$ have Lipschitz constants
  $L_{\nabla F},L_S$ on this ball. Write
  \[
    d_i=\varepsilon+\sqrt{v_i},\quad e=S(x)-v,\quad\text{\normalfont and}\quad
    V=F(x)+\frac1{2a}\sum_i\frac{m_i^2}{d_i}.
  \]
  On this compact region, the $d_i$ are bounded above and away from zero,
and the $v_i$ are positive. Since the reciprocal denominators are smooth
on $U$ and $\nabla F,S$ are locally Lipschitz, $W$ is locally Lipschitz.
  Here $\dot d_i=b(S_i-v_i)/(2\sqrt{v_i})$, so
  \[
    \frac{d}{dt}\frac{m_i^2}{2ad_i}
       =\frac{m_i(\partial_iF-m_i)}{d_i}
        -\frac b{4a}\frac{m_i^2(S_i-v_i)}{\sqrt{v_i}d_i^2}.
  \]
  Adding $\dot F=-\lambda\langle x,\nabla F(x)\rangle-\sum_i\partial_iF(x)m_i/d_i$
  cancels the mixed terms and gives
  \[
\dot V
      =-\lambda\langle x,\nabla F(x)\rangle-\sum_i\frac{m_i^2}{d_i}
        +\frac b{4a}\sum_i
          \frac{m_i^2(v_i-S_i(x))}{\sqrt{v_i}d_i^2} \le -\sum_i\frac{\varepsilon m_i^2}{d_i^2}
                  +\lambda M_x\|\nabla F(x)\|,
\]
  where $M_x$ bounds $\|x\|$ on the enlarged compact region.
  Here $S_i>0$ and $b\le4a$ give
\[
-\frac{m_i^2}{d_i}
+\frac b{4a}\frac{m_i^2(v_i-S_i)}{\sqrt{v_i}d_i^2}
\le-\frac{m_i^2}{d_i^2}
 \left[\varepsilon+\left(1-\frac b{4a}\right)\sqrt{v_i}\right] \le-\frac{\varepsilon m_i^2}{d_i^2}.
\]
  If $d_{\max}$ bounds all $d_i$, put $\mu=\varepsilon/d_{\max}^2>0$.

  For $G=\nabla F$ or $G=S$ and $s<t$ along the orbit,
\[
\|G(x(t))-G(x(s))\|
\le L_G\|x(t)-x(s)\|
\le L_G\int_s^t\|\dot x(r)\|\,dr .
\]
Thus both compositions are locally absolutely continuous. Almost everywhere,
\[
\left\|\frac d{dt}\nabla F(x(t))\right\|
 \le L_{\nabla F}\|\dot x\|,\quad
\left\|\frac d{dt}S(x(t))\right\|\le L_S\|\dot x\|,\quad\text{\normalfont and}\quad
\|\dot x\|\le\varepsilon^{-1}\|m\|+\lambda M_x.
\]
Using $\dot m=a(\nabla F-m)$ and Young's inequality, we obtain
\[
\begin{aligned}
\frac d{dt}(-\langle\nabla F,m\rangle)
&=-\left\langle\frac d{dt}\nabla F(x(t)),m\right\rangle
  -a\|\nabla F\|^2+a\langle\nabla F,m\rangle\\
&\le L_{\nabla F}(\varepsilon^{-1}\|m\|+\lambda M_x)\|m\|
 -a\|\nabla F\|^2+\frac a2(\|\nabla F\|^2+\|m\|^2)\\
&=-\frac a2\|\nabla F\|^2
 +\left(\frac{L_{\nabla F}}\varepsilon+\frac a2\right)\|m\|^2
 +\lambda L_{\nabla F}M_x\|m\|.
\end{aligned}
\]
Since $\dot v=be$, we have $\dot e+be=\frac{d}{dt}S(x(t))$
almost everywhere. Taking the inner product with $2e$ and using
$2L_S\|e\|\|m\|/\varepsilon
\le b\|e\|^2+L_S^2\|m\|^2/(b\varepsilon^2)$ gives
  \[
\frac{d}{dt}\|e\|^2
    \le-2b\|e\|^2
      +2L_S\|e\|\left(\varepsilon^{-1}\|m\|+\lambda M_x\right) \le-b\|e\|^2+\frac{L_S^2}{b\varepsilon^2}\|m\|^2
      +2\lambda L_SM_x\|e\|.
\]

  Set
  \[
    C_m=\frac{L_{\nabla F}}{\varepsilon}+\frac a2+
                   \frac{L_S^2}{b\varepsilon^2},\quad
    0<\delta\le\frac{\mu}{2C_m},\quad\text{\normalfont and}\quad
    c_0=\min\{\mu/2,\delta a/2,\delta b\}.
  \]
  Adding the estimates gives
\[
\begin{aligned}
\dot W
&\le-(\mu-\delta C_m)\|m\|^2
-\frac{\delta a}{2}\|\nabla F\|^2-\delta b\|e\|^2
+\lambda M_x\bigl(\|\nabla F\|+\delta L_{\nabla F}\|m\|
+2\delta L_S\|e\|\bigr)\\
&\le-c_0R^2
+\lambda M_x\bigl(\|\nabla F\|+\delta L_{\nabla F}\|m\|
+2\delta L_S\|e\|\bigr),
\end{aligned}
\]
where $\mu-\delta C_m\ge\mu/2$ by the choice of $\delta$.
  Cauchy--Schwarz gives
  \[
    \|\nabla F\|+\delta L_{\nabla F}\|m\|+2\delta L_S\|e\|
    \le (1+\delta^2L_{\nabla F}^2+4\delta^2L_S^2)^{1/2}R.
  \]
  Thus the first inequality in~\eqref{eq:quantitative-energy-estimates}
  holds when $B_0\ge M_x(1+\delta^2L_{\nabla F}^2+4\delta^2L_S^2)^{1/2}$. The explicit formula also gives
  \[
|W-F|\le\frac{\|m\|^2}{2a\varepsilon}
          +\frac\delta2(\|\nabla F\|^2+\|m\|^2)+\delta\|e\|^2 \le\left(\frac1{2a\varepsilon}+\delta\right)R^2,
\]
  which proves the second with the same sufficiently large $B_0$.

  Finally, let $L\subset K$ be nonempty and compact, containing the
complete orbit through each of its points. Suppose
$R(z)>B_0\lambda/c_0$ at a maximum or minimum $z\in L$ of $W$. By continuity of $R$, on a sufficiently short orbit
segment through $z$ we have, almost everywhere,
\[
\dot W\le-c_0R^2+B_0\lambda R =-R(c_0R-B_0\lambda) <0.
\]
Since $\Phi_\lambda(\pm h,z)\in L$, maximality or minimality,
respectively, would give
\[
\begin{aligned}
0&\le W(z)-W(\Phi_\lambda(-h,z))
=\int_{-h}^0\dot W(\Phi_\lambda(t,z))\,dt<0,\\
0&\le W(\Phi_\lambda(h,z))-W(z)
=\int_0^h\dot W(\Phi_\lambda(t,z))\,dt<0.
\end{aligned}
\]
Thus $R(z)\le B_0\lambda/c_0$ at both extrema.
\end{proof}

\subsection{Proof of Lemma~\ref{lem:quantitative-band}}\label{guide:appendix-7-7}
\label{app:quantitative-band}
We estimate the energy at its maximum and minimum on the limit set.

\begin{proof}
Write $\pi_x(x,m,v)=x$. Fix a nonempty compact $K\subset U$ and a
closed ball $X$ with $\pi_xK\subset\operatorname{int}X$ and
$X\cap\operatorname{Crit}(F)\ne\varnothing$. We work with
$0<\lambda\le\lambda_K\le1$, choosing $\lambda_K$ sufficiently small
below. Let $w:[0,\infty)\to U$
be continuous, with compact closure of its range in $U$, satisfying
\eqref{eq:deterministic-flow-tracking}, and put $L=\omega(w)\subset K$.
The geometric exponent satisfies $0<\theta\le1/2$.
Lemma~\ref{lem:quantitative-geometry} gives
\[
F(\operatorname{Crit}(F)\cap X)=\{c_1,\ldots,c_J\},
\]
and constants $s_0,k>0$ such that
\[
|F(x)-c_j|\le k\|\nabla F(x)\|^{1/(1-\theta)}
\quad\text{\normalfont for some }j,\qquad
x\in X,\quad\|\nabla F(x)\|\le s_0.
\]
It also shows that $F(\operatorname{Crit}(F))$ is countable and hence
has empty interior. Lemma~\ref{lem:compact-limit-geometry} therefore gives
  \[
    \max_{(x,m,v)\in L}
       \operatorname{dist}\bigl(F(x),\{c_1,\ldots,c_J\}\bigr)
       \longrightarrow0\quad(\lambda\downarrow0),
  \]
  uniformly over all paths satisfying these conditions. Indeed, for
$z=(x,m,v)\in L$ choose a nearest $(q,0,S(q))\in\mathcal E_0$, which
is closed and nonempty. For sufficiently small $\lambda$,
\[
\|x-q\|\le\operatorname{dist}(z,\mathcal E_0)
\le\rho_K(\lambda)<\operatorname{dist}(\pi_xK,X^c).
\]
Thus $q\in X$ and $F(q)\in\{c_1,\ldots,c_J\}$. The distance of $F(x)$
to these values is at most $|F(x)-F(q)|$, which tends uniformly to zero
by uniform continuity of $F$ on $X$.

  Choose disjoint open intervals of one fixed small radius around
  the distinct $c_j$. For every sufficiently small $\lambda$ and
  every such limit set $L$, the connectedness recalled in
  Appendix~\ref{sec:compact-limit-geometry} implies that the image
  $F(\pi_xL)$ lies in one interval, whose center we call
  $c_*$ for this fixed $\lambda$ and $L$. When $J=1$, the choice is immediate.

  Take the energy $W$ and constants $c_0,B_0$ from
Lemma~\ref{lem:quantitative-energy}, and let $z_-,z_+$ minimize and
maximize $W$ on $L$. By~\eqref{eq:limit-set-invariance}, the
solution through each point of $L$ stays in $L$ for all real times. Thus
Lemma~\ref{lem:quantitative-energy} gives
  $R(z_\pm)\le B_0\lambda/c_0$. Reduce $\lambda_K$ so that
$B_0\lambda_K/c_0\le s_0$. Then
$\|\nabla F(\pi_xz_\pm)\|\le R(z_\pm)\le s_0$, so
  Lemma~\ref{lem:quantitative-geometry} gives
$|F(\pi_xz_\pm)-c_j|\le C\lambda^{1/(1-\theta)}$ for some critical value $c_j$.
When $J>1$, take the interval radius less than one third of the
minimum separation between distinct critical values and reduce
$\lambda_K$ so that $C\lambda_K^{1/(1-\theta)}$ is below that radius.
Writing $\eta$ for this radius, at either extremum we have
\[
|c_j-c_*|
\le |c_j-F(\pi_xz_\pm)|+|F(\pi_xz_\pm)-c_*|
<2\eta<\min_{i\ne j}|c_i-c_j|.
\]
Hence $c_j=c_*$.
For $J=1$ these values coincide automatically.
Consequently, for $0<\lambda\le1$ and $1/(1-\theta)\le2$,
\[
|W(z_\pm)-c_*|
\le |W(z_\pm)-F(\pi_xz_\pm)|
      +|F(\pi_xz_\pm)-c_*| \le C\lambda^2+C\lambda^{1/(1-\theta)} \le C'\lambda^{1/(1-\theta)}.
\]
Consequently,
\[
\max_{z\in L}W(z)-\min_{z\in L}W(z)
=W(z_+)-W(z_-)\le2C'\lambda^{1/(1-\theta)},
\]
which proves \eqref{eq:quantitative-band}. The constants and the decay
threshold depend only on $K,X,F,S,a,b,\varepsilon$ and the geometric
bounds; they are independent of $\lambda$ and $w$.
\end{proof}

\subsection{Proof of Lemma~\ref{lem:quantitative-excursion}}\label{guide:appendix-7-8}
\label{app:quantitative-excursion}
We bound the position displacement by the energy decrease while the residual is greater than a threshold.

\begin{proof}
  Fix $0<\lambda\le1$, a compact $K\subset U$, and a nonempty compact
set $L\subset K$ containing the complete $f_\lambda$ orbit through each
of its points. Let $X$ be a closed ball containing the position
projection of $K$ and meeting $\operatorname{Crit}(F)$. Take
$W,c_0,B_0$ from Lemma~\ref{lem:quantitative-energy} on $K$; the residual is
\[
R(x,m,v)^2=\|\nabla F(x)\|^2+\|m\|^2+\|S(x)-v\|^2.
\]
  Cover $\operatorname{Crit}(F)\cap X$ by finitely many neighborhoods
$U_j$ where the assumed distance bound holds with constants $M_j$.
Take $M=\max_j M_j$. If $D_X:=X\setminus\bigcup_j U_j$ is nonempty,
choose
\[
0<s_0<\min_{x\in D_X}\|\nabla F(x)\|.
\]
The minimum is positive by compactness and avoidance of critical points.
If $D_X$ is empty, take any $s_0>0$. Then
  \[
    \operatorname{dist}(x,\operatorname{Crit}(F))\le M\|\nabla F(x)\|^\gamma
    \qquad(x\in X,\ \|\nabla F(x)\|\le s_0).
  \]
  Let $M_x=\max_{(x,m,v)\in K}\|x\|$ and choose
$A\ge\max\{1,2B_0/c_0\}$. Fix a residual threshold $u$ with
$A\lambda\le u\le s_0$, and take a solution $z(t)=(x(t),m(t),v(t))$
defined for all $t\in\mathbb R$ and staying in $L$, with an arbitrary
point of $L$ at time zero. Write $R(t)=R(z(t))$ and $W(t)=W(z(t))$.
If $R(0)\le u$, then
$\|\nabla F(x(0))\|\le u\le s_0$,
  so its position is within $Mu^\gamma$ of $\operatorname{Crit}(F)$.

  Otherwise let $(t_-,t_+)$ be the component of $\{t:R(t)>u\}$
  containing zero. On this interval, $R>u\ge A\lambda$ gives
$\lambda\le R/A$, and the residual controls $\|m\|\le R$. Thus
  \begin{gather*}
-\dot W\ge c_0R^2-B_0\lambda R \ge\left(c_0-\frac{B_0}{A}\right)R^2 \ge\frac{c_0}{2}R^2, \\
\|\dot x\|\le\varepsilon^{-1}\|m\|+\lambda M_x \le(\varepsilon^{-1}+M_x/A)R =:C_xR.
\end{gather*}
  The left endpoint $t_-$ is finite. Otherwise, for every $T>0$,
\[
\left(\max_{z\in L}W(z)-\min_{z\in L}W(z)\right)\ge W(-T)-W(0) \ge\frac{c_0u^2}{2}T,
\]
contradicting boundedness of $W$ on $L$ as $T\to\infty$.
  Continuity gives $R(t_-)=u$, and hence
  $\|\nabla F(x(t_-))\|\le u\le s_0$; the distance at this point
  is at most $Mu^\gamma$. The distance traveled between these times satisfies
  \begingroup
\[
\begin{aligned}
\int_{t_-}^{0}\|\dot x\|\,dt
&\le C_x\int_{t_-}^{0}R\,dt
\le\frac{C_x}{u}\int_{t_-}^{0}R^2\,dt\\
&\le\frac{2C_x}{c_0u}\int_{t_-}^0(-\dot W)\,dt
=\frac{2C_x}{c_0u}\bigl(W(t_-)-W(0)\bigr)\\
&\le\frac{2C_x}{c_0u}\left(\max_{z\in L}W(z)-\min_{z\in L}W(z)\right).
\end{aligned}
\]
\endgroup
  Thus
\[
\operatorname{dist}(x(0),\operatorname{Crit}(F))
\le\operatorname{dist}(x(t_-),\operatorname{Crit}(F))
     +\int_{t_-}^0\|\dot x(t)\|\,dt \le Mu^\gamma+\frac{2C_x}{c_0u}\left(\max_{z\in L}W(z)-\min_{z\in L}W(z)\right).
\]
The initial point in $L$ was arbitrary, so this proves the estimate.
\end{proof}

\clearpage
\section{Construction and simulation of the nonconvergent example}
\label{app:counterexample-details}
\noindent
\begin{minipage}[t]{0.54\linewidth}
\vspace{0pt}
We verify the assumptions for Example~\ref{ex:isolated-corank-two},
give the simulation parameters, and prove nonconvergence.

\subsection{Objective and second moments}\label{guide:appendix-8-1}
\label{app:counterexample-construction}
\label{sec:quartic-cycle-model}
Recall $c=(1,1,1)$ and let $n=c/\sqrt3$. Decompose $x-c$ into
its components parallel and perpendicular to $n$:
\[
x-c=\zeta n+y,\qquad \zeta\in\mathbb R,\quad y\perp n.
\]
Here $\zeta=n^\top(x-c)$ is the signed displacement along $n$.
Figure~\ref{fig:counterexample-geometry} illustrates these coordinates.

The matrix $P=I-nn^\top$ is the orthogonal projection onto $n^\perp$,
so $D_xy=P$ in the derivative calculations below. Define $Ju=n\times u$. On $n^\perp$, $J$ rotates vectors by
$90^\circ$ in the orientation determined by $n$. The objective
in~\eqref{eq:counterexample-objective} and its gradient are
\[
\begin{gathered}
F(x)=\tfrac14\|y\|^4+\tfrac12\zeta^2,\\
\nabla F(x)=\|y\|^2y+\zeta n.
\end{gathered}
\]
The second-moment field in~\eqref{eq:visualization-second-moments} is
\[
\begin{gathered}
S_i(x)=1+(\partial_iF(x))^2+3\bigl(1-3y_i+(Jy)_i\bigr)^2,\\
i=1,2,3.
\end{gathered}
\]
\end{minipage}\hfill
\begin{minipage}[t]{0.43\linewidth}
\vspace{0pt}
\captionsetup{type=figure}
\centering
\includegraphics[width=\linewidth]{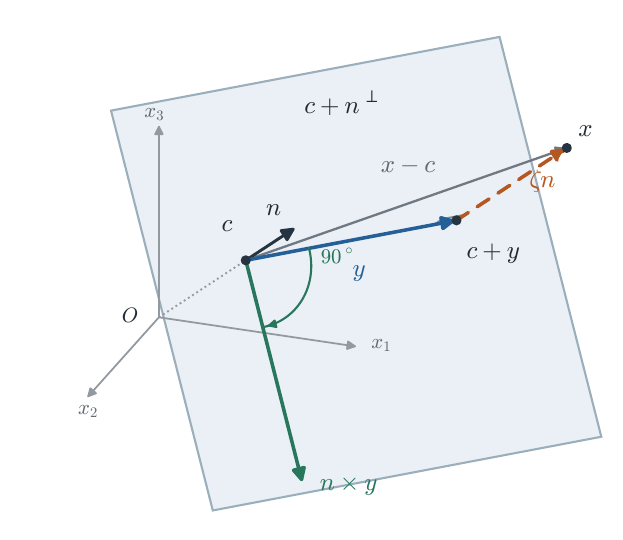}
\caption{The decomposition $x-c=y+\zeta n$.
The point $c+y$ is the orthogonal projection of $x$ onto the shaded
plane $c+n^\perp$. The vector $n\times y$ is obtained by rotating
$y$ through $90^\circ$ in this plane, oriented by $n$.
The point $c$ is the unique minimum of the objective.}
\label{fig:counterexample-geometry}
\end{minipage}
\par\medskip

\paragraph{Derivatives used in the instability proof.}
We now compute $H(x)$ and $\mathrm D S(x)$ for the linearization in
Appendix~\ref{app:counterexample-proof}.

The two components of the decomposition satisfy
\[
D_xy=P,\qquad D_x\zeta=n^\top.
\]
Since $y^\top P=y^\top$, differentiating $\nabla F(x)=\|y\|^2y+\zeta n$ gives
\begingroup
\[
H(x):=\nabla^2F(x)
=\|y\|^2P+2yy^\top P+nn^\top
=\|y\|^2P+2yy^\top+nn^\top.
\]
\endgroup
Since $J$ is linear and independent of $x$, we have
\[
D_x(-3y+Jy)=-3D_xy+J D_xy=-3P+JP.
\]
Using $Jn=n\times n=0$, we get
\[
JP=J(I-nn^\top)=J-(Jn)n^\top=J.
\]
Thus $D_x(-3y+Jy)=-3P+J$.
At $x=hc$ with $h\in\mathbb R$, we have $y=0$ and
$\zeta=\sqrt3(h-1)$, so $\nabla F(hc)=(h-1)c$ and $H(hc)=nn^\top$. Differentiating the two squares
in $S_i$ in a direction $u$ gives
\[
\mathrm D S_i(hc)[u]=2(h-1)(nn^\top u)_i+6((-3P+J)u)_i.
\]
Thus
\begingroup
\begin{equation}\label{eq:counterexample-DS}
\mathrm D S(hc)u=-18u+6Ju \quad (u\perp n),\qquad
\mathrm D S(hc)n=2(h-1)n.
\end{equation}
\endgroup
The rotated displacement $Jy$ in the second moments produces the
rotation term proportional to $J$ in~\eqref{eq:counterexample-expanding-block}.
We choose the initial displacement along
$(1,-1,0)/\sqrt2\in\ker\nabla^2F(c)=n^\perp$, within the plane
where this expansion and rotation occur.

\paragraph{Growth and critical points.}
To verify the assumptions for this example, first note that
orthogonality gives $\|x-c\|^2=\|y\|^2+\zeta^2$. If $F(x)\le M$, then
$\|x-c\|^2\le 2\sqrt M+2M$, so $F$ is coercive. Both terms in $F$
vanish together only at $c$, its unique minimizer. Strict convexity
follows from the strictly convex terms $\|y\|^4$ on $n^\perp$ and
$\zeta^2$ on $\mathbb R$. Since $\nabla F(x)=\|y\|^2y+\zeta n$, its two
orthogonal terms vanish together only at $c$. Thus the critical set
is $\{c\}$ and the critical-value set is $\{0\}$.

\paragraph{Regularity and moment domination.}
The fields $F$ and $S$ are polynomial and hence smooth. The formula
for $S_i$ gives $S_i(x)\ge 1+(\partial_iF(x))^2$, which gives strict positivity
and~\eqref{eq:moment-domination}. With $a=b=1$, we have $0<b\le4a$.
This verifies Assumption~\ref{ass:long-time}.

\paragraph{Critical-manifold condition.}
The only critical component is the singleton $\{c\}$. Evaluating the
Hessian formula at $y=0$ gives $H(c)=nn^\top$, with two-dimensional
kernel $n^\perp$. This verifies Assumption~\ref{ass:critical-manifold}.

\subsection{Simulation setup}\label{guide:appendix-8-2}
\label{app:counterexample-simulation}
We integrate the full nine-dimensional bias-corrected initial-value problem
\eqref{eq:normalized}--\eqref{eq:normalized-initial-value}, with
$F,S$ as defined in Appendix~\ref{app:counterexample-construction}, $a=b=1$, and $\varepsilon=10^{-8}$.
All three runs use
\begingroup
\[
e_1=(1,-1,0)/\sqrt2,\qquad x_0=c+0.04e_1,\qquad
\lambda\in\{10^{-3},3\times10^{-4},10^{-4}\}.
\]
\endgroup
We compute the nearby equilibrium position $p_\lambda$ by solving
\eqref{eq:counterexample-scalar-equilibrium} for $h_\lambda$ on $[0,1]$
using Brent's method and setting $p_\lambda=h_\lambda c$:
\begingroup
\[
p_{10^{-3}}\approx0.998003991014\,c,\qquad
p_{3\times10^{-4}}\approx0.999400359754\,c,\qquad
p_{10^{-4}}\approx0.999800039990\,c.
\]
\endgroup

\subsection{Proof of Example~\ref{ex:isolated-corank-two}}\label{guide:appendix-8-3}
\label{app:counterexample-proof}

We first state the local equilibrium and spectral estimates used in
the nonconvergence argument; their proof follows that argument.

\begin{lemma}[Local equilibrium and spectral splitting]
\label{lem:counterexample-local-splitting}
Let $F,S$ be given by~\eqref{eq:counterexample-objective}
and~\eqref{eq:visualization-second-moments}, set $a=b=1$, and fix
$0<\varepsilon\le1$.
For all sufficiently small $\lambda>0$, there is a
unique equilibrium in a fixed neighborhood of $q_0=(c,0,4c)$, denoted by
\[
q_\lambda=(p_\lambda,\nabla F(p_\lambda),S(p_\lambda)).
\]
There are invertible linear maps
$Q_\lambda:n^\perp\times\mathbb R^7\to\mathbb R^9$, continuous in
$\lambda$ down to zero, which separate the linearized equation into
two blocks:
\[
M_\lambda:=Df_\lambda(q_\lambda)\quad\text{\normalfont and}\quad
Q_\lambda^{-1}M_\lambda Q_\lambda
=\begin{pmatrix}A_\lambda&0\\0&B_\lambda\end{pmatrix}.
\]
In these coordinates, the linearized equations for the unstable and stable
components are, respectively,
\[
\dot u=A_\lambda u,\qquad \dot w=B_\lambda w.
\]
For a trajectory $z(t)$, write
\[
\xi(t):=(u(t),w(t))=Q_\lambda^{-1}(z(t)-q_\lambda),
\qquad u(t)\in n^\perp,\quad w(t)\in\mathbb R^7.
\]
We choose the coordinates so that, for some $\nu>0$ independent of
small $\lambda$,
\begin{equation}\label{eq:counterexample-linear-bounds}
\langle u,A_\lambda u\rangle\ge\frac\lambda4\|u\|^2,
\qquad \langle w,B_\lambda w\rangle\le-\nu\|w\|^2.
\end{equation}
We choose the coordinate map $Q_\lambda$ so that
\begin{equation}\label{eq:counterexample-coordinate-normalization}
Q_0(u,0)=(u,0,\mathrm D S(c)u).
\end{equation}
\end{lemma}

We apply Lemma~\ref{lem:counterexample-local-splitting} to prove the
nonconvergence example.

\begin{restatement}{Example}{ex:isolated-corank-two}
Let $F,S$ be defined in~\eqref{eq:counterexample-objective}
and~\eqref{eq:visualization-second-moments}. Set $a=b=1$ and fix
$0<\varepsilon\le1$. There exist constants $s,\lambda_{\mathrm{ex}}>0$,
which may depend on $\varepsilon$, such that the following holds.
Fix the single initial position
$x_0=c+s(1,-1,0)/\sqrt2$. For every $0<\lambda<\lambda_{\mathrm{ex}}$, the
solution of~\eqref{eq:normalized}--\eqref{eq:normalized-initial-value}
with this $x_0$ is bounded, but $x^\lambda(t)$ does not converge as
$t\to\infty$. The functions $F,S$ and the initial position are
independent of $\lambda$.
\end{restatement}

\begin{proof}
\noindent\textbf{Step 1: Local coordinates.}
Set $a=b=1$, fix $0<\varepsilon\le1$, and let $e_1=(1,-1,0)/\sqrt2$.
The assumptions verified in Appendix~\ref{app:counterexample-construction}
allow us to apply Theorem~\ref{thm:adamw-wp} to obtain the global
corrected solution. Let $q_\lambda$ and $(u,w)$ be the equilibrium and
coordinates from Lemma~\ref{lem:counterexample-local-splitting}.
We will choose one initial position with $\|w(T)\|<\|u(T)\|$ and use
the differential inequalities in Step~6 to rule out convergence to
$q_\lambda$.

\noindent\textbf{Step 2: The trapping set.}
As in Lemma~\ref{lem:finite-time-trap}, we use an energy sublevel set
and compare the perturbed trajectory with the zero-decay flow on
successive unit intervals.

Choose $0<\eta<1$. On $\{F\le\ell\}$, the two terms of $F$ give
$\|x-c\|\le(4\ell)^{1/4}+\sqrt{2\ell}$. Since $S_i(c)=4$, continuity
therefore lets us choose $\ell>0$ so that $|S_i(x)-4|<\eta/2$ on this
sublevel. Define
\[
B:=\{V\le\ell,\ 4-\eta\le v_i\le4+\eta\},
\]
where $V$ is defined in~\eqref{eq:adam-reference-energy}.
On this set,
\[
\|x-c\|\le(4\ell)^{1/4}+\sqrt{2\ell},\quad
\|m\|^2\le2(\varepsilon+\sqrt{4+\eta})\ell,\quad\text{\normalfont and}\quad \|v-4c\|\le\sqrt3\eta.
\]
By taking $\eta$ and then $\ell$ sufficiently small, we can place $B$
in an arbitrarily small neighborhood of $q_0$. Choose a compact
neighborhood $K$ of $B$ inside the uniqueness neighborhood in
Lemma~\ref{lem:counterexample-local-splitting}. We specify how small
$K$ must be after deriving the error estimate below.
Along the zero-decay flow, $V$ is nonincreasing by
\eqref{eq:adam-energy-dissipation}. Thus $F\le V\le\ell$, and
$\dot v_i=S_i(x)-v_i$ points strictly inward at the two moment bounds
because $|S_i(x)-4|<\eta/2$. The energy bound and the inward
directions at the moment bounds keep the flow in $B$. More explicitly,
for a zero-decay trajectory starting in $B$,
\begingroup
\[
|v_i(t)-4|
\le e^{-t}|v_i(0)-4|
  +\int_0^t e^{-(t-r)}|S_i(x(r))-4|\,dr
\le e^{-t}\eta+(1-e^{-t})\eta/2<\eta
\qquad(t>0).
\]
\endgroup
If $V(1)=\ell$,
then $V$ was constant on $[0,1]$. The dissipation inequality forces
$m=0$ throughout this interval, so $\dot m=\nabla F(x)=0$ and $x=c$.
This gives $V=0$, contradicting $\ell>0$. Hence compactness and
continuity give $\Phi_0(1,B)\Subset\operatorname{int}B$.

The set $B$ therefore satisfies the flow conditions of
Lemma~\ref{lem:unit-interval-trap}. Apply it with the compact
neighborhood $K$ to obtain $\kappa>0$.

\noindent\textbf{Step 3: Entry into $B$.}
Since $a=b=1$, subtracting the zero-decay autonomous field from the
corrected field at the same state $z=(x,m,v)$ gives
\[
\dot z-f_0(z)=
\begin{pmatrix}
-\lambda x\\
\dfrac{e^{-t}}{1-e^{-t}}(\nabla F(x)-m)\\
\dfrac{e^{-t}}{1-e^{-t}}(S(x)-v)
\end{pmatrix}.
\]
For $t\ge1$, we have $1-e^{-t}\ge1-e^{-1}$, and the three state-dependent
factors are bounded on $K$. Hence
\[
\|\dot z-f_0(z)\|\le C_K(\lambda+e^{-t})\qquad(t\ge1).
\]
For any such $K$, confinement applies from time $T$ onward if
$z(T)\in B$, $C_Ke^{-T}<\kappa/2$, and $C_K\lambda<\kappa/2$.
Let $z^{s,\lambda}$ solve~\eqref{eq:normalized} with initial value
\[
z^{s,\lambda}(0)=\bigl(c+se_1,\nabla F(c+se_1),S(c+se_1)\bigr).
\]
The trajectory $z^{0,0}$ is constant at $q_0$. For a fixed $T>0$, take a ball of radius
$M$ containing the position paths on $[0,T]$ for $|s|\le1$ and
$\lambda\ge0$, using~\eqref{eq:wp-position}. Write
$E(t)=\sup_{0\le r\le t}\|x^{s,\lambda}(r)-c\|$. Subtracting the position equations and integrating gives
\begingroup
\[
x^{s,\lambda}(t)-c
=se_1-\int_0^t\bigl(R_{x^{s,\lambda}}(r)-R_c(r)\bigr)\,dr
\quad-\lambda\int_0^t x^{s,\lambda}(r)\,dr.
\]
\endgroup
Here $R_c=0$ since $\nabla F(c)=0$. By Lemma~\ref{lem:wp-lipschitz},
$\|R_{x^{s,\lambda}}(r)-R_c(r)\|\le L_T E(r)$ for some $L_T>0$.
Taking the supremum over times in $[0,t]$ and using
$\|x^{s,\lambda}(r)\|\le M$ gives the first bound below.
For the second, bound $|s|+\lambda Mt$ by $|s|+\lambda MT$ on
$[0,T]$ and apply Gronwall's inequality.
\[
E(t)\le |s|+\lambda Mt+L_T\int_0^t E(r)\,dr,
\qquad
E(T)\le (|s|+\lambda MT)e^{L_TT}.
\]
In~\eqref{eq:moment-integrals}, both moments use the weight
\[
k_T(r)=\frac{e^{-(T-r)}}{1-e^{-T}},
\qquad \int_0^T k_T(r)\,dr=1.
\]
Since $\nabla F(c)=0$ and $S(c)=4c$, the Lipschitz bounds on the ball give
\[
\begin{aligned}
\|m^{s,\lambda}(T)\|
&\le\int_0^T k_T(r)\|\nabla F(x^{s,\lambda}(r))-\nabla F(c)\|\,dr
\le L_{\nabla F}E(T),\\
\|v^{s,\lambda}(T)-4c\|
&\le\int_0^T k_T(r)\|S(x^{s,\lambda}(r))-S(c)\|\,dr
\le L_SE(T).
\end{aligned}
\]
Thus $z^{s,\lambda}(T)\to q_0$ as $|s|+\lambda\to0$. Since
$q_0\in\operatorname{int}B$, a ball about $q_0$ is contained in $B$.
For sufficiently small $|s|+\lambda$, the state $z^{s,\lambda}(T)$
lies in that ball.
We use this entry estimate after fixing $K$ and $T$ below.

\noindent\textbf{Step 4: Error bounds and parameter choices.}
Write $f_\lambda^{\mathrm{corr}}(t,z)$ for the right-hand side
of~\eqref{eq:normalized}. For fixed $\lambda$, differentiating
$\xi=Q_\lambda^{-1}(z-q_\lambda)$ and using the block decomposition
in Lemma~\ref{lem:counterexample-local-splitting} gives
\[
\begin{pmatrix}\dot u\\\dot w\end{pmatrix}
=\begin{pmatrix}A_\lambda u\\B_\lambda w\end{pmatrix}
+\begin{pmatrix}E_u(t,\xi)\\E_w(t,\xi)\end{pmatrix},
\]
where, with $z=q_\lambda+Q_\lambda\xi$,
\[
\begin{pmatrix}E_u(t,\xi)\\E_w(t,\xi)\end{pmatrix}
=Q_\lambda^{-1}\bigl[
 f_\lambda(z)-M_\lambda(z-q_\lambda)
 +f_\lambda^{\mathrm{corr}}(t,z)-f_\lambda(z)
\bigr].
\]
Since $f_\lambda(q_\lambda)=0$ and
$Df_\lambda(q_\lambda)=M_\lambda$, Taylor's formula bounds the first
difference by $C\|z-q_\lambda\|^2$ uniformly for small $\lambda$.
The second difference has zero position component. Its moment
components are $\nabla F(x)-m$ and $S(x)-v$, each multiplied by
$e^{-t}/(1-e^{-t})$. Both residuals vanish at $q_\lambda$, so this
difference is bounded by $Ce^{-t}\|z-q_\lambda\|$ for $t\ge1$.
The maps $Q_\lambda$ and $Q_\lambda^{-1}$ are uniformly bounded near
$\lambda=0$. Thus $\|z-q_\lambda\|=\|Q_\lambda\xi\|\le C\|\xi\|$, and
the two remainder bounds give
$\|(E_u,E_w)\|\le C(\|\xi\|^2+e^{-t}\|\xi\|)$. Since $\|\xi\|\le\|u\|+\|w\|$,
\[
\max\{\|E_u(t,\xi)\|,\|E_w(t,\xi)\|\}
\le \delta_t(\|u\|+\|w\|),
\qquad \delta_t:=C(\|\xi\|+e^{-t}).
\]
Fix a bound $M_Q$ for $\|Q_\lambda^{-1}\|$ for all sufficiently
small $\lambda\ge0$. The constant $C$ in the error estimate can be
fixed on a neighborhood of $q_0$ before shrinking $K$. While $z(t)\in K$,
\[
\|\xi(t)\|\le M_Q\bigl(\sup_{z\in K}\|z-q_0\|+\|q_\lambda-q_0\|\bigr).
\]
In the confinement construction, choose $K$ so small that
$CM_Q\sup_{z\in K}\|z-q_0\|\le\nu/32$. After choosing $K$, take $T$ large
enough for the confinement estimate and for $Ce^{-T}\le\nu/16$.
Then restrict $\lambda$ further so that $CM_Q\|q_\lambda-q_0\|\le\nu/32$,
using $q_\lambda\to q_0$. For this fixed $T$, the entry argument gives
$z^{s,\lambda}(T)\in B$ when $|s|+\lambda$ is sufficiently small.
Confinement then gives $z^{s,\lambda}(t)\in K$ for all $t\ge T$, so
\[
\delta_t=C(\|\xi(t)\|+e^{-t})
\le \nu/32+\nu/32+\nu/16=\nu/8.
\]

\noindent\textbf{Step 5: A fixed initial position.}
Our goal is to find a single initial position $x_0=c+se_1$ such that
$\|w^{s,\lambda}(T)\|<\|u^{s,\lambda}(T)\|$ for every sufficiently small
$\lambda>0$. Here $(u^{s,\lambda},w^{s,\lambda})$
are the coordinates of $z^{s,\lambda}$.

First take $\lambda=0$. For the trajectory from $c+se_1$,
Lemma~\ref{lem:hp-zero-decay-energy} gives
\[
\frac14\|y(t)\|^4+\frac12\zeta(t)^2
=F(x^{s,0}(t))\le F(c+se_1)=\frac{s^4}{4}.
\]
Thus $\|y(t)\|\le |s|$ and $|\zeta(t)|\le s^2/\sqrt2$. Substituting
these bounds into $\nabla F(x)=\|y\|^2y+\zeta n$ gives
\[
\|\nabla F(x^{s,0}(t))\|
\le \|y(t)\|^3+|\zeta(t)|
\le |s|^3+\frac{s^2}{\sqrt2}\le 2s^2
\qquad (|s|\le1).
\]
For $0<t\le T$, the first-moment formula in~\eqref{eq:moment-integrals}
and the gradient bound give
\begingroup
\[
\|m^{s,0}(t)\|
\le \frac{a}{A(t)}\int_0^t e^{-a(t-r)}
     \|\nabla F(x^{s,0}(r))\|\,dr
\le \sup_{0\le r\le t}\|\nabla F(x^{s,0}(r))\|
     \frac{a}{A(t)}\int_0^t e^{-a(t-r)}\,dr
\le 2s^2,
\]
\endgroup
since the last kernel integral with its prefactor equals one. Similarly,
$S_i\ge1$ gives
\[
v_i^{s,0}(t)
\ge \frac{b}{B(t)}\int_0^t e^{-b(t-r)}\,dr=1.
\]
Both bounds extend to $t=0$ by continuity. At zero decay, the position
equation and $\varepsilon+\sqrt{v_i^{s,0}}\ge1$ then give
\begingroup
\[
\|x^{s,0}(t)-c-se_1\|
=\left\|\int_0^t \dot x^{s,0}(r)\,dr\right\|
\le \int_0^t \|\dot x^{s,0}(r)\|\,dr
\le \int_0^t \|m^{s,0}(r)\|\,dr
\le 2ts^2\le 2Ts^2.
\]
\endgroup
For fixed $T$ and $\varepsilon$, $O_T(s^2)$ denotes a scalar or vector
remainder bounded in norm by $C_Ts^2$ as $s\to0$, with $C_T$ independent
of $s$. Taylor expansion of $S$ at $c$ gives, uniformly for $0\le r\le T$,
\[
S(x^{s,0}(r))
=S(c)+\mathrm D S(c)(x^{s,0}(r)-c)+O_T(s^2)
=S(c)+s\mathrm D S(c)e_1+O_T(s^2).
\]
Averaging this identity in~\eqref{eq:moment-integrals} at time $T$ gives
$v^{s,0}(T)=S(c)+s\mathrm D S(c)e_1+O_T(s^2)$.
Combining the three estimates at $t=T$ and using
\eqref{eq:counterexample-coordinate-normalization} gives
\begingroup
\[
z^{s,0}(T)-q_0
=s\bigl(e_1,0,\mathrm D S(c)e_1\bigr)+O_T(s^2)
=Q_0(se_1,0)+O_T(s^2).
\]
\endgroup
Applying $Q_0^{-1}$ gives
\[
\bigl(u^{s,0}(T),w^{s,0}(T)\bigr)
=Q_0^{-1}\bigl(z^{s,0}(T)-q_0\bigr)
=(se_1,0)+O_T(s^2).
\]
Taking the two components gives
\begin{equation}\label{eq:counterexample-initial-tangent}
u^{s,0}(T)=se_1+O_T(s^2)\quad\text{\normalfont and}\quad w^{s,0}(T)=O_T(s^2).
\end{equation}
For sufficiently small $s>0$, the estimates give
$\|u^{s,0}(T)\|\ge s/2$ and $\|w^{s,0}(T)\|<s/2$. Fix such an $s$ small
enough that $z^{s,0}(T)\in\operatorname{int}B$. Comparing
$x^{s,\lambda}$ with $x^{s,0}$ as in Step 3 gives
\[
\sup_{0\le t\le T}\|x^{s,\lambda}(t)-x^{s,0}(t)\|
\le \lambda MT e^{L_TT},
\]
since the two initial positions coincide. The moment estimates of Step 3
then give $z^{s,\lambda}(T)\to z^{s,0}(T)$. Since
$q_\lambda\to q_0$ and $Q_\lambda^{-1}\to Q_0^{-1}$, we obtain
$(u^{s,\lambda}(T),w^{s,\lambda}(T))\to(u^{s,0}(T),w^{s,0}(T))$.
Thus $z^{s,\lambda}(T)\in B$ and $\|w^{s,\lambda}(T)\|<\|u^{s,\lambda}(T)\|$
for every sufficiently small $\lambda>0$.

\noindent\textbf{Step 6: The expanding component prevents convergence.}

Fix $s$ and $\lambda$ as in Step~5 and omit the corresponding
superscripts. By that step, $\|w(T)\|<\|u(T)\|$. Suppose equality occurs
later and let
\[
t_*:=\inf\{t>T:\|w(t)\|=\|u(t)\|\}.
\]
By continuity, $t_*>T$, equality holds at $t_*$, and
$\|w(t)\|<\|u(t)\|$ for $T\le t<t_*$. On this interval, $\|u(t)\|>0$,
and the error estimate in Step 4 and~\eqref{eq:counterexample-linear-bounds} give
\begingroup
\[
\frac d{dt}\|u\|
=\frac{\langle u,A_\lambda u+E_u(t,\xi)\rangle}{\|u\|}
\ge\frac\lambda4\|u\|-\delta_t(\|u\|+\|w\|)
\ge(\lambda/4-2\delta_t)\|u\|\ge-\frac\nu4\|u\|.
\]
\endgroup
Integrating from $T$ to $t<t_*$ and letting $t\uparrow t_*$ gives,
by continuity,
\[
\|u(t_*)\|\ge \|u(T)\|e^{-\nu(t_*-T)/4}>0.
\]
Thus $\|w(t_*)\|=\|u(t_*)\|>0$. At $t_*$, we have
\[
\begin{aligned}
\frac12\frac d{dt}(\|w\|^2-\|u\|^2)
&=\langle w,B_\lambda w+E_w\rangle
  -\langle u,A_\lambda u+E_u\rangle\\
&\le-\nu\|w\|^2-\frac\lambda4\|u\|^2
     +\delta_t(\|u\|+\|w\|)^2\\
&=\left(-\nu-\frac\lambda4+4\delta_t\right)\|u\|^2
\le-\frac\nu2\|u\|^2<0.
\end{aligned}
\]
On the other hand, $g(t):=\|w(t)\|^2-\|u(t)\|^2$ satisfies
$g(t_*)=0$ and $g(t_*-h)<0$ for $0<h<t_*-T$. Hence
\[
g'(t_*)=\lim_{h\downarrow0}\frac{g(t_*)-g(t_*-h)}{h}\ge0,
\]
contradicting the strict inequality above. Thus $\|w(t)\|<\|u(t)\|$
for every $t\ge T$.

If $z(t)\to q_\lambda$, then $\xi(t)\to0$ and
$\delta_t\to0$. Choose $T_1\ge T$ so that $\delta_t\le\lambda/16$
for every $t\ge T_1$. The same differential inequality gives
\[
\frac d{dt}\|u(t)\|\ge\left(\frac\lambda4-2\delta_t\right)\|u(t)\| \ge\frac\lambda8\|u(t)\|,\qquad t\ge T_1.
\]
Since $\|u(T_1)\|>0$, integration gives
\[
\|u(t)\|\ge \|u(T_1)\|e^{\lambda(t-T_1)/8} \longrightarrow\infty\qquad(t\to\infty),
\]
contradicting $u(t)\to0$.

Having ruled out $z(t)\to q_\lambda$, we now show that convergence of
$x(t)$ would imply $z(t)\to q_\lambda$.
Suppose $x(t)\to p$. For either input $G=\nabla F$ or $G=S$,
choose $R$ so that $\|G(x(r))-G(p)\|\le\epsilon$ for $r\ge R$.
Writing $k_t(r)=e^{-(t-r)}/(1-e^{-t})$ and
$M_R=\sup_{0\le r\le R}\|G(x(r))-G(p)\|$, we have for $t>R$
\begingroup
\[
\left\|\int_0^t k_t(r)G(x(r))\,dr-G(p)\right\|
\le M_R\int_0^R k_t(r)\,dr+\epsilon\int_R^t k_t(r)\,dr
\le \frac{M_Re^{-(t-R)}}{1-e^{-t}}+\epsilon .
\]
\endgroup
Letting $t\to\infty$ and then $\epsilon\downarrow0$ in
\eqref{eq:moment-integrals} gives $m(t)\to\nabla F(p)$ and
$v(t)\to S(p)$. The position equation then gives
\[
\dot x(t)\longrightarrow d:=-\frac{\nabla F(p)}{\varepsilon+\sqrt{S(p)}}-\lambda p.
\]
Since $x(t)$ converges,
\[
0=\lim_{t\to\infty}(x(t+1)-x(t))
=\lim_{t\to\infty}\int_t^{t+1}\dot x(r)\,dr=d.
\]
The vanishing limit of $\dot x$ means that
$(p,\nabla F(p),S(p))$ is an equilibrium of $f_\lambda$. It belongs
to $K$ by confinement and closedness, so local uniqueness identifies
it with $q_\lambda$, contradicting the preceding argument.
The full trajectory is bounded by Step 4 and continuity on $[0,T]$.

\end{proof}

\paragraph{Verification of the local equilibrium and spectral estimates.}
\begin{proof}[Proof of Lemma~\ref{lem:counterexample-local-splitting}]
\noindent\textbf{Step 1: Uniqueness of the nearby equilibrium.}
At an equilibrium, $m=\nabla F(x)$, $v=S(x)$ and
$\nabla F(x)+\lambda D(x)x=0$, where
$D(x)=\diag(\varepsilon+\sqrt{S_i(x)})$ as
in~\eqref{eq:adamw-equilibrium-set}. Write $x=hc+y$, $y\perp n$.
At $x=hc$, we have $y=0$ and $\zeta=\sqrt3(h-1)$, so
$\nabla F(hc)=\zeta n=(h-1)c$. Substituting into
\eqref{eq:visualization-second-moments} gives
\[
S_i(hc)=1+(h-1)^2+3(1-3\cdot0+0)^2=4+(h-1)^2,
\qquad i=1,2,3.
\]
Write $d(h):=\varepsilon+\sqrt{S_1(hc)}$ for the common denominator.
To show that an equilibrium must have $y=0$, take the inner product of
its position equation with $Jy$. Since $\nabla F(hc+y)=\|y\|^2y+(h-1)c$
and $Jy\perp y,n$, only the decay term contributes. We expand that term
at $hc$ in the direction $y\perp n$.
Taking $u=y$ in~\eqref{eq:counterexample-DS} gives
\[
S(hc+y)=S(hc)+\mathrm D S(hc)y+O(\|y\|^2)
       =S(hc)-18y+6Jy+O(\|y\|^2).
\]
Taking the square root in coordinate $i$, the $i$th diagonal entry of
$D(hc+y)$ is
\[
d(h)+\frac{(-18y+6Jy)_i}{2\sqrt{S_i(hc)}}+O(\|y\|^2),
\]
uniformly for $h$ near one. Since $c_i=1$, multiplication gives
\[
\begin{aligned}
[D(hc+y)(hc+y)]_i
&=\left[d(h)+\frac{-18y_i+6(Jy)_i}{2\sqrt{S_1(hc)}}
+O(\|y\|^2)\right](h+y_i)\\
&=hd(h)+d(h)y_i
+\frac{h[-18y_i+6(Jy)_i]}{2\sqrt{S_1(hc)}}+O(\|y\|^2).
\end{aligned}
\]
Using
\[
\langle Jy,c\rangle=\langle Jy,y\rangle=0,\qquad
\langle Jy,(-18I+6J)y\rangle=6\|Jy\|^2=6\|y\|^2,
\]
we obtain
\begin{gather*}
D(hc+y)(hc+y)
=d(h)hc+\left[d(h)I+\frac{h}{2\sqrt{S_1(hc)}}(-18I+6J)\right]y
 +O(\|y\|^2), \\
\langle Jy,D(hc+y)(hc+y)\rangle
=\frac{3h}{\sqrt{S_1(hc)}}\|y\|^2+O(\|y\|^3).
\end{gather*}
Since $\|Jy\|=\|y\|$ and $S_1(c)=4$, for some uniform $C>0$,
\begingroup
\[
\langle Jy,D(hc+y)(hc+y)\rangle
\ge \left(\frac{3h}{\sqrt{S_1(hc)}}-C\|y\|\right)\|y\|^2,\qquad
\frac{3h}{\sqrt{S_1(hc)}}\longrightarrow\frac32
\qquad(h\to1).
\]
\endgroup
Restrict to $h,y$ with $3h/\sqrt{S_1(hc)}\ge5/4$ and $C\|y\|\le1/4$.
The equilibrium equation then gives
\begingroup
\[
0=\langle Jy,\nabla F(hc+y)+\lambda D(hc+y)(hc+y)\rangle
=\lambda\langle Jy,D(hc+y)(hc+y)\rangle
\ge\lambda(\tfrac54-\tfrac14)\|y\|^2=\lambda\|y\|^2.
\]
\endgroup
Thus $y=0$ for $\lambda>0$, so $x=hc$. Substitution gives
\[
0=\nabla F(hc)+\lambda D(hc)(hc)=(h-1+\lambda h d(h))c.
\]
Since $c\ne0$, this is equivalent to
\begin{equation}\label{eq:counterexample-scalar-equilibrium}
h-1+\lambda h d(h)=0.
\end{equation}
At $(h,\lambda)=(1,0)$, the left-hand side is zero and its
$h$-derivative is one. The implicit function theorem gives a unique
smooth branch $h_\lambda$ near $1$, with $h_0=1$. Hence the unique
nearby equilibrium is
$q_\lambda=(p_\lambda,\nabla F(p_\lambda),S(p_\lambda))$ with
$p_\lambda=h_\lambda c$.

\noindent\textbf{Step 2: The two-dimensional expanding block.}
We construct a two-dimensional invariant plane and prove the growth
bound used in Step 6 of the nonconvergence proof.
Write $h=h_\lambda$, $d=d(h)$, and
$k=-\lambda h/(2\sqrt{S_1(hc)}\,d)$. At $q_\lambda$, the derivative of the
$i$th position equation with respect to $v_i$ is
\[
\frac{m_i}{2\sqrt{v_i}(\varepsilon+\sqrt{v_i})^2}
=\frac{h-1}{2\sqrt{S_i(hc)}\,d^2}
=-\frac{\lambda h}{2\sqrt{S_i(hc)}\,d}=k.
\]
For $j\ne i$, the $i$th position equation does not depend on $v_j$, so
\[
\frac{\partial}{\partial v_j}
\left(-\frac{m_i}{\varepsilon+\sqrt{v_i}}-\lambda x_i\right)=0.
\] Since $y=0$ at $p_\lambda$,
$\nabla^2F(p_\lambda)=nn^\top$. Thus the full linearization is
\[
M_\lambda=
\begin{pmatrix}
-\lambda I&-d^{-1}I&kI\\
nn^\top&-I&0\\
\mathrm D S(p_\lambda)&0&-I
\end{pmatrix}.
\]
Write the linearized displacement as $(\delta x,\delta m,\delta v)$.
The matrices $nn^\top$ and $\mathrm D S(p_\lambda)$ preserve $n^\perp$ and
$\mathbb R n$, so displacements in these two subspaces evolve independently.
When $\delta x,\delta m,\delta v\in n^\perp$, we have
$nn^\top\delta x=0$ and
$\mathrm D S(p_\lambda)\delta x=(-18I+6J)\delta x$ by
\eqref{eq:counterexample-DS}. The equations on $n^\perp$ are therefore
\[
\begin{aligned}
\frac{d}{dt}(\delta x)&=-\lambda\delta x-d^{-1}\delta m+k\delta v,\\
\frac{d}{dt}(\delta m)&=-\delta m,\\
\frac{d}{dt}(\delta v)&=(-18I+6J)\delta x-\delta v.
\end{aligned}
\]
The second equation gives $\delta m(t)=e^{-t}\delta m(0)$, so
$\delta m(0)=0$ implies $\delta m(t)=0$. At $\lambda=0$, with
$\delta x=u\in n^\perp$ and $\delta m=0$, the remaining equations are
\begingroup
\[
\frac{d}{dt}(\delta x)=0,\qquad
\frac{d}{dt}(\delta v)=\mathrm D S(c)u-\delta v.
\]
\endgroup
Thus $(u,0,\mathrm D S(c)u)$ is stationary for the linearized equations at
$\lambda=0$. For small $\lambda$, we seek a plane of the same form
$(u,0,T_\lambda u)$, with $T_\lambda$ to be determined.
If $A_\lambda$ is the matrix acting on $u$ in this plane, invariance
requires
\begingroup
\[
M_\lambda(u,0,T_\lambda u)
=\bigl(({-\lambda I}+kT_\lambda)u,0,
        (\mathrm D S(p_\lambda)-T_\lambda)u\bigr)
=(A_\lambda u,0,T_\lambda A_\lambda u).
\]
\endgroup
Comparing the first and third components and eliminating $T_\lambda$
gives
\begingroup
\[
\begin{aligned}
A_\lambda+\lambda I&=kT_\lambda,\qquad
T_\lambda(I+A_\lambda)=\mathrm D S(p_\lambda)|_{n^\perp},\\
(A_\lambda+\lambda I)(I+A_\lambda)
&=kT_\lambda(I+A_\lambda)=k\,\mathrm D S(p_\lambda)|_{n^\perp}.
\end{aligned}
\]
\endgroup
At $\lambda=0$, $k=0$ and the matrix equation reduces to $A+A^2=0$.
Its derivative at $A=0$ is $H\mapsto H$.
The implicit function theorem gives a smooth solution with $A_0=0$,
so $A_\lambda=O(\lambda)$. Multiplying out the matrix equation gives
$A_\lambda+\lambda I=k\,\mathrm D S(p_\lambda)|_{n^\perp}
-A_\lambda^2-\lambda A_\lambda$.
Using $h=h_\lambda=1+O(\lambda)$ and $S_1(hc)=4+(h-1)^2$, we get
\[
k=-\frac{\lambda h}{2\sqrt{S_1(hc)}\,d(h)}
=-\frac{\lambda}{4(\varepsilon+2)}+O(\lambda^2).
\]
Substitute this and $\mathrm D S(p_\lambda)|_{n^\perp}=-18I+6J$:
\begingroup
\begin{equation}\label{eq:counterexample-expanding-block}
A_\lambda
=-\lambda I+\frac{\lambda}{4(\varepsilon+2)}(18I-6J)+O(\lambda^2)
=\frac{\lambda}{2(\varepsilon+2)}
\bigl((5-2\varepsilon)I-3J\bigr)+O(\lambda^2).
\end{equation}
\endgroup
Finally, $I+A_\lambda$ is invertible for small $\lambda$ because
$A_0=0$. Set $T_\lambda=\mathrm D S(p_\lambda)(I+A_\lambda)^{-1}$ on
$n^\perp$. Then
\[
\begin{aligned}
T_\lambda(I+A_\lambda)&=\mathrm D S(p_\lambda)|_{n^\perp},\\
kT_\lambda
&=k\,\mathrm D S(p_\lambda)|_{n^\perp}(I+A_\lambda)^{-1}\\
&=(A_\lambda+\lambda I)(I+A_\lambda)(I+A_\lambda)^{-1}
=A_\lambda+\lambda I.
\end{aligned}
\]
Thus the plane is invariant.
Since $\langle u,Ju\rangle=0$ and
$(5-2\varepsilon)/(2(\varepsilon+2))\ge1/2$, for small $\lambda>0$
the linearized flow on this plane satisfies
\[
\frac12\frac{d}{dt}\|u\|^2
=\langle u,A_\lambda u\rangle
\ge (\lambda/2-C\lambda^2)\|u\|^2
\ge\frac{\lambda}{4}\|u\|^2.
\]

\noindent\textbf{Step 3: The seven-dimensional contracting block.}
We choose coordinates on the seven-dimensional stable subspace to
obtain the uniform decay estimate used in Step 6.
At $\lambda=0$, the equations on $n^\perp$ from Step 2 have coefficient
matrix, in the order $(\delta m,\delta x,\delta v)$,
\[
\begin{pmatrix}
-I&0&0\\
-(\varepsilon+2)^{-1}I&0&0\\
0&-18I+6J&-I
\end{pmatrix}.
\]
Each block acts on the two-dimensional space $n^\perp$, so the diagonal
blocks give four eigenvalues $-1$ and two zeros.
For displacements along $n$, write
$(\delta x,\delta m,\delta v)=(x_n n,m_n n,v_n n)$,
where $x_n,m_n,v_n$ are real scalars. Since $nn^\top n=n$ and
$\mathrm D S(c)n=0$ by~\eqref{eq:counterexample-DS}, their equations at
$\lambda=0$ are
\[
\begin{aligned}
\frac{d}{dt}x_n&=-\frac{m_n}{\varepsilon+2},\\
\frac{d}{dt}m_n&=x_n-m_n,\\
\frac{d}{dt}v_n&=n^\top\mathrm D S(c)(x_n n)-v_n=-v_n.
\end{aligned}
\]
The $v_n$ equation contributes the eigenvalue $-1$. The $(x_n,m_n)$ block has matrix and characteristic determinant
\begingroup
\[
N=\begin{pmatrix}0&-1/(\varepsilon+2)\\1&-1\end{pmatrix},\qquad
\det(\mu I-N)=\det\begin{pmatrix}\mu&1/(\varepsilon+2)\\-1&\mu+1\end{pmatrix}
=\mu(\mu+1)+1/(\varepsilon+2).
\]
\endgroup
The roots have sum $-1$ and positive product, so both are negative
if real and both have real part $-1/2$ if nonreal. The seven stable eigenvalues remain in the left half-plane for small
$\lambda$ by continuity. Let $E_\lambda^u$ be the plane from Step 2 and $E_\lambda^s$ the
stable subspace. A nonzero intersection would be invariant and carry
an eigenvalue common to both restrictions. Since their spectra are
disjoint,
\[
E_\lambda^u\cap E_\lambda^s=\{0\},\qquad
\dim E_\lambda^u+\dim E_\lambda^s=2+7=9,
\qquad E_\lambda^u\oplus E_\lambda^s=\mathbb R^9.
\]
The stable eigenvalues remain separated from the two near zero.
Apply Kato's theorem~\citep[Chapter~II, \S5.1, Theorem~5.1]{Kato1995}
to each distinct stable eigenvalue of $M_0$ and sum the resulting
projections to obtain $\Pi_\lambda\to\Pi_0$.
For a real basis $(b_j)_{j=1}^7$ of the stable subspace at zero,
\[
b_j(\lambda):=\Pi_\lambda b_j\longrightarrow\Pi_0 b_j=b_j,
\qquad j=1,\ldots,7.
\]
These vectors lie in the real stable subspace and remain linearly
independent for small $\lambda$. Since that subspace has dimension
seven, they form a basis.
With the projected basis, define
\[
Q_\lambda(u,\widetilde w)
=(u,0,T_\lambda u)+\sum_{j=1}^7\widetilde w_j b_j(\lambda).
\]
The two invariant subspaces are complementary, so $Q_\lambda$ is
invertible and
\[
Q_\lambda^{-1}M_\lambda Q_\lambda
=\begin{pmatrix} A_\lambda&0\\0&B_\lambda\end{pmatrix}.
\]
Continuity of $M_\lambda$ and $Q_\lambda$ gives continuity of $B_\lambda$.
At zero, $A_0=0$ and $p_0=c$, so
$Q_0(u,0)=(u,0,\mathrm D S(c)u)$ as in~\eqref{eq:counterexample-coordinate-normalization}.
All eigenvalues of $B_0$ have strictly negative real parts; a matrix
with this property is called Hurwitz. The Lyapunov equation theorem
\citep[Theorem~4.6]{Khalil2002} gives a symmetric positive definite
matrix $P_{\mathrm s}$ such that
\[
B_0^\top P_{\mathrm s}+P_{\mathrm s}B_0=-I.
\]
By continuity of $B_\lambda$, for sufficiently small $\lambda$,
\[
\begin{aligned}
B_\lambda^\top P_{\mathrm s}+P_{\mathrm s}B_\lambda
&=-I+(B_\lambda-B_0)^\top P_{\mathrm s}
+P_{\mathrm s}(B_\lambda-B_0),\\
2\|P_{\mathrm s}\|\,\|B_\lambda-B_0\|&\le\tfrac12.
\end{aligned}
\]
Set $w=P_{\mathrm s}^{1/2}\widetilde w$ and
$\nu=1/(4\|P_{\mathrm s}\|)$. Then
\[
\begin{aligned}
B_\lambda^\top P_{\mathrm s}+P_{\mathrm s}B_\lambda
&\le -I+2\|P_{\mathrm s}\|\|B_\lambda-B_0\|I\le -I/2,\\
\frac12\frac{d}{dt}\|w\|^2
&=\langle w,P_{\mathrm s}^{1/2}B_\lambda P_{\mathrm s}^{-1/2}w\rangle\\
&=\tfrac12\widetilde w^\top
(B_\lambda^\top P_{\mathrm s}+P_{\mathrm s}B_\lambda)\widetilde w\\
&\le -\tfrac14\|\widetilde w\|^2\le-\nu\|w\|^2.
\end{aligned}
\]
The new coordinate map is $Q_\lambda(u,P_{\mathrm s}^{-1/2}w)$.
Keep $Q_\lambda,B_\lambda$ for the rescaled map and stable matrix.
This gives~\eqref{eq:counterexample-linear-bounds}.
\end{proof}

\end{document}